\pdfoutput=1
\RequirePackage{silence}
\documentclass[a4paper,11pt]{amsart}
\usepackage[hmarginratio={1:1},vmarginratio={1:1},lmargin=60.0pt,tmargin=60.0pt]{geometry}

\allowdisplaybreaks

\usepackage[numbers]{natbib}
\usepackage[utf8]{inputenc}

\usepackage{latexsym,exscale,mathtools,textcomp}
\usepackage{amssymb,amsmath,amsthm,amsfonts,mathrsfs,bbm,enumitem,stmaryrd}
\usepackage[table]{xcolor}
\usepackage{graphicx}
\usepackage{mathtools}
\usepackage{ytableau}

\usepackage{xparse}

\usepackage{dynkin-diagrams}

\setlist[enumerate]{itemsep=0.15cm,label=\emph{\upshape(\alph*)}}
\setlist[enumerate,2]{itemsep=0.15cm,label=\emph{\upshape(\roman*)}}

\usepackage{array}
\newcolumntype{C}{>{$}c<{$}}

\definecolor{mygray}{gray}{0.6}
\definecolor{mygraydark}{gray}{0.4}
\definecolor{mygraylight}{gray}{0.85}
\definecolor{spinach}{RGB}{46,139,87}
\definecolor{tomato}{RGB}{255,99,71}
\definecolor{orchid}{RGB}{143,40,194}
\definecolor{neon}{RGB}{77,77,255}
\definecolor{pumpkin}{RGB}{224,180,80}
\definecolor{citron}{RGB}{190,180,90}

\definecolor{lava}{RGB}{207,16,32}
\definecolor{cream}{RGB}{255,253,208}
\definecolor{verdigris}{RGB}{67,179,174}
\definecolor{Black}{RGB}{0,0,0}
\definecolor{mydarkblue}{RGB}{10,10,170}
\definecolor{darkspinach}{RGB}{20,70,20}
\definecolor{darktomato}{RGB}{155,40,30}
\definecolor{darkorchid}{RGB}{50,10,100}
\definecolor{darklava}{RGB}{150,8,16}

\usepackage{todonotes}

\let\emph\relax
\DeclareTextFontCommand{\emph}{\bfseries\em}

\newcommand{\placeholder}{{}_{-}}

\renewcommand{\dots}{\text{...}}
\newcommand{\mystrut}{\rule[-0.2\baselineskip]{0pt}{0.9\baselineskip}}

\let\<=\langle
\let\>=\rangle

\renewcommand{\dots}{\text{...}}
\renewcommand{\vdots}{\rotatebox{90}{\text{...}}}
\renewcommand{\ddots}{\raisebox{0.25cm}{\rotatebox{-45}{\text{...}}}}

\DeclarePairedDelimiterX{\set}[1]{\{}{\}}{\setargs{#1}}
\NewDocumentCommand{\setargs}{>{\SplitArgument{1}{|}}m}{\setargsaux#1}
\NewDocumentCommand{\setargsaux}{mm}
{\IfNoValueTF{#2}{#1} {#1\,\delimsize|\,\mathopen{}#2}}%{#1\:;\:#2}

\newcommand{\ie}{\text{i.e.}}
\newcommand{\eg}{\text{e.g.}}
\newcommand{\cf}{\text{cf.}}
\newcommand{\etc}{\text{etc.}}
\newcommand{\aka}{\text{a.k.a.}}
\newcommand{\ver}{\text{verbatim}}

\newcommand{\muta}{\text{mutatis mutandis}}

\newcommand{\C}{\mathbb{C}}

\newcommand{\N}{\mathbb{Z}_{\geq 0}}

\newcommand{\Z}{\mathbb{Z}}
\newcommand{\K}{\mathbb{K}}
\newcommand{\F}{\mathbb{F}}
\newcommand{\chark}{\mathrm{char}(\mathbb{K})}
\newcommand{\charkk}{\mathrm{char}^{0=\infty}(\mathbb{K})}

\newcommand{\rk}[1][\K]{\mathrm{rk}_{#1}}
\newcommand{\mrk}[1][\K]{\mathrm{rank}_{#1}}
\renewcommand{\dim}[1][\K]{\mathrm{dim}_{#1}}

\newcommand{\Aut}{\mathrm{Aut}}
\newcommand{\End}{\mathrm{End}}
\newcommand{\Hom}{\mathrm{Hom}}

\newcommand{\id}{id}

\newcommand{\algebra}[1][A]{\mathscr{#1}}

\newcommand{\Pcal}{\mathcal{P}}
\newcommand{\apex}{\mathcal{P}^{ap}}
\newcommand{\Tcal}{\mathcal{T}}
\newcommand{\Bcal}{\mathcal{B}}
\newcommand{\cellbasis}[1][\mathscr{A}]{B_{\mathscr{A}}}
\newcommand{\calg}[1][\lambda]{\mathscr{A}_{#1}}
\newcommand{\sandorder}[1][\Pcal]{<_{#1}}
\newcommand{\rsandorder}[1][\Pcal]{>_{#1}}
\newcommand{\gsandorder}[1][\Pcal]{\geq_{#1}}
\newcommand{\lsandorder}[1][\Pcal]{\leq_{#1}}
\newcommand{\sand}[1][\lambda]{\mathscr{H}_{#1}}
\newcommand{\sandbasis}[1][\lambda]{B_{#1}}

\newcommand{\smatrix}[1][{\lambda,K}]{\mathrm{S}^{#1}}
\newcommand{\gmatrix}[1][\lambda]{\mathrm{G}^{#1}}

\newcommand{\dmod}[1][{\lambda,K}]{\Delta(#1)}
\newcommand{\modd}[1][{\lambda,K}]{\nabla(#1)}
\newcommand{\lmod}[1][{\lambda,K}]{L(#1)}

\newcommand{\usummand}{\inplus}

\newcommand{\rad}[1][\jcell]{\mathrm{Rad}}

\newcommand{\lcell}{\mathcal{L}}
\newcommand{\rcell}{\mathcal{R}}
\newcommand{\jcell}{\mathcal{J}}
\newcommand{\hcell}{\mathcal{H}}

\newcommand{\jideal}[1][\lambda]{\algebra^{>_{lr}\lambda}}

\newcommand{\module}[1][M]{{#1}}

\newcommand{\jb}{\mathcal{J}_{b}}
\newcommand{\jt}{\mathcal{J}_{t}}

\newcommand{\Zv}{\Z[v,v^{-1}]}
\newcommand{\hecke}[1][W]{\mathrm{H}_{#1}}
\newcommand{\klbasis}[1][W]{B_{#1}}
\newcommand{\pklbasis}[1][W]{B_{#1}^{p}}

\newcommand{\monoid}[1][S]{\mathrm{#1}}
\newcommand{\monoidunit}{id}
\newcommand{\onemon}{\mathrm{1}}
\newcommand{\xmon}[1][n]{\mathrm{D}_{#1}}
\newcommand{\group}[1][G]{\mathrm{#1}}

\newcommand{\tmon}[1][n]{\mathrm{T}_{#1}}
\newcommand{\ptmon}[1][n]{\mathrm{T}_{#1}^{p}}
\newcommand{\sym}[1][n]{\mathrm{S}_{#1}}
\newcommand{\psym}[1][n]{\mathrm{S}_{#1}^{p}}

\newcommand{\stirling}[2]{\begin{Bsmallmatrix}#1\\#2\end{Bsmallmatrix}}
\newcommand{\bstirling}[2]{\begin{Bmatrix}#1\\#2\end{Bmatrix}}

\newcommand{\pamon}[1][n]{\mathrm{Pa}_{#1}}
\newcommand{\robrmon}[1][n]{\mathrm{RoBr}_{#1}}
\newcommand{\romon}[1][n]{\mathrm{Ro}_{#1}}
\newcommand{\ppamon}[1][n]{\mathrm{Pa}_{#1}^{p}}
\newcommand{\probrmon}[1][n]{\mathrm{RoBr}_{#1}^{p}}
\newcommand{\brmon}[1][n]{\mathrm{Br}_{#1}}
\newcommand{\tlmon}[1][n]{\mathrm{TL}_{#1}}
\newcommand{\momon}[1][n]{\mathrm{Mo}_{#1}}
\newcommand{\promon}[1][n]{\mathrm{Ro}_{#1}^{p}}
\newcommand{\pbrmon}[1][n]{\mathrm{Br}_{#1}^{p}}

\newcommand{\para}{\delta}

\usepackage[all]{xy}
\usepackage{tikz}
\usetikzlibrary{cd}
\usetikzlibrary{decorations}
\usetikzlibrary{decorations.markings}
\usetikzlibrary{decorations.pathreplacing}
\usetikzlibrary{decorations.pathmorphing}
\usetikzlibrary{arrows.meta,shapes,positioning,matrix,calc}
\usetikzlibrary{shapes.callouts}

\tikzset{
anchorbase/.style={baseline={([yshift=#1]current bounding box.center)}},
anchorbase/.default={-0.5ex},
tinynodes/.style={font=\tiny,text height=0.25ex,text depth=0.05ex},
smallnodes/.style={font=\scriptsize,text height=0.75ex,text depth=0.15ex},
mor/.style={line width=0.75,color=black,fill=cream},
mor2/.style={line width=0.75,color=black,fill=tomato},
mor3/.style={line width=0.75,color=black,fill=spinach},
usual/.style={line width=1.2,color=black},
crossline/.style={preaction={draw=white,line width=5.0pt,-},preaction={draw=black,line width=0.9pt,-}},
dot/.style = {
decoration={markings,
post length=0.25mm,
pre length=0.25mm,
mark=at position #1 with {\node[circle,radius=0.15cm,inner sep=-1.2pt,color=black,fill=black]{};}
},
postaction={decorate}
},
dot/.default=1,
}
\tikzstyle directed=[postaction={decorate,decoration={markings,
mark=at position #1 with {\arrow[line width=0.25mm, black]{>}}}}]

\usepackage{aliascnt,etoolbox}
\def\NewTheorem#1{%
\newaliascnt{#1}{equation}%
\newtheorem{#1}[#1]{#1}%
\aliascntresetthe{#1}%
\expandafter\def\csname #1autorefname\endcsname{#1}%
}
\def\equationautorefname~#1\null{(#1)\null}

\numberwithin{equation}{subsection}

\NewTheorem{Proposition}
\NewTheorem{Theorem}
\NewTheorem{Corollary}
\AtEndEnvironment{Corollary}{\null\hfill$\square$}%
\NewTheorem{Lemma}
\theoremstyle{definition}
\NewTheorem{Definition}
\NewTheorem{Notation}
\NewTheorem{Example}
\AtEndEnvironment{Example}{\null\hfill$\Diamond$}%
\theoremstyle{remark}
\NewTheorem{Remark}
\NewTheorem{Assumption}
\NewTheorem{Hypothesis}
\NewTheorem{Conjecture}
\NewTheorem{Question}
\NewTheorem{Observation}
\NewTheorem{Fact}
\NewTheorem{Conclusion}
\NewTheorem{lemma}

\usepackage[hypertexnames=false]{hyperref}
\usepackage{bookmark}
\hypersetup{
pdftoolbar=true,
pdfmenubar=true,
pdffitwindow=false,
pdfstartview={FitH},
pdftitle={Sandwich cellularity and a version of cell theory},
pdfauthor={Daniel Tubbenhauer},
pdfsubject={},
pdfcreator={Daniel Tubbenhauer},
pdfproducer={Daniel Tubbenhauer},
pdfkeywords={},
pdfnewwindow=true,
colorlinks=true,
linkcolor=mydarkblue,
citecolor=teal,
filecolor=magenta,
urlcolor=orchid,
linkbordercolor=lava,
citebordercolor=teal,
urlbordercolor=orchid,
linktocpage=true
}

\newcommand{\nnfootnote}[1]{%
\begin{NoHyper}
\renewcommand\thefootnote{}\footnote{#1}%
\addtocounter{footnote}{-1}%
\end{NoHyper}
}

\def\makeautorefname#1#2{\csdef{#1autorefname}{#2}}

\makeautorefname{section}{Section}%
\makeautorefname{subsection}{Section}%
\makeautorefname{subsubsection}{Section}%

\begin{document}
\title[Sandwich cellularity and a version of cell theory]{Sandwich cellularity and a version of cell theory}
\author[D. Tubbenhauer]{Daniel Tubbenhauer}

\address{D.T.: The University of Sydney, School of Mathematics and Statistics F07, Office Carslaw 827, NSW 2006, Australia, \href{http://www.dtubbenhauer.com}{www.dtubbenhauer.com}, https://orcid.org/0000-0001-7265-5047}
\email{daniel.tubbenhauer@sydney.edu.au}

\begin{abstract}
We explain how the theory of sandwich cellular algebras 
can be seen as a version of cell theory for algebras. 
We apply this theory to many examples such as Hecke algebras, 
and various monoid and diagram algebras.
\end{abstract}

\nnfootnote{\textit{Mathematics Subject Classification 2020.} Primary: 16G10, 16G99; Secondary: 18M30, 20C08, 20M30.}
\nnfootnote{\textit{Keywords.} Cellular algebras, representations of monoids and semigroups, Kazhdan--Lusztig theory and Hecke algebras, diagram algebras.}

\addtocontents{toc}{\protect\setcounter{tocdepth}{1}}

\maketitle

\tableofcontents

%%%%%%%%%%%%%%%%%%%%%%%%%%%%%%%%%%%%%%%%%

\section{Introduction}\label{S:Introduction}

%%%%%%%%%%%%%%%%%%%%%%%%%%%%%%%%%%%%%%%%%

One of the main tools to study the representation theory 
of algebras is the notion of a \emph{cellular algebra} due to Graham and Lehrer \cite{GrLe-cellular}, and cellular algebras are nowadays
ubiquitous in representation theory.
\emph{Sandwich cellular algebras} 
can be seen as a common and strict generalization of 
Graham--Lehrer cellular algebras 
and affine cellular algebras as in \cite{KoXi-affine-cellular}, as 
well as partial generalizations of 
Kazhdan--Lusztig bases, monoid and diagram algebras. All of these, 
in a precise sense, fit under the umbrella of sandwich cellularity.

These sandwich cellular algebras are certain algebras equipped 
with a sandwich cell datum which in turn 
gives rise to the notion of \emph{cells} for these algebras. 
Cells partition sandwich cellular algebras in the same way 
as \emph{Green's relations} \cite{Gr-structure-semigroups} 
partition monoids, or more generally semigroups (we stay 
with monoids in this paper for simplicity), 
and \emph{Kazhdan--Lusztig cells} \cite{KaLu-reps-coxeter-groups} 
partition Hecke algebras. (The latter is the reason why we use 
the name cells.) Sandwich cellular structures 
were also successfully 
applied very early on, albeit in disguise, for example
in the study of Brauer algebras as in 
\cite{Br-gen-matrix-algebras} and \cite{FiGr-canonical-cases-brauer}.

The analogy between sandwich cellular algebras and monoids respectively 
Hecke algebras goes even further. Parts of a sandwich cell datum 
are \emph{sandwich cellular bases}
and \emph{sandwiched algebras}. 
For cellular algebras all sandwiched algebras are trivial, and the 
sandwiched algebras are the main new ingredient in the theory.
These sandwiched algebras play the role of Green's $H$-groups 
\cite{Gr-structure-semigroups}
and asymptotic Hecke algebras associated to intersections of 
left and right cells as {\eg} in \cite{Lu-leading-coeff-hecke}.
The sandwich cellular bases are versions of Kazhdan--Lusztig bases 
in the theory of sandwich cellularity.

The most important theorem regarding
sandwich cellular algebras is \emph{$H$-reduction}, see 
\autoref{T:SandwichCMP}. $H$-reduction classifies simple modules
of sandwich cellular algebras by using their 
cells and the sandwiched algebras.
In the theory of monoids $H$-reduction is known as the celebrated 
\emph{Clifford--Munn--Ponizovski\u{\i} theorem} 
that classifies simple modules of monoids by Green's $J$-classes 
and the simple modules of the $H$-groups, see \cite{GaMaSt-irreps-semigroups} 
or \cite{St-rep-monoid} for modern expositions. In the theory of Kazhdan--Lusztig cells 
the corresponding theorem does not have a name (as far as we know), 
and is weaker than the Clifford--Munn--Ponizovski\u{\i} theorem, see \autoref{S:KL} for details.
A variant of the $H$-reduction for sandwich cellular algebras 
is the theorem with the same name in categorification, see 
{\eg} \cite{MaMaMiZh-h-reduction} and \cite{MaMaMiTuZh-bireps}, 
where we got the nomenclature from. Note however that this categorical 
$H$-reduction is similar in spirit but different in nature, 
which becomes evident when comparing \cite{MaMaMiTuZh-soergel-2reps} with \autoref{S:KL}.  

As we will see, all algebras are sandwich cellular with potentially many different cell structures. Hence, the main point is 
the cell structure itself and not the property of being sandwich cellular. The important task is thus to find a useful cell structure. What useful means is 
difficult to gauge and example dependent. But what one should roughly have in mind here is a fine cell structure with many cells and well-understood sandwiched algebras.

Examples of sandwich cellular algebras with such 
cell structures include many algebras from various diverse fields 
of mathematics:
all cellular algebras with useful 
cell datum, some Hecke algebras and their ($p$\hspace{0.01cm}-)Kazhdan--Lusztig bases, 
many monoid algebras such as transformation monoids, 
many diagram algebras such as Brauer algebras, 
KLR(W) algebras, potentially weighted, seem to have natural sandwich cellular structures, 
see {\eg} 
\cite{Bo-many-cellular-structures}, \cite{MaTu-klrw-algebras} and \cite{MaTu-klrw-algebras-bad}, 
diagram algebras that appear in the study 
of higher genus knot invariants \cite{RoTu-homflypt-handlebody}, \cite{TuVa-handlebody} have useful sandwich cell structures. 
Note that some, but not all, of the named 
algebras are cellular but the sandwich cell structure is 
often much easier to find and to work with, as we 
will see.

Sandwich cellular algebras have been around for a long time, 
however, in various disguises and often implicit in the literature.
They originate from at least four perspectives as mentioned above. In historical order, they appeared via the
Clifford--Munn--Ponizovski\u{\i} theorem in monoid theory,
in the study of the Brauer algebra, 
they are related to
Kazhdan--Lusztig cells, 
and are generalizations of cellular algebras. 
In this work we will draw connections between these 
different fields by taking all of these perspectives at once.

This paper is organized as follows:
\begin{enumerate}[label=$\bullet$]

\item We give a concise summary and 
advance the theory of sandwich cellular algebras 
at the same time, 
see \autoref{S:Sandwich}. In \autoref{SS:SandwichProperties} we give a reformulation 
of the original definition that is useful in practice and 
explains our choice of using cell theory in the title. 
We also make the connection to cellular algebras and monoid theory precise.

\item In \autoref{S:KL} we apply the theory of 
sandwich cellular algebras to Hecke algebras of finite Coxeter 
type and their Kazhdan--Lusztig and $p$-Kazhdan--Lusztig bases.

\item In \autoref{S:DAlgebras} we study various diagram algebras 
from the viewpoint of sandwich cellular algebras. This includes 
diagram algebras without antiinvolution such as transformation 
and planar transformation monoids, as well as diagram algebras with 
an antiinvolution such as Brauer and Temperley--Lieb algebras.

\item In all the examples we study, we classify simple modules 
using $H$-reduction and also compute some dimensions of 
simple modules, which is often doable with the general theory 
of sandwich cellular algebras at hand.

\end{enumerate}

Quite a few, but not all, results in this paper have been obtained 
before, sometimes a long time ago. However, our point is 
that all of them fit under the umbrella of 
sandwich cellularity.
At the end of each section 
we, for convenience, systematically collect references 
and explain how the results in this paper 
compare to known theorems in the literature.

\begin{Remark}\label{R:IntroductionGitHub}
Parts of this paper are based on computer calculations. For the reader that wants to run these calculations themselves, we have 
collected all relevant material on GitHub \cite{Tu-github}.
\end{Remark}

\noindent\textbf{Acknowledgments.}
We thank the organizers of QUACKS (Quantum groups, Categorification, Knot invariants, and Soergel bimodules) 2022 
for the opportunity to present the theory of cells during that event,
which ultimately led to this paper.
Special thanks to James East for pointers to the monoid/semigroup literature and valuable help with GAP, to Robert Spencer for discussions about determinants of Gram matrices, and to 
Alexis Langlois-R{\'e}millard for many helpful comments. We thank the referee for a careful reading and their helpful suggestions for improving the paper.

We are supported, in part, by the Australian Research Council as well as
by the aphorism ``Publish or perish''.

%%%%%%%%%%%%%%%%%%%%%%%%%%%%%%%%%%%%%%%%%

\section{Cell theory for algebras}\label{S:Sandwich}

%%%%%%%%%%%%%%%%%%%%%%%%%%%%%%%%%%%%%%%%%

We now recall the notion of a sandwich cellular 
algebra, and discuss a modification of 
sandwich cellularity which is 
useful in practice 
(in particular, for all algebras 
in this paper).

%%%%%%%%%%%%%%%%%%%%%%%%%%%%%%%%%%%%%%%%%

\subsection{Sandwich cellular algebras}\label{S:SandwichDefinition}

%%%%%%%%%%%%%%%%%%%%%%%%%%%%%%%%%%%%%%%%%

We start with notation:

\begin{Notation}\label{N:SandwichGroundRing}
We use the following conventions.
\begin{enumerate}

\item We fix a commutative unital ring $\K$. This is our ground ring throughout, and {\eg} ranks $\rk$ are with respect to $\K$ if not specified otherwise. We sometimes need $\K$ to be a field, {\eg}
for classification results. Hence the notation $\K$.

\item Throughout this section let $\algebra$ denote an associative unital $\K$-algebra. Algebras in this paper are always assumed to be 
associative and unital.

\item Whenever 
we define an order, say $\sandorder[\Pcal]$, then 
we will also use $\rsandorder[\Pcal]$, $\lsandorder[\Pcal]$
or $\gsandorder[\Pcal]$, having the evident meanings. We also write 
{\eg} $\algebra^{\rsandorder[\Pcal]\lambda}$ which stands for 
the $\K$-submodule of $\algebra$ spanned by
\begin{gather*}
\set{c_{U,n,V}^{\mu}|\mu\in\Pcal,\mu\rsandorder[\Pcal]\lambda,U\in\Tcal(\mu),V\in\Bcal(\mu),n\in\sandbasis[\mu]}
\end{gather*}

\item Unless otherwise specified, modules will always be left modules. In diagrammatic terms these are given by acting from the top, and we use
\begin{gather*}
ab
=
\begin{tikzpicture}[anchorbase,smallnodes,rounded corners]
\node[mor,draw,minimum width=0.5cm,minimum height=0.5cm] at(0,0){\raisebox{-0.05cm}{$b$}};
\node[mor,draw,minimum width=0.5cm,minimum height=0.5cm] at(0,0.5){\raisebox{-0.05cm}{$a$}};
\end{tikzpicture}
.
\end{gather*}
We sometimes use right modules and bimodules, and we will stress whenever that is the case.

\end{enumerate}	
\end{Notation}

\begin{Remark}\label{R:SandwichColors}
As in \autoref{N:SandwichGroundRing}, we will use colors in this 
paper. The colors are a visual aid only and the paper is 
readable in black-and-white without restrictions.
\end{Remark}

The following definition is a modification of \cite[Definition 2A.2]{MaTu-klrw-algebras-bad}, see also \cite[Section 2]{TuVa-handlebody} 
for the same definition using a slightly different formulation.

\begin{Definition}\label{D:SandwichCellularAlgebra}
A \emph{sandwich cell datum} for $\algebra$ is a quadruple 
$\big(\Pcal,(\Tcal,\Bcal),(\sand,\sandbasis),C\big)$, where:
\begin{itemize}

\item $\Pcal=(\Pcal,\sandorder[\Pcal])$ is a poset (the \emph{middle poset} 
with \emph{sandwich order} $\sandorder[\Pcal]$),

\item $\Tcal=\bigcup_{\lambda\in\Pcal}\Tcal(\lambda)$ and $\Bcal=\bigcup_{\lambda\in\Pcal}\Bcal(\lambda)$ are collections of finite sets (the \emph{top/bottom sets}),

\item For $\lambda\in\Pcal$ we have
algebras $\sand[\lambda]$
(the \emph{sandwiched algebras}) and bases $\sandbasis[\lambda]$ of $\sand[\lambda]$,

\item  $C\colon\coprod_{\lambda\in\Pcal}\Tcal(\lambda)\times\sandbasis[\lambda]\times \Bcal(\lambda)\to\algebra;(T,m,B)\mapsto c_{T,m,B}^{\lambda}$ is an injective map,

\end{itemize}
such that:
\begin{enumerate}[label=\upshape(AC${}_{\arabic*}$\upshape)]

\item The set $\cellbasis=\set[\big]{c_{T,m,B}^{\lambda}|\lambda\in\Pcal,T\in\Tcal(\lambda),B\in\Bcal(\lambda),m\in\sandbasis[\lambda]}$
is a basis of $\algebra$.
(We call $\cellbasis$ a \emph{sandwich cellular basis}.)

\item For all $x\in\algebra$ 
there exist scalars $r_{TU}^{x}\in\K$ that do not depend
on $B$ or on $m$, such that
\begin{gather}\label{Eq:SandwichCellCondition}
xc_{T,m,B}^{\lambda}\equiv
\sum_{U\in\Tcal(\lambda),n\in\sandbasis}r_{TU}^{x}c_{U,n,B}^{\lambda}\pmod{\algebra^{\rsandorder[\Pcal]\lambda}}.
\end{gather}
Similarly for right multiplication by $x$.

\item There exists a free 
$\algebra$-$\sand[\lambda]$-bimodule $\dmod[\lambda]$,
a free  $\sand[\lambda]$-$\algebra$-bimodule
$\modd[\lambda]$, and an $\algebra$-bimodule isomorphism
\begin{gather}\label{Eq:SandwichCellAlgebra}
\calg=\algebra^{\gsandorder[\Pcal]\lambda}/\algebra^{\rsandorder[\Pcal]\lambda}\cong\dmod[\lambda]\otimes_{\sand[\lambda]}\modd[\lambda].
\end{gather}
We call $\calg$ the \emph{cell algebra}, and $\dmod[\lambda]$ and 
$\modd[\lambda]$ left and right \emph{cell modules}.

\end{enumerate}
The algebra $\algebra$ is a \emph{sandwich cellular algebra} if it
has a sandwich cell datum.
\end{Definition}

\begin{Definition}\label{D:SandwichInvolutive}
In the setup of \autoref{D:SandwichCellularAlgebra} assume that $\Tcal(\lambda)=\Bcal(\lambda)$ for all $\lambda\in\Pcal$, and that there is an antiinvolution $(\placeholder)^{\star}\colon\algebra\to\algebra$ 
and order two bijections $(\placeholder)^{\star}\colon\sandbasis\to\sandbasis$  such that:
\begin{enumerate}[label=\upshape(AC${}_{4}$\upshape)]

\item We have $(c_{T,m,B}^{\lambda})^{\star}\equiv c_{B,m^{\star},T}^{\lambda}\pmod{\algebra^{\rsandorder[\Pcal]\lambda}}$.

\end{enumerate}
In this case we call the sandwich cell datum 
\emph{involutive} and write $\big(\Pcal,\Tcal,(\sand,\sandbasis),C,(\placeholder)^{\star}\big)$ 
for it.
\end{Definition}

Diagrammatically, although not quite accurate, we think of \autoref{Eq:SandwichCellCondition} and \upshape(AC${}_{4}$\upshape) as
\begin{gather*}
\scalebox{0.85}{$\begin{tikzpicture}[anchorbase,scale=1]
\draw[mor] (0,-0.5) to (0.25,0) to (0.75,0) to (1,-0.5) to (0,-0.5);
\node at (0.5,-0.25){$B^{\prime}$};
\draw[mor] (0,1) to (0.25,0.5) to (0.75,0.5) to (1,1) to (0,1);
\node at (0.5,0.75){$T^{\prime}$};
\draw[mor] (0.25,0) to (0.25,0.5) to (0.75,0.5) to (0.75,0) to (0.25,0);
\node at (0.5,0.25){$m^{\prime}$};
\draw[mor] (0,1) to (0.25,1.5) to (0.75,1.5) to (1,1) to (0,1);
\node at (0.5,1.25){$B$};
\draw[mor] (0,2.5) to (0.25,2) to (0.75,2) to (1,2.5) to (0,2.5);
\node at (0.5,2.25){$T$};
\draw[mor] (0.25,1.5) to (0.25,2) to (0.75,2) to (0.75,1.5) to (0.25,1.5);
\node at (0.5,1.75){$m$};
\end{tikzpicture}
\equiv
\underbrace{\begin{tikzpicture}[anchorbase,scale=1]
\draw[mor] (0,1) to (0.25,0.5) to (0.75,0.5) to (1,1) to (0,1);
\node at (0.5,0.75){$T^{\prime}$};
\draw[mor] (0,1) to (0.25,1.5) to (0.75,1.5) to (1,1) to (0,1);
\node at (0.5,1.25){$B$};
\end{tikzpicture}}_{\in\sand}
\cdot
\begin{tikzpicture}[anchorbase,scale=1]
\draw[mor] (0,-0.5) to (0.25,0) to (0.75,0) to (1,-0.5) to (0,-0.5);
\node at (0.5,-0.25){$B^{\prime}$};
\draw[mor] (0,1) to (0.25,0.5) to (0.75,0.5) to (1,1) to (0,1);
\node at (0.5,0.75){$T$};
\draw[mor] (0.25,0) to (0.25,0.5) to (0.75,0.5) to (0.75,0) to (0.25,0);
\node at (0.5,0.25){\scalebox{0.55}{$mm^{\prime}$}};
\end{tikzpicture}
\pmod{\algebra^{\rsandorder[\Pcal]\lambda}}
,\quad
\left(
\begin{tikzpicture}[anchorbase,scale=1]
\draw[mor] (0,-0.5) to (0.25,0) to (0.75,0) to (1,-0.5) to (0,-0.5);
\node at (0.5,-0.25){$B$};
\draw[mor] (0,1) to (0.25,0.5) to (0.75,0.5) to (1,1) to (0,1);
\node at (0.5,0.75){$T$};
\draw[mor] (0.25,0) to (0.25,0.5) to (0.75,0.5) to (0.75,0) to (0.25,0);
\node at (0.5,0.25){$m$};
\end{tikzpicture}
\right)^{\star}
\equiv
\begin{tikzpicture}[anchorbase,scale=1]
\draw[mor] (0,-0.5) to (0.25,0) to (0.75,0) to (1,-0.5) to (0,-0.5);
\node at (0.5,-0.25){$T$};
\draw[mor] (0,1) to (0.25,0.5) to (0.75,0.5) to (1,1) to (0,1);
\node at (0.5,0.75){$B$};
\draw[mor] (0.25,0) to (0.25,0.5) to (0.75,0.5) to (0.75,0) to (0.25,0);
\node at (0.5,0.25){\scalebox{0.95}{$m^{\star}$}};
\end{tikzpicture}$}
\pmod{\algebra^{\rsandorder[\Pcal]\lambda}}
.
\end{gather*}

\begin{Notation}\label{N:SandwichLanguage}
We will also say, for example, that $\algebra$ itself is sandwich cellular, although this is a bit misleading since an algebra can be sandwich cellular with different sandwich cell data, see {\eg} \autoref{E:SandwichCellularMonoid}.
\end{Notation}

\begin{Example}\label{E:SandwichCellularAlgebra}
Two crucial special cases of involutive sandwich cellular algebras are:
\begin{enumerate}

\item If all sandwiched algebras are isomorphic to $\K$, then 
$\big(\Pcal,\Tcal,(\sand\cong\K,\sandbasis=\{1\}),
C,(\placeholder)^{\star}\big)$ is a cell 
datum and $\algebra$ is a \emph{cellular 
algebra} in the sense of \cite{GrLe-cellular}. (Strictly speaking $\algebra$ 
is cellular in the sense of \cite{GoGr-cellularity-jones-basic} due to 
the weakened condition \upshape(AC${}_{4}$\upshape) on the antiinvolution. 
Throughout, we will always use this weaker version of cellularity.)

\item If all sandwiched algebras are commutative, then 
$\big(\Pcal,\Tcal,(\sand,\sandbasis),C,(\placeholder)^{\star}\big)$ is an affine cell datum and $\algebra$ is an \emph{affine cellular algebra} in the sense of \cite{KoXi-affine-cellular}.

\end{enumerate}
The mnemonic (in particular for readers familiar with diagram algebras) is:
\begin{gather*}
\text{cellular: }
c_{T,1,B}^{\lambda}
\leftrightsquigarrow
\scalebox{0.8}{$\begin{tikzpicture}[anchorbase,scale=1]
\draw[mor] (0,-0.5) to (0.25,0) to (0.75,0) to (1,-0.5) to (0,-0.5);
\node at (0.5,-0.25){$B$};
\draw[mor] (0,0.5) to (0.25,0) to (0.75,0) to (1,0.5) to (0,0.5);
\node at (0.5,0.25){$T$};
\draw[thick,->] (1.5,0)node[right]{$\sand\cong\K$} to (1,0);
\end{tikzpicture}$}
,\quad
\text{affine cellular: }
c_{T,m,B}^{\lambda}
\leftrightsquigarrow
\scalebox{0.8}{$\begin{tikzpicture}[anchorbase,scale=1]
\draw[mor] (0,-0.5) to (0.25,0) to (0.75,0) to (1,-0.5) to (0,-0.5);
\node at (0.5,-0.25){$B$};
\draw[mor] (0,1) to (0.25,0.5) to (0.75,0.5) to (1,1) to (0,1);
\node at (0.5,0.75){$T$};
\draw[mor] (0.25,0) to (0.25,0.5) to (0.75,0.5) to (0.75,0) to (0.25,0);
\node at (0.5,0.25){$m$};
\draw[thick,->] (1.5,0.25)node[right]{commutative $\sand$} to (1,0.25);
\end{tikzpicture}$}
,
\end{gather*}
\begin{gather}\label{Eq:SandwichMnemonic}
\text{sandwich cellular: }
c_{T,m,B}^{\lambda}
\leftrightsquigarrow
\scalebox{0.8}{$\begin{tikzpicture}[anchorbase,scale=1]
\draw[mor] (0,-0.5) to (0.25,0) to (0.75,0) to (1,-0.5) to (0,-0.5);
\node at (0.5,-0.25){$B$};
\draw[mor] (0,1) to (0.25,0.5) to (0.75,0.5) to (1,1) to (0,1);
\node at (0.5,0.75){$T$};
\draw[mor] (0.25,0) to (0.25,0.5) to (0.75,0.5) to (0.75,0) to (0.25,0);
\node at (0.5,0.25){$m$};
\draw[thick,->] (1.5,0.25)node[right]{general $\sand$} to (1,0.25);
\end{tikzpicture}$}
,
\end{gather}
where we assume the existence of an antiinvolution for the 
top two pictures.
\end{Example}

The comparison of sandwich cellular and cellular algebras is as follows.

\begin{Proposition}\label{P:SandwichRefine}
An involutive $\big(\Pcal,\Tcal,(\sand,\sandbasis),C,(\placeholder)^{\star}\big)$
sandwich cell datum is a cell datum if and only if $\sand\cong\K$ 
for all $\lambda\in\Pcal$.

Moreover, an involutive sandwich cellular algebra 
$\algebra$ such that all $\sand$ are cellular (with
the same antiinvolution) is cellular with 
a refined sandwich cell datum. Conversely, if at
least one $\sand$ is non-cellular, then $\algebra$ is non-cellular.
\end{Proposition}

\begin{proof}
The first claim is immediate, for the second
see \cite[Proposition 2.9]{TuVa-handlebody}.
\end{proof}

\begin{Example}\label{E:SandwichCellularMonoid}
For a group $\group$ let $\algebra=\K[\group]$.
Then the group element basis is a sandwich cellular basis.
For $\group$ being the symmetric group or the dihedral group associated to a polygon with an odd number of edges we will 
see a very different sandwich cellular basis in \autoref{S:KL}.

We stress that the group element basis is not a cellular or affine cellular basis in general. For example, if $\group$ is any noncommutative group, then the group element basis is neither cellular nor affine cellular basis but still sandwich cellular.

In \autoref{P:SandwichMonoid} we give a 
condition (that can likely be weakened) on a monoid that 
ensures that the monoid element basis is a sandwich cellular basis.
\end{Example}

We will discuss more examples, and also return to \autoref{E:SandwichCellularMonoid}, later on.

\begin{Proposition}\label{P:SandwichAll}
Any algebra has the structure of a sandwich cellular 
algebra. Moreover, if $\K$ is an algebraically closed field,
then there exists a sandwich cell datum with $\Tcal(\lambda)=\Bcal(\lambda)$ and 
$\sand[\lambda]\cong\K$ for all $\lambda\in\Pcal$.
\end{Proposition}

\begin{proof}
For the first statement take all ingredients to be trivial, {\eg} $\Pcal=\set{0}$, $\Tcal(0)=\Bcal(0)=\set{0}$, $\sand[0]=\algebra$ {\etc} After comparison of definitions, the second claim is \cite[Corollary B]{CoZh-borelic-pairs}.
\end{proof}

\begin{Remark}\label{R:SandwichAll}
By \autoref{P:SandwichAll}, the main point is 
not to prove that an algebra is sandwich cellular, but rather 
to find a useful sandwich cell datum. Note also that 
(the proofs of) \autoref{P:SandwichRefine} and 
\autoref{P:SandwichAll} imply that any algebra is cellular 
over an algebraically closed field
if one does not assume the existence of an antiinvolution.
\end{Remark}

\begin{Notation}\label{N:SandwichAdjoin}
The $\algebra$-bimodule $\calg$ has an induced multiplication, 
but might not have a unit in this multiplication. 
If it does not have a unit, then we always adjoin a unit whenever 
we want to see it as an algebra.
\end{Notation}

An \emph{apex} of an $\algebra$-module $\module$, if it exists, is a maximal $\lambda\in\Pcal$ such that
\begin{gather*}
\K\set[\big]{c_{T,m,B}^{\lambda}|T\in\Tcal(\lambda),B\in\Bcal(\lambda),m\in\sandbasis[\lambda]}
\end{gather*}
does not annihilate $\module$.
The same notion is used for $\calg$ instead of $\algebra$.

\begin{Example}\label{E:SandwichApex}
The cell algebra $\calg$ has at most two apexes, one for the unit and $\lambda$.
\end{Example}

\begin{Lemma}\label{L:SandwichApexExists}
Every simple $\algebra$-module has an apex $\lambda\in\Pcal$, 
and similarly for 
simple $\calg$-modules.
\end{Lemma}

\begin{proof}
See \cite[Lemma 2.15]{TuVa-handlebody}.	
\end{proof}

For a set $X$ of modules we write $X\,/\cong$ for the set 
of isomorphism classes obtained from $X$ by 
identifying isomorphic modules.
The main theorem about sandwich cellular algebras is the following version
of the \emph{Clifford--Munn--Ponizovski\u{\i}} theorem:

\begin{Theorem}(\emph{$H$-reduction})\label{T:SandwichCMP}
Let $\K$ be a field.
\begin{enumerate}

\item If $\lambda\in\Pcal$ is an apex and $\sand$ is Artinian, then we have bijections
\begin{gather*}
\begin{aligned}
\text{Start}\colon&
\\
&\mystrut
\\
\text{$J$-reduction}\colon&
\\
&\mystrut
\\
\text{$H$-reduction}\colon&
\end{aligned}
\begin{gathered}
\left\{
\text{simple $\algebra$-modules with apex $\lambda$}
\right\}/\cong
\\
\xleftrightarrow{1:1}
\\
\left\{
\text{simple $\calg$-modules with apex $\lambda$}
\right\}/\cong
\\
\xleftrightarrow{1:1}
\\
\left\{
\text{simple $\sand$-modules}
\right\}/\cong
.
\end{gathered}
\end{gather*}

\item All $\algebra$-modules with apex $\lambda$ have composition 
factors of apex $\mu$ with $\mu\lsandorder[\Pcal]\lambda$.

\end{enumerate}

\end{Theorem}

\begin{proof}
See \cite[Theorem 2.16]{TuVa-handlebody}.
\end{proof}

\begin{Notation}\label{N:SandwichCMP}	
Whenever \autoref{T:SandwichCMP} applies, we will write 
$\apex\subset\Pcal$ for the set of apexes, and
$\lmod[{\lambda,K}]$ for the simple $\algebra$-modules 
associated to $\lambda\in\apex$ and a simple $\sand$-module $K$.
If $\sand\cong\K$, then we write $\lmod[\lambda]$ for the unique 
simple $\algebra$-module of apex $\lambda$.
\end{Notation} 

\begin{Remark}\label{R:SandwichCMP}
The bijections in \autoref{T:SandwichCMP} can be made explicit 
and $\lmod[{\lambda,K}]$ is the head of 
the induction of $K$ to $\algebra$, see 
\cite[Theorem 2.16]{TuVa-handlebody} for details.
\end{Remark}

\begin{Notation}\label{N:SandwichFinite}
As already indicated by the Artinian condition in \autoref{T:SandwichCMP}, 
there are some technicalities when working with 
infinite dimensional algebras, see 
\cite[Section 2]{TuVa-handlebody} for a more detailed treatment. To simplify the exposition, we will 
from now on assume that our algebras are finite dimensional.
\end{Notation}

\begin{Remark}\label{R:SandwichHistory}
Sandwich cellular algebras are inspired by
the theory of cellular algebras on the one hand, see {\eg}
\cite{GrLe-cellular}, \cite{HuMa-klr-basis},
\cite{KoXi-affine-cellular}, \cite{AnStTu-cellular-tilting}, \cite{EhTu-relcell}, \cite{GuWi-almost-cellular} 
or \cite{TuVa-handlebody}, and ideas coming from monoid 
representation theory on the other hand, see {\eg} 
\cite{Gr-structure-semigroups}, \cite{GaMaSt-irreps-semigroups} or 
\cite{KhSiTu-monoidal-cryptography}.
There was also a huge influence from Kazhdan--Lusztig (KL) theory 
and based algebras \cite{KaLu-reps-coxeter-groups}, \cite{Lu-leading-coeff-hecke}, \cite{KiMa-based-algebras}.
Sandwich cellular algebras have appeared in disguise 
in, for example, \cite{Br-gen-matrix-algebras}, \cite{FiGr-canonical-cases-brauer}, \cite{KoXi-cellular-inflation-morita} 
and \cite{KoXi-brauer}.
We learned the idea underlying 
sandwich cellular algebras from \cite{GuWi-almost-cellular}.
\end{Remark}

%%%%%%%%%%%%%%%%%%%%%%%%%%%%%%%%%%%%%%%%%

\subsection{Green's theory of cells in algebras}\label{SS:SandwichProperties}

%%%%%%%%%%%%%%%%%%%%%%%%%%%%%%%%%%%%%%%%%

We now discuss a special case of sandwich cellular
algebras for which we can 
reformulate \autoref{D:SandwichCellularAlgebra} to be closer to 
Green's classical theory of cells {\aka} Green's relations, 
and KL theory.

\begin{Remark}\label{R:SandwichLoss}
The main slogan of cell theory is the following.
In contrast to groups, the multiplication in an algebra 
often destroys information. For example, if $b=ca$, then $b$ can be obtained 
from $a$ by left multiplication, and we can say that 
$b$ is left bigger than $a$. In a group we can go back by
$c^{-1}b=a$ so $a$ is also left bigger than $b$, 
but this is not always possible in an algebra.
Cells can then be thought of as keeping 
track of the information loss during multiplication. 
The reader is encouraged to keep 
this analogy in mind while reading the text below.
\end{Remark}

Recall that $\algebra$ denotes an algebra as in \autoref{N:SandwichGroundRing}. Fix a basis $\cellbasis$ of 
$\algebra$. Everything below depends on the choice of this basis, and 
we will use the \emph{(based) pair} $(\algebra,\cellbasis)$.

For $a,b,c\in\cellbasis$ 
we write $b\usummand ca$ if, when 
$ca$ is expanded in terms of $\cellbasis$, 
$b$ appear with nonzero coefficient in $ca$.  
We define preorders
on $\cellbasis$ by
\begin{gather*}
(a\leq_{l}b)\Leftrightarrow
\exists c:b\usummand ca
,\quad
(a\leq_{r}b)\Leftrightarrow
\exists d:b\usummand ad
,\quad
(a\leq_{lr}b)\Leftrightarrow
\exists c,d:b\usummand cad
.
\end{gather*}
We call these left, right and two-sided \emph{cell orders}.

\begin{Remark}\label{R:SandwichOrder}
Our convention for the 
cell orders is the same as the most common convention used
in the theory of cellular algebras but the opposite of the one usually used in monoid theory.
\end{Remark}

\begin{Definition}\label{D:SandwichCells}
We define equivalence relations by
\begin{gather*}
(a\sim_{l}b)\Leftrightarrow
(a\leq_{l}b\text{ and }b\leq_{l}a)
,\;
(a\sim_{r}b)\Leftrightarrow
(a\leq_{r}b\text{ and }b\leq_{r}a)
,\;
(a\sim_{lr}b)\Leftrightarrow
(a\leq_{lr}b\text{ and }b\leq_{lr}a)
.
\end{gather*}
The respective equivalence classes are called 
left, right respectively two-sided \emph{cells}.
\end{Definition}

We also say \emph{$J$-cells} instead of two-sided cells, following 
the notation in \cite{Gr-structure-semigroups}.

\begin{Definition}\label{D:SandwichHCells}
An \emph{$H$-cell} $\hcell=\hcell(\lcell,\rcell)=\lcell\cap\rcell$ is an intersection of a left cell $\lcell$ and a right cell $\rcell$.
\end{Definition}

The following is easy:

\begin{Lemma}\label{L:SandwichCellsOrders}
The left, right and $J$-cell orders induce
preorders on the left, right and $J$-cells, respectively. 
Similarly, we can also compare elements and cells via the cell orders.\qed
\end{Lemma}

We will use the same symbols for the various preorders.

\begin{Example}\label{E:SandwichMonoidsAndGroups}
Let $\algebra=\K\monoid$ for a finite monoid 
$\monoid$ (all monoids we use are finite, 
and we drop that adjective), and fix $\cellbasis=\monoid$, the monoid basis.

\begin{enumerate}

\item The above recovers Green's relations. That is, the preorders simplify to
\begin{gather*}
(a\leq_{l}b)\Leftrightarrow
\exists c:b=ca
,\quad
(a\leq_{r}b)\Leftrightarrow
\exists d:b=ad
,\quad
(a\leq_{lr}b)\Leftrightarrow
\exists c,d:b=cad
.
\end{gather*}
In monoid theory the corresponding cells 
are called $L$, $R$, $J$ and $H$-\emph{classes}.

\item Let $\group\subset\monoid$ be the group of invertible elements. 
We have $\group\leq_{l}a$, $\group\leq_{r}a$ 
and $\group\leq_{lr}a$ for all $a\in\monoid$. This can be seen 
by $a=(ag^{-1})g$ and similar calculations for the right and $J$-orders.
In particular, if $\group=\monoid$ (thus, $\monoid$ is a group), then
there is only one cell, which 
is a left, right, $J$ and $H$-cell at the same time.

\end{enumerate}
If not specified otherwise, when referring to monoids 
or groups, then we will always fix $\cellbasis$ to be the monoid basis.
\end{Example}

\begin{Remark}\label{R:SandwichCells}
We think of cells as a matrix decomposition of $\cellbasis$:
\begin{gather}\label{Eq:SandwichCellMatrix}
\begin{tikzpicture}[baseline=(A.center),every node/.style=
{anchor=base,minimum width=1.4cm,minimum height=1cm}]
\matrix (A) [matrix of math nodes,ampersand replacement=\&] 
{
\hcell_{11} \& \hcell_{12} 
\& \hcell_{13} \& \hcell_{14}
\\
\hcell_{21} \& \hcell_{22} 
\& \hcell_{23} \& \hcell_{24}
\\
\hcell_{31} \& \hcell_{32} 
\& \hcell_{33} \& \hcell_{34}
\\
};
\draw[fill=blue,opacity=0.25] (A-3-1.north west) node[blue,left,xshift=0.15cm,yshift=-0.5cm,opacity=1] 
{$\rcell$} to (A-3-4.north east) 
to (A-3-4.south east) to (A-3-1.south west) to (A-3-1.north west);
\draw[fill=red,opacity=0.25] (A-1-3.north west) node[red,above,xshift=0.7cm,opacity=1] 
{$\lcell$} to (A-3-3.south west) 
to (A-3-3.south east) to (A-1-3.north east) to (A-1-3.north west);
\draw[very thick,black,dotted] (A-1-2.north west) to (A-3-2.south west);
\draw[very thick,black,dotted] (A-1-3.north west) to (A-3-3.south west);
\draw[very thick,black,dotted] (A-1-4.north west) to (A-3-4.south west);
\draw[very thick,black,dotted] (A-2-1.north west) to (A-2-4.north east);
\draw[very thick,black,dotted] (A-3-1.north west) to (A-3-4.north east);
\draw[very thick,black] (A-1-1.north west) node[black,above,xshift=-0.5cm] {$\jcell$} to 
(A-1-4.north east) to (A-3-4.south east) 
to (A-3-1.south west) to (A-1-1.north west);
\draw[very thick,black,->] ($(A-1-1.north west)+(-0.4,0.4)$) to (A-1-1.north west);
\draw[very thick,black,->] ($(A-3-4.south east)+(0.5,0.2)$) 
node[right]{$\hcell(\lcell,\rcell)=\hcell_{33}$} 
to[out=180,in=0] ($(A-3-3.south east)+(0,0.2)$);
\end{tikzpicture}
,\quad
\begin{tikzpicture}[anchorbase,scale=1]
\node at (0,0) {\raisebox{0.5cm}{\scalebox{-1}[-1]{\includegraphics[height=3.5cm]{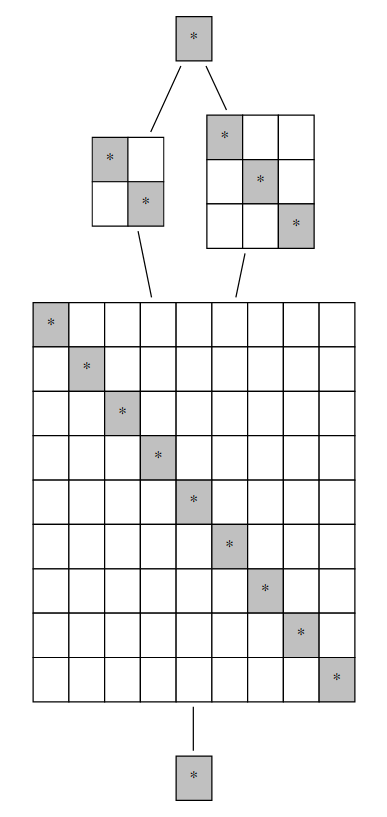}}}};
\end{tikzpicture}
.
\end{gather}
In the left picture we use 
matrix notation for the twelve $H$-cells in $\jcell$ so that $\hcell_{ij}$ is the intersection of the $i$th right with the $j$th left cell of $\jcell$. The $J$-cell $\jcell$ is then a matrix with entries from $H$-cells.
Moreover, the left cells are column vectors and the right cells are row vectors of $\jcell$.

The right picture is a GAP output using 
the package Semigroups, see also \cite{Tu-github}. The illustration shows the cell structure 
of a given, not further specified, monoid in the matrix notion. The $J$-order 
is illustrated as well. The shaded $H$-cells are 
strictly idempotent in the sense 
of \autoref{D:SandwichIdempotent}.
\end{Remark}

We denote the indexing set for the $J$-cells by $\Pcal$, and denote the 
$J$-cells by $\jcell_{\lambda}$ for 
$\lambda\in\Pcal$. 

\begin{Notation}\label{N:SandwichJNotation}
We use a subscript 
$\lambda$ to indicate that we are working in a fixed $J$-cell $\jcell_{\lambda}$.
\end{Notation}

The next lemma is immediate and will be used throughout.

\begin{Lemma}\label{L:SandwichIdeal}
$\jideal=
\K\big\{\bigcup_{\mu\in\Pcal,\mu>_{lr}\lambda}\jcell_{\mu}\big\}$	
is a two-sided ideal in $\algebra$.\qed
\end{Lemma}

\begin{Example}\label{E:SandwichIdeal}
If $\jcell_{\lambda}$ is the $J$-cell 
of size $3$-$3$ on the right-hand side in \autoref{Eq:SandwichCellMatrix}, then 
$\jideal$ is the two-sided ideal supported on the size $9$-$9$ 
and the top $J$-cells. Note that we can ignore the $J$-cell 
of size $2$-$2$ as the multiplication will never get us from 
$\jcell_{\lambda}$ into it.
\end{Example}

\begin{Notation}
An equation {\etc} holds \emph{up to higher order terms} if 
it holds modulo $\jideal$.
\end{Notation}

\begin{Definition}\label{D:SandwichIdempotent}
If the $\K$-linear span of a $J$-cell
$\jcell$ contains a pseudo-idempotent 
up to higher order terms, 
that is, $e\in\K\jcell$ with 
$e^{2}=s(e)\cdot e\pmod{\jideal}$ for $s(e)\in\K\setminus\{0\}$, 
then we call $\jcell$ \emph{idempotent}. 
We say $\jcell$ is \emph{strictly idempotent} if 
$\jcell$ itself contains a pseudo-idempotent up to higher order terms. 
We use the same terminology for $H$-cells.
\end{Definition}

Note that pseudo-idempotents can be rescaled to idempotents 
if $s(e)$ is invertible, which is always true, for example, if $\K$ is a field.
We write $\jcell_{\lambda}(e)$ and 
$\hcell_{\lambda}(e)\subset\jcell_{\lambda}(e)$ for idempotent cells.

\begin{Example}\label{E:SandwichMonoidIdempotent}
All idempotent $H$ and $J$-cells are strictly idempotent 
for monoids. The notion then corresponds 
to having an idempotent and is called 
\emph{regular} in monoid theory.
\end{Example}

A \emph{bottom $J$-cell} is a $J$-cell 
that is minimal in the $J$-order. 
Similarly, a \emph{top $J$-cell} is a $J$-cell 
that is maximal in the $J$-order.

\begin{Proposition}\label{P:SandwichHCells}
For the pair $(\algebra,\cellbasis)$ we have:
\begin{enumerate}
\item Every $H$-cell is contained in some $J$-cell, and every $J$-cell is a disjoint union of $H$-cells.

\item If $s(e)\in\K\setminus\{0\}$ is invertible, then $\sand=\K\hcell_{\lambda}(e)/\jideal$ is an algebra with identity $\frac{1}{s(e)}e$. This algebra is a subalgebra of $\algebra^{\leq_{lr}\lambda}=\algebra/\jideal$ and 
$\algebra^{\lambda}=\algebra^{\leq_{lr}\lambda}\cap\K\jcell_{\lambda}$.

\item Assume that $1\in\cellbasis$. The pair $(\algebra,\cellbasis)$ has a unique bottom $J$-cell, and this $J$-cell is strictly idempotent.

\item The pair $(\algebra,\cellbasis)$ has a unique top $J$-cell, which is an $H$-cell at the same time.

\end{enumerate}
Point (a) works analogously for left and right cells in a fixed $J$ cell.
\end{Proposition}

If they exist, then we write $\jb$ and $\jt$ for the bottom 
and top $J$-cell, respectively.

\begin{proof}
\textit{(a)+(b).} These follow by construction.

\textit{(c).} The bottom $J$-cell is easy to find: 
it is the subset of invertible basis elements, and the 
unit is the idempotent in it.

\textit{(d).} If $\jcell$ and $\jcell^{\prime}$ are maximal $J$-cells, then 
$\jcell=\jcell\jcell^{\prime}=\jcell^{\prime}$ by maximality.
Existence of a maximal $J$ cell follows from the finiteness of 
$\algebra$.
\end{proof}

\begin{Definition}\label{D:SandwichInvolutiveTwo}
The pair $(\algebra,\cellbasis)$ is called \emph{involutive} 
if it admits an order two bijection
${\placeholder}^{\star}\colon\cellbasis\to\cellbasis$ 
that gives rise to an antiinvolution on $\algebra$.
\end{Definition}

\begin{Example}\label{E:SandwichInvolutive}
Groups are involutive with the inversion operation 
being the antiinvolution.
\end{Example}

In the involutive setting all $J$-cells are \emph{square}, meaning 
they have the same number of left and right cells:

\begin{Lemma}\label{L:SandwichInvolutive}
In the involutive setting, ${\placeholder}^{\star}$ gives 
rise to mutually inverse bijections $\{\text{left cells}\}\leftrightarrow\{\text{right cells}\}$
that preserve containment in $J$-cells.
\end{Lemma}

\begin{proof}
Because $(b=ca)^{\star}\Rightarrow(b^{\star}=a^{\star}c^{\star})$.
\end{proof}

The connection to sandwich cellular algebras is given by:

\begin{Theorem}\label{T:SandwichIsSandwich}
For any sandwich cellular algebra $\algebra$ 
we get a pair $(\algebra,\cellbasis)$ for which the above theory applies 
with the following cell structure:
\begin{enumerate}[label=\emph{\upshape(\roman*)}]

\item The basis $\cellbasis$ is the 
sandwich cellular basis.

\item The poset can be taken (potentially changing the order) to 
be $\Pcal=(\Pcal,<_{lr})$.

\item The set $\Tcal(\lambda)$ indexes the right cells within $\jcell_{\lambda}$, and $\Bcal(\lambda)$ indexes the left cells within $\jcell_{\lambda}$.

\item All left, right and 
$H$-cells within one $J$-cell are of the same size.

\item If $\jcell_{\lambda}$ is idempotent, then 
the respective sandwiched algebra and cell algebra are isomorphic to
$\sand$ and $\algebra^{\lambda}$.

\end{enumerate}
\end{Theorem}

\begin{proof}
The first three points are just reformulations 
of \upshape(AC${}_{1}$\upshape) and \upshape(AC${}_{2}$\upshape), 
while the fourth and fifth conditions follow from \upshape(AC${}_{3}$\upshape).
\end{proof}

\begin{Notation}\label{N:SandwichIsSandwich}
We will use \autoref{T:SandwichIsSandwich} to associate a pair $(\algebra,\cellbasis)$ 
to any (involutive) sandwich cellular algebra. 
We call such a pair a(n involutive) \emph{sandwich pair}.
\end{Notation}

If $(\algebra,\cellbasis)$ is a sandwich pair, then 
$H$-reduction \autoref{T:SandwichCMP} applies.

\begin{Remark}\label{R:SandwichIsSandwich}
Part (iv) of \autoref{T:SandwichIsSandwich} is crucial for the $H$-reduction \autoref{T:SandwichCMP} to work. See 
\autoref{E:KLDihedralTwo} for an explicit (counter)example.
\end{Remark}

We now collect a few 
numerical properties of cells. We denote 
by {\eg} $|\jcell_{\lambda}|$ the number of elements in the 
cells and by $\#\lcell_{\lambda}$ {\etc} the number 
of such cells, measured within 
one fixed $J$-cell.

\begin{Lemma}\label{L:SandwichNumericalCells}
For a sandwich pair $(\algebra,\cellbasis)$
we have:
\begin{enumerate}

\item The number of $J$-cells is $|\Pcal|$ and
$\rk(\algebra)=|\cellbasis|=\sum_{\lambda\in\Pcal}|\jcell_{\lambda}|$.

\item $|\jcell_{\lambda}|=\#\lcell_{\lambda}\cdot|\hcell_{\lambda}|
\cdot\#\rcell_{\lambda}$ and $|\jcell_{\lambda}|\cdot|\hcell_{\lambda}|=|\lcell_{\lambda}|\cdot|\rcell_{\lambda}|$.

\item $|\lcell_{\lambda}|=|\hcell_{\lambda}|\cdot\#\rcell_{\lambda}$.

\item $|\rcell_{\lambda}|=|\hcell_{\lambda}|\cdot\#\lcell_{\lambda}$.

\end{enumerate}
\end{Lemma}

\begin{proof}
(a) is clear.
Using the notation from \autoref{Eq:SandwichMnemonic}:
as in \cite[Section 2]{TuVa-handlebody}, that is, as free $\K$-modules we have
\begin{gather}\label{Eq:SandwichFixedDiagrams}
\begin{gathered}
\scalebox{0.8}{$\K\lcell_{\lambda}\cong\K\Tcal(\lambda)\otimes_{\K}\sand
\leftrightsquigarrow
\begin{tikzpicture}[anchorbase,scale=1]
\draw[mor,orchid] (0,-0.5) to (0.25,0) to (0.75,0) to (1,-0.5) to (0,-0.5);
\node at (0.5,-0.25){$B$};
\draw[mor] (0,1) to (0.25,0.5) to (0.75,0.5) to (1,1) to (0,1);
\node at (0.5,0.75){$T$};
\draw[mor] (0.25,0) to (0.25,0.5) to (0.75,0.5) to (0.75,0) to (0.25,0);
\node at (0.5,0.25){$m$};
\end{tikzpicture}
,\quad
\K\rcell_{\lambda}\cong 
\sand\otimes_{\K}\K\Bcal(\lambda)
\leftrightsquigarrow
\begin{tikzpicture}[anchorbase,scale=1]
\draw[mor] (0,-0.5) to (0.25,0) to (0.75,0) to (1,-0.5) to (0,-0.5);
\node at (0.5,-0.25){$B$};
\draw[mor,orchid] (0,1) to (0.25,0.5) to (0.75,0.5) to (1,1) to (0,1);
\node at (0.5,0.75){$T$};
\draw[mor] (0.25,0) to (0.25,0.5) to (0.75,0.5) to (0.75,0) to (0.25,0);
\node at (0.5,0.25){$m$};
\end{tikzpicture}$}
,
\\
\scalebox{0.8}{$\mathcal{J}_{\lambda}\cong \K\Tcal(\lambda)\otimes_{\K}\sand\otimes_{\K}\K\Bcal(\lambda)
\leftrightsquigarrow
\begin{tikzpicture}[anchorbase,scale=1]
\draw[mor] (0,-0.5) to (0.25,0) to (0.75,0) to (1,-0.5) to (0,-0.5);
\node at (0.5,-0.25){$B$};
\draw[mor] (0,1) to (0.25,0.5) to (0.75,0.5) to (1,1) to (0,1);
\node at (0.5,0.75){$T$};
\draw[mor] (0.25,0) to (0.25,0.5) to (0.75,0.5) to (0.75,0) to (0.25,0);
\node at (0.5,0.25){$m$};
\end{tikzpicture}
,\quad
\K\hcell_{\lambda}\cong
\sand
\leftrightsquigarrow
\begin{tikzpicture}[anchorbase,scale=1]
\draw[mor,orchid] (0,-0.5) to (0.25,0) to (0.75,0) to (1,-0.5) to (0,-0.5);
\node at (0.5,-0.25){$B$};
\draw[mor,orchid] (0,1) to (0.25,0.5) to (0.75,0.5) to (1,1) to (0,1);
\node at (0.5,0.75){$T$};
\draw[mor] (0.25,0) to (0.25,0.5) to (0.75,0.5) to (0.75,0) to (0.25,0);
\node at (0.5,0.25){$m$};
\end{tikzpicture}$}
.
\end{gathered}
\end{gather}
We shaded the parts that one can think of as being fixed. 
This follows from \upshape(AC${}_{2}$\upshape)
and \upshape(AC${}_{3}$\upshape) and implies (b), (c) and (d).
\end{proof}

\begin{Proposition}\label{P:SandwichCellsSemisimple}
Let $\K$ be a field and $(\algebra,\cellbasis)$ be a sandwich pair.
The following are equivalent.
\begin{enumerate}[label=\emph{\upshape(\roman*)}]

\item The algebra $\algebra$ is semisimple.

\item All $J$-cells are idempotent and square, all $\sand$ are semisimple and
all $\lmod[{\lambda,K}]$ are isomorphic as $\algebra^{\leq_{lr}\lambda}$-modules to $\mathrm{Ind}_{\sand}^{\algebra^{\leq_{lr}\lambda}}K$.

\end{enumerate}
\end{Proposition}

\begin{proof}
We will use \autoref{L:SandwichNumericalCells} in this proof.

Assume condition (ii) holds. Then we get 
$|\jcell_{\lambda}|=\#\lcell_{\lambda}\cdot|\hcell_{\lambda}|
\cdot\#\rcell_{\lambda}=(\#\lcell_{\lambda})^{2}\cdot|\hcell_{\lambda}|$, because $\jcell_{\lambda}$ is square. Moreover, since 
$\sand$ is semisimple we get $\dim(\sand)=|\hcell_{\lambda}|=\sum_{\text{simples}}\dim(K)^{2}$, and combining this 
with the final assumption gives $|\jcell_{\lambda}|=\sum_{\text{simples}}\dim\big(\lmod[{\lambda,K}]\big)^{2}$.
Since all $J$-cells are idempotent we get 
$\dim(\algebra)=\sum_{\lambda,\text{simples}}\dim\lmod[{\lambda,K}]^{2}$ which shows that $\algebra$ is semisimple.

Assume that (i) holds. Then \autoref{R:SandwichCMP} shows that 
we can read the argument above backwards, showing that (i) implies (ii).
\end{proof}

\begin{Remark}\label{R:SandwichCellTheory}
The discussion above is new but of course strongly inspired 
by very similar construction known throughout the literature, 
see {\eg} \cite{Gr-structure-semigroups}, \cite{KaLu-reps-coxeter-groups},
\cite{KiMa-based-algebras}.
\end{Remark}

%%%%%%%%%%%%%%%%%%%%%%%%%%%%%%%%%%%%%%%%%

\subsection{\texorpdfstring{$J$}{J}-reduction and \texorpdfstring{$H$}{H}-reduction}\label{SS:SandwichJHReduction}

%%%%%%%%%%%%%%%%%%%%%%%%%%%%%%%%%%%%%%%%%

Only sandwich pairs will satisfy $H$-reduction in general, but 
(the much weaker) $J$-reduction works for all pairs:

\begin{Proposition}\label{P:SandwichJReduction}
Let $\K$ be a field and consider the pair $(\algebra,\cellbasis)$.
\begin{enumerate}

\item $J$-reduction as in \autoref{T:SandwichCMP} holds 
for the pair $(\algebra,\cellbasis)$.

\item $H$-reduction as in \autoref{T:SandwichCMP} does not hold
in general for the pair $(\algebra,\cellbasis)$.

\end{enumerate}
\end{Proposition}

\begin{proof}
\textit{(a).} One can copy the arguments in 
\cite[Lemma 2.15]{TuVa-handlebody}.

\textit{(b).} See 
\autoref{E:KLDihedralTwo} for a counterexample.
\end{proof}

\begin{Remark}\label{R:SandwichJRedcution}
The $J$-reduction for monoids is classical, and for 
algebras it is certainly well-known but hard to track down. 
See however \cite[Chapter 4]{Lu-characters-reductive-groups} 
for essentially the same construction in the setting of Hecke 
algebras. In the generality of \autoref{P:SandwichJReduction} $J$-reduction 
appears to be new.
\end{Remark}

%%%%%%%%%%%%%%%%%%%%%%%%%%%%%%%%%%%%%%%%%

\subsection{Sandwich and Gram matrices}\label{SS:SandwichGram}

%%%%%%%%%%%%%%%%%%%%%%%%%%%%%%%%%%%%%%%%%

With reference to \autoref{Eq:SandwichMnemonic} and \autoref{Eq:SandwichFixedDiagrams},
consider the following, purely symbolic, equations:
\begin{gather}\label{Eq:SandwichPairing}
\left\langle
\scalebox{0.8}{$\begin{tikzpicture}[anchorbase,scale=1]
\draw[mor,orchid] (0,-0.5) to (0.25,0) to (0.75,0) to (1,-0.5) to (0,-0.5);
\node at (0.5,-0.25){$B^{\prime}$};
\draw[mor] (0,1) to (0.25,0.5) to (0.75,0.5) to (1,1) to (0,1);
\node at (0.5,0.75){$T$};
\draw[mor] (0.25,0) to (0.25,0.5) to (0.75,0.5) to (0.75,0) to (0.25,0);
\node at (0.5,0.25){$m$};
\end{tikzpicture}$}
,
\scalebox{0.8}{$\begin{tikzpicture}[anchorbase,scale=1]
\draw[mor] (0,-0.5) to (0.25,0) to (0.75,0) to (1,-0.5) to (0,-0.5);
\node at (0.5,-0.25){$B$};
\draw[mor,orchid] (0,1) to (0.25,0.5) to (0.75,0.5) to (1,1) to (0,1);
\node at (0.5,0.75){$T^{\prime}$};
\draw[mor] (0.25,0) to (0.25,0.5) to (0.75,0.5) to (0.75,0) to (0.25,0);
\node at (0.5,0.25){$m^{\prime}$};
\end{tikzpicture}$}
\right\rangle
=
\scalebox{0.8}{$\begin{tikzpicture}[anchorbase,scale=1]
\draw[mor] (0,1) to (0.25,0.5) to (0.75,0.5) to (1,1) to (0,1);
\node at (0.5,0.75){$T$};
\draw[mor,orchid] (0,2.5) to (0.25,2) to (0.75,2) to (1,2.5) to (0,2.5);
\node at (0.5,2.25){$T^{\prime}$};
\draw[mor] (0.25,0) to (0.25,0.5) to (0.75,0.5) to (0.75,0) to (0.25,0);
\node at (0.5,0.25){$m$};
\draw[mor] (0,1) to (0.25,1.5) to (0.75,1.5) to (1,1) to (0,1);
\node at (0.5,1.25){$B$};
\draw[mor,orchid] (0,-0.5) to (0.25,0) to (0.75,0) to (1,-0.5) to (0,-0.5);
\node at (0.5,-0.25){$B^{\prime}$};
\draw[mor] (0.25,1.5) to (0.25,2) to (0.75,2) to (0.75,1.5) to (0.25,1.5);
\node at (0.5,1.75){$m^{\prime}$};
\end{tikzpicture}$}
\equiv
\underbrace{r_{TB}\cdot m^{\prime}m}_{\in\sand}
\pmod{\jideal}
.
\end{gather}
The scalar $r_{TB}\in\K$ is essentially given by \autoref{Eq:SandwichCellCondition}.
This construction using half-diagrams and pairings is
standard in quantum algebra and we will exploit it for
our purposes. We omit the colored boxes from illustrations 
as they are fixed and 
do not play any crucial role.

\begin{Remark}\label{R:SandwichHCellValues}
As already indicated in \autoref{Eq:SandwichPairing}, 
pairings in the sandwich setting take naturally place in 
$\sand$ and not in $\K$. The comparison to \cite{GrLe-cellular} is that 
all pairing therein still take values in $\sand$ 
but we have $\sand\cong\K$ in the cellular setting anyway.
\end{Remark}

As shown in \cite[Lemma 2.12]{TuVa-handlebody}, the natural 
multiplication on the $\algebra$-$\algebra$-bimodule $\calg$ in \autoref{Eq:SandwichCellAlgebra} is determined by a bilinear map 
$\phi^{\lambda}\colon\modd[\lambda]\otimes_{\algebra}\dmod[\lambda]\to\sand$.
One can think of this map as being given by \autoref{Eq:SandwichPairing}.
This map in turn gives rise to a bilinear form, denoted using the same symbol,
$\phi^{\lambda}\colon\dmod[\lambda]\to\Hom_{\K}(\modd[\lambda],\K)$
that we can extend via $\placeholder\otimes_{\sand}K$, where $K$ is any
$\sand$-module:
\begin{gather*}
\phi^{\lambda}\otimes_{\sand}K
=\phi^{\lambda}\otimes_{\sand}\id_{K}
\colon\dmod[\lambda]\otimes_{\sand}K\to
\Hom_{\K}\big(\modd[\lambda],\K\big)\otimes_{\sand}K.
\end{gather*}

\begin{Lemma}\label{L:SandwichMatrix}
Let $(\algebra,\cellbasis)$ be a
sandwich pair. The form
$\phi^{\lambda}$ is determined by a 
$(|\hcell_{\lambda}|\cdot\#\rcell_{\lambda})$-$(\#\lcell_{\lambda}\cdot|\hcell_{\lambda}|)$-matrix with values in $\sand$ given by \autoref{Eq:SandwichPairing}.
Moreover, if $\sand$ is semisimple, then the same is true for $\phi^{\lambda}\otimes_{\sand}K$ with a modification (see in the proof).
\end{Lemma}

\begin{proof}
Take $\#\rcell_{\lambda}$ many representatives 
of the free right $\sand$-module $\dmod[\lambda]$ and 
$\#\lcell_{\lambda}$ many representatives 
of the free left $\sand$-module $\modd[\lambda]$.
Then \autoref{Eq:SandwichPairing} defines the element in $\sand$ 
which determines the pairing. Note hereby that we still have $m$ and $m^{\prime}$
in these pictures which increases the size of the matrix to 
$(|\hcell_{\lambda}|\cdot\#\rcell_{\lambda})$ times $(\#\lcell_{\lambda}\cdot|\hcell_{\lambda}|)$.

With respect to \autoref{Eq:SandwichPairing}, the entries are then as follows:
\begin{gather*}
\scalebox{0.8}{$\begin{tikzpicture}[anchorbase,scale=1]
\draw[mor] (0,1) to (0.25,0.5) to (0.75,0.5) to (1,1) to (0,1);
\node at (0.5,0.75){$T$};
\draw[mor] (0.25,0) to (0.25,0.5) to (0.75,0.5) to (0.75,0) to (0.25,0);
\node at (0.5,0.25){$m$};
\draw[mor] (0,1) to (0.25,1.5) to (0.75,1.5) to (1,1) to (0,1);
\node at (0.5,1.25){$B$};
\draw[mor] (0.25,1.5) to (0.25,2) to (0.75,2) to (0.75,1.5) to (0.25,1.5);
\node at (0.5,1.75){$m^{\prime}$};
\end{tikzpicture}$}
=
\begin{cases}
r_{TB}\cdot m^{\prime}m&\text{if the product is in }\sand\text{ up to higher order terms},
\\
0&\text{else}.
\end{cases}
\end{gather*}

For general $K$, semisimplicity ensures that
every simple $\sand$-module 
$K$ has at least one associated idempotent $e_{K}\in\sand$ and 
the same argument as before works, but using the elements of the form
\begin{gather}\label{Eq:SandwichMatrixPicture}
\scalebox{0.8}{$\begin{tikzpicture}[anchorbase,scale=1]
\draw[mor] (0,1) to (0.25,0.5) to (0.75,0.5) to (1,1) to (0,1);
\node at (0.5,0.75){$T$};
\draw[mor] (0.25,0) to (0.25,0.5) to (0.75,0.5) to (0.75,0) to (0.25,0);
\node at (0.5,0.25){$m$};
\draw[mor] (0,1) to (0.25,1.5) to (0.75,1.5) to (1,1) to (0,1);
\node at (0.5,1.25){$B$};
\draw[mor] (0.25,1.5) to (0.25,2) to (0.75,2) to (0.75,1.5) to (0.25,1.5);
\node at (0.5,1.75){$m^{\prime}$};
\draw[mor,spinach] (0,-0.5) to (0,0) to (1,0) to (1,-0.5) to (0,-0.5);
\node at (0.5,-0.25){$e_{K}$};
\end{tikzpicture}$}
.
\end{gather}
These are obtained from the 
previous ones by putting idempotents 
(denoted above as a box) at the bottom of the diagram used 
in the first part of the proof.
\end{proof}

\begin{Definition}\label{D:SandwichMatrix}
Under the assumptions 
from \autoref{L:SandwichMatrix}, the \emph{sandwich matrix} $\smatrix$ (associated to $\lambda\in\apex$ 
and a simple $\sand$-module $K$) 
is the $(|\hcell_{\lambda}|\cdot\#\rcell_{\lambda})$-$(\#\lcell_{\lambda}\cdot|\hcell_{\lambda}|)$-matrix 
with values in $\sand$ given by \autoref{L:SandwichMatrix}.

For a sandwich pair $(\algebra,\cellbasis)$, the \emph{Gram matrix} $\gmatrix$ 
(associated to $\lambda\in\apex$) is the 
$(\#\rcell_{\lambda})$-$(\#\lcell_{\lambda})$-matrix 
with values in $\K$ defined by
\begin{gather*}
(\gmatrix)_{i,j}=
\begin{cases}
s(e)&\text{if $\hcell_{ij}$ is idempotent with eigenvalue $s(e)$},
\\
0&\text{else}.
\end{cases}
\end{gather*}
\end{Definition}

In other words, $\gmatrix$ is the matrix with entries 
being the eigenvalues of pseudo-idempotents.

\begin{Example}\label{E:SandwichMatrix}
For admissible monoids the sandwich matrix is a 
variation of the matrix with the 
same name in monoid theory, see {\eg} \cite[Section 5.4]{St-rep-monoid}.
Precisely, the matrix given in \cite[Section 5.4]{St-rep-monoid} would be 
the analog of the matrix for $\phi^{\lambda}$.

The Gram matrices for the monoid with cell structure as in \autoref{Eq:SandwichCellMatrix} are identity matrices of sizes 
$1$-$1$ (twice), $3$-$3$, $2$-$2$ and $9$-$9$.
\end{Example}

\begin{Proposition}\label{P:SandwichCellsGram}
Let $\K$ be a field and $(\algebra,\cellbasis)$ be a
sandwich pair.
\begin{enumerate}

\item Under the assumptions 
from \autoref{L:SandwichMatrix},
for $\lambda\in\apex$ and a 
simple $\sand$-modules $K$, we have
\begin{gather*}
\dim\big(\lmod[{\lambda,K}]\big)=
\mrk(\smatrix).
\end{gather*}	

\item For $\sand\cong\K$ and
all $\lambda\in\apex$ we have
\begin{gather*}
\dim\big(\lmod[{\lambda}]\big)=
\mrk(\gmatrix).
\end{gather*}

\end{enumerate}
\end{Proposition}

\begin{proof}
\textit{(a).} This follows from \cite[Lemma 2.13]{TuVa-handlebody} 
and the proof of \autoref{L:SandwichMatrix}.

\textit{(b).} For $\sand\cong\K$, the Gram matrix is 
the sandwich matrix, so (a) applies.
\end{proof}

As we will see in \autoref{T:DAlgebrasTMonTSimples} and in contrast to cellular algebras, the Gram matrix 
is in general not enough to determine the simple dimensions and 
one needs to know the (in general much bigger) sandwich matrices.

\begin{Remark}\label{R:SandwichMatrix}
As we wrote above, the idea of using pairings and matrices 
to determine information about modules is everywhere in quantum 
algebra and related fields. Its a bit hard to track down, 
but see \cite[Section 5.2]{ClPr-semigroups-1} for an early reference.
The above is the version in the 
theory of sandwich cellular algebras.
\end{Remark}

%%%%%%%%%%%%%%%%%%%%%%%%%%%%%%%%%%%%%%%%%

\subsection{Monoids and sandwich cellularity}\label{SS:SandwichMonoids}

%%%%%%%%%%%%%%%%%%%%%%%%%%%%%%%%%%%%%%%%%

We call a monoid $\monoid$ \emph{admissible} if all 
$J$-cells are either idempotent 
with $H$-cells of arbitrary size, or 
have $H$-cells of size one.

\begin{Proposition}\label{P:SandwichMonoid}
For any admissible monoid $\monoid$ (with an antiinvolution) we have a(n involutive) sandwich pair $(\K\monoid,\monoid)$.
\end{Proposition}

\begin{proof}
The crucial fact about monoids 
we need is that left, right and 
$H$-cells in
one $J$-cells are always of the same size, see 
\cite[Theorem 1]{Gr-structure-semigroups}. 
In particular, we can let $\Pcal$ be the poset coming 
from Green cells, we can let $\Tcal$ and $\Bcal$ be indexed 
by right and left cells, and $\sandbasis=\hcell$ to be any $H$-cell 
in $\jcell_{\lambda}$. This choice satisfies \upshape(AC${}_{1}$\upshape) and \upshape(AC${}_{2}$\upshape), by construction.

We need to work a bit more for \upshape(AC${}_{3}$\upshape). 
As free $\K$-modules we can take $\dmod[\lambda]=\K\lcell$, $\modd[\lambda]=\K\rcell$ 
and $\sand=\K\hcell$ for any choice of cells. As a free $\K$-module 
we thus get $\calg\cong\dmod[\lambda]\otimes_{\K}
\sand\otimes_{\K}\modd[\lambda]$. If the $H$-cells in $\jcell$ 
are of size one, then this construction satisfies \upshape(AC${}_{3}$\upshape).

The final case is when 
$\jcell_{\lambda}(e)$ is idempotent, but has arbitrary sized $H$-cells.
In this case $\jcell_{\lambda}(e)$ is strictly idempotent by \autoref{E:SandwichMonoidIdempotent}. It then follows that
$\hcell_{\lambda}(e)$ is a subgroup by
\cite[Theorem 7]{Gr-structure-semigroups}. We can take $\sand=\K\hcell_{\lambda}(e)$ and $\dmod[\lambda]=\K\lcell$, $\modd[\lambda]=\K\rcell$ for the left and right cells defining $\hcell_{\lambda}(e)$ via their intersection. It is then easy to see 
(again with reference to 
\cite{Gr-structure-semigroups}) 
that \upshape(AC${}_{3}$\upshape) holds for these choices
since $\hcell_{\lambda}(e)$ acts freely on its defining left and right cell.

Finally, having an antiinvolution on $\monoid$ clearly 
implies \upshape(AC${}_{4}$\upshape).
\end{proof}

\begin{Example}\label{E:SandwichSchutzerberger}
For admissible monoids the cell modules $\dmod[\lambda]$ are the classical 
Sch{\"u}tzenberger modules with their origin in 
\cite{Sc-cell-modules}.
\end{Example}

\begin{Proposition}\label{P:SandwichMonoidCellular}
For any admissible monoid $\monoid$ with an antiinvolution 
the sandwich pair $(\K\monoid,\monoid)$ can be refined into a cellular 
pair if and only if all $\K\hcell(e)$ are cellular algebras with a compatible antiinvolution.
\end{Proposition}

\begin{proof}
This follows directly from \autoref{P:SandwichRefine} and the proof of 
\autoref{P:SandwichMonoid}.
\end{proof}

\begin{Remark}\label{R:SandwichAdmissible}
The admissibility condition in \autoref{P:SandwichMonoid} 
can be removed if one works with algebras that 
are potentially not unital. In fact, the theory 
of sandwich cellular algebras that are potentially not unital
runs in parallel to the one in this paper; only the technical 
details will change.
\end{Remark}

\begin{Remark}\label{R:SandwichMonoid}
The discussion in this section is new, but the 
question whether monoid algebras are cellular is 
well-studied. See {\eg} \cite{Ea-cellular-semigroups} 
which is very similar to 
\autoref{P:SandwichMonoidCellular}, but our point is
that this follows from the general theory of sandwich cellular 
algebras. Our general theory can also be used to reprove 
the main results in \cite{Wi-cellular-twisted-semigroups}
or \cite{GuXi-cellular-twisted-semigroups}, which in turn 
generalize \cite{Ea-cellular-semigroups}.
\end{Remark}

%%%%%%%%%%%%%%%%%%%%%%%%%%%%%%%%%%%%%%%%%

\section{Kazhdan--Lusztig bases and sandwich cellularity}\label{S:KL}

%%%%%%%%%%%%%%%%%%%%%%%%%%%%%%%%%%%%%%%%%

The KL bases of Hecke algebras partially motivated the 
definition of sandwich cellularity. In this section 
we discuss the relation between these two notions. 

%%%%%%%%%%%%%%%%%%%%%%%%%%%%%%%%%%%%%%%%%

\subsection{The classical story}\label{SS:KLClassical}

%%%%%%%%%%%%%%%%%%%%%%%%%%%%%%%%%%%%%%%%%

Throughout this section, let $W=(W,S)$ be a (connected finite) \emph{Coxeter system}.
We often identify $W$ with its Coxeter diagram.

\begin{Definition}\label{D:KLHecke}
For a parameter $v$ 
let $\hecke=\hecke[{(W,S)}]$ denote the $\Zv$-algebra 
generated by $\{H_{s}|s\in S\}$ subject to the braid relations
and $H_{s}^{2}=1+(v^{-1}-v)H_{s}$.
\end{Definition}

\begin{Notation}\label{N:KLHecke}
Throughout this section $\K$ will denote a specialization of $\Zv$. 
Everything below is defined over $\K$ because the main ingredients 
are defined over $\Zv$.
\end{Notation}

By \cite[Theorem 1.1]{KaLu-reps-coxeter-groups} there exists a distinguished 
basis of $\hecke$ that we call the \emph{KL basis} and that we denote by 
$\klbasis=\klbasis[{(W,S)}]=\{b_{w}\mid w\in W\}$. We will consider the pair $(\hecke,\klbasis)$.

\begin{Remark}\label{R:KLConvetions}
The precise conventions for $\hecke$ and its distinguished basis 
will only be of importance when we 
do explicit calculations in dihedral type. To not distract the 
reader from the main points, we specify our
conventions in \autoref{SS:KLDihedral} after the main 
statements.
\end{Remark}

\begin{Example}\label{E:KLDihedralOne}
For the pair $(\hecke,\klbasis)$ the cells are the KL cells 
from \cite{KaLu-reps-coxeter-groups}. 
Explicitly, take $\K=\C(v)$ and
consider the dihedral type Coxeter system 
$I_{2}(n)$ determined 
by $\dynkin[Coxeter,gonality=n]I{}$ with 
the left node called $1$ and the right node called $2$.
Then we have the cell structure
\begin{gather*}
\noalign{\global\arrayrulewidth=0.5mm}
n=4\text{ (type $B_{2}$)}\colon
\xy
(0,0)*{\begin{gathered}
\begin{tabular}{C}
\cellcolor{mydarkblue!25}
b_{1212}
\end{tabular}
\\[2pt]
\begin{tabular}{C|C}
\arrayrulecolor{tomato}
\cellcolor{mydarkblue!25}b_{1},b_{121} & \cellcolor{spinach!25}b_{12}
\\
\hline
\cellcolor{spinach!25}b_{21} & \cellcolor{mydarkblue!25}b_{2},b_{212}
\\
\end{tabular}
\\[2pt]
\begin{tabular}{C}
\cellcolor{mydarkblue!25}b_{\emptyset}
\end{tabular}
\end{gathered}};
(-22,8.6)*{\jt};
(-22,0)*{\jcell_{m}};
(-22,-8.6)*{\jb};
\endxy
\,,\quad
n=5\text{ (type $H_{2}$)}\colon
\xy
(0,0)*{\begin{gathered}
\begin{tabular}{C}
\cellcolor{mydarkblue!25}
b_{12121}
\end{tabular}
\\[2pt]
\begin{tabular}{C|C}
\arrayrulecolor{tomato}
\cellcolor{mydarkblue!25}b_{1},b_{121} & \cellcolor{spinach!25}b_{12},b_{1212}
\\
\hline
\cellcolor{spinach!25}b_{21},b_{2121} & \cellcolor{mydarkblue!25}b_{2},b_{212}
\\
\end{tabular}
\\[2pt]
\begin{tabular}{C}
\cellcolor{mydarkblue!25}b_{\emptyset}
\end{tabular}
\end{gathered}};
(-24,8.6)*{\jt};
(-24,0)*{\jcell_{m}};
(-24,-8.6)*{\jb};
\endxy
\,.
\end{gather*}
The pattern for general even and $n$ odd is the same: in the even 
case the diagonal $H$-cells have one more element than the off-diagonal 
$H$-cells, while they have the same size for $n$ odd.

As before, the shaded $H$-cells indicate 
strictly idempotent (diagonals) and idempotent but not 
strictly idempotent 
(off-diagonals) $H$-cells, respectively, in the sense 
of \autoref{D:SandwichIdempotent}. This is where 
we use that $\K=\C(v)$ in 
this example, the rest works for any $\K$.
\end{Example}

For the following theorem recall the \emph{bar involution} of $\hecke$, 
see \cite[Introduction]{KaLu-reps-coxeter-groups}.

\begin{Theorem}\label{T:KLSandwich}
For a Coxeter system $W$ we have: 
\begin{enumerate}

\item The pair $(\hecke,\klbasis)$ is a sandwich pair if and only if $W$ is of type $A$ or type $I_{2}(n)$ for $n$ odd.

\item If $W$ is of type $A$ 
or type $I_{2}(n)$ for $n$ odd, then the sandwich pair $(\hecke,\klbasis)$ 
is involutive with the bar involution.

\item For type $I_{2}(n)$ the involutive sandwich pair $(\hecke,\klbasis)$ 
comes neither from a cellular nor an affine 
cellular algebra, but can be refined into a cellular pair.

\end{enumerate}
\end{Theorem}

\begin{proof}
\textit{(a).} If the Coxeter graph is not of type $A$ 
or type $I_{2}(n)$ for $n$ odd, then one will always find $H$-cells of different sizes within one $J$-cell, see {\eg} \cite{Lu-characters-reductive-groups} 
for Weyl types and 
\cite[Section 8]{MaMaMiTuZh-soergel-2reps} 
for the other types. So these types cannot give a sandwich pair $(\hecke,\klbasis)$ by \autoref{T:SandwichIsSandwich}.
Conversely, in type $A$ the KL basis is a cellular basis, as follows from
\cite{KaLu-reps-coxeter-groups}, and 
for type $I_{2}(n)$ with $n$ odd sandwich cellularity can be 
easily verified by hand using the calculations in \autoref{SS:KLDihedral}.

\textit{(b).} This can be proven using
the results in \cite{KaLu-reps-coxeter-groups}.

\textit{(c).} This follows by (a) and 
\autoref{P:SandwichRefine} (and some care with the antiinvolution).
\end{proof}

The following example shows that the assumption of 
$H$-cells being of the same size within one $J$-cell
is crucial for \autoref{T:SandwichCMP} to work.

\begin{Example}\label{E:KLDihedralTwo}
In \autoref{E:KLDihedralOne} consider the non-quantum case 
over $\K=\C$, that is the group algebras of dihedral groups equipped with 
the KL basis $\klbasis[In]$ (we write $In$ for $I_{2}(n)$). 

We consider $n=4$ and $n=5$, the dihedral groups $\group[D]_{4}$ 
and $\group[D]_{5}$ of orders $8$ and $10$. The character tables are 
(the conjugacy classes index the columns, the simple characters the rows and $\phi$ is the golden ratio):
\begin{gather}\label{Eq:KLCharTable}
\group[D]_{4}\colon
\scalebox{0.8}{$\begin{tabular}{C||C|C|C|C|C}
& \chi^{t} & \chi^{b} & \chi^{m,1} & \chi^{m,2} & \chi^{m,3} \\
\hline
\hline
\chi_{t} & 1 & 1 & 1 & 1 & 1 \\
\hline
\chi_{b} & 1 & 1 & -1 & 1 & -1 \\
\hline
\chi_{m,1} & 1 & 1 & 1 & -1 & -1 \\
\hline
\chi_{m,2} & 1 & 1 & -1 & -1 & 1 \\
\hline
\chi_{m,3} & 2 & -2 & 0 & 0 & 0 \\
\end{tabular}$}
\,,\quad
\group[D]_{5}\colon
\scalebox{0.8}{$\begin{tabular}{C||C|C|C|C}
& \chi^{t} & \chi^{b} & \chi^{m1} & \chi^{m2} \\
\hline
\hline
\chi_{t} & 1 & 1 & 1 & 1 \\
\hline
\chi_{b} & 1 & -1 & 1 & 1 \\
\hline
\chi_{m1} & 2 & 0 & \phi & \phi^{2} \\
\hline
\chi_{m2} & 2 & 0 & \phi^{2} & \phi \\
\end{tabular}$}
\,.
\end{gather}
Now, $(\C[\group[D]_{5}],\klbasis[I5])$ 
is a sandwich pair, 
with sandwiched algebras $\sand[b]\cong\C$, $\sand[m]\cong\C[\Z/2\Z]$ and $\sand[t]\cong\C$. Thus, $H$-reduction \autoref{T:SandwichCMP} gives the expected classification 
of simple $\group[D]_{5}$-modules as in \autoref{Eq:KLCharTable} where the subscripts indicate the associated cells. 

In contrast, $(\C[\group[D]_{4}],\klbasis[I4])$ is not a sandwich pair since the algebras one gets are 
$\sand[b]\cong\C$, $\sand[m]\cong\C[\Z/2\Z]$ or $\sand[m]\cong\C$, and $\sand[t]\cong\C$. The $H$-reduction for the 
middle cell fails: neither $\C[\Z/2\Z]$ nor $\C$ give the expected count of three simple $\group[D]_{4}$-modules.
\end{Example}

For $\K$
let $p\in\N$ be minimal such that $p\cdot 1=0\in\K$, and let $p=\infty$ 
if no such $p$ exists.
That is, $p=\chark$ but we consider 
$\chark=0$ as $p=\infty$. We will write $\charkk$ if we 
want to see $\chark=0$ as $p=\infty$.

We use the English convention 
for Young diagrams, {\ie} we draw left-justified rows of boxes 
successively down the page.
Let $\mathrm{P}(\lambda|p)$ 
the set of \emph{$p$-restricted Young diagrams} with $\lambda$ 
boxes. That is, each element of 
$\mathrm{P}(\lambda|p)$ is a Young diagram with 
$\lambda$ boxes such that the difference between the length of two 
consecutive columns is $<p$ (empty columns are of length zero).
Note that the set $\mathrm{P}(\lambda|p)$ for $p>\lambda$ is the 
set of all Young diagrams with $\lambda$ boxes.
The set $\mathrm{P}(\lambda|p)=
\big(\mathrm{P}(\lambda|p),<_{d}\big)$ is a poset 
when equipped with the dominance order $<_{d}$ (with conventions 
specified by \autoref{E:KLpYoung}).

\begin{Example}\label{E:KLpYoung}
For example,
\begin{gather*}
\ytableausetup{centertableaux,boxsize=0.6em}
p>3\colon
\mathrm{P}(3|p)=
\left\{
\ydiagram{3}<_{d}\ydiagram{2,1}<_{d}\ydiagram{1,1,1}
\right\}
,\quad
p\in\{2,3\}\colon
\mathrm{P}(3|p)=
\left\{
\ydiagram{3}<_{d}\ydiagram{2,1}
\right\}
,
\end{gather*}
are the sets of $p$-restricted partitions of $3$ for all
possible $p$.
\end{Example}

For $a\in\N$, let $[a]=\frac{v^{a}-v^{-a}}{v-v^{-1}}\in\Zv$ 
denote the usual quantum numbers. Below we will 
need the condition that $\K=\F(q)$ for $\F$ an algebraically closed 
field and $\K$ is such that quantum numbers $[a]$ evaluated at $q$ are invertible. We call 
such fields \emph{admissible}.

\begin{Example}\label{E:KLFields}
The fields 
$\K=\C(v)$ or $\K=\C$, where $q=v$ and $q=1$ 
respectively, are examples of admissible fields.
\end{Example}

\begin{Theorem}\label{T:KLSimples}
Consider the case where $(\hecke,\klbasis)$ is an involutive sandwich pair. 
\begin{enumerate}

\item Let $W$ be of type $A$.
\begin{enumerate}

\item The $J$-cells of $\hecke[A]$
are given by $\Pcal=\big(\mathrm{P}(\lambda|\infty),<_{d}\big)$.

\item The left cells and right cells of $\hecke[A]$ 
are determine by the Robinson--Schensted correspondence.

\item A $J$-cell of $\hecke[A]$ is strictly idempotent if and only if 
its corresponding partition is $p$-restricted. For strictly idempotent $J$-cells and $\K$ a field we have
\begin{gather*}
\sand\cong\K.
\end{gather*}

\item Let $\K$ be a field. The set of apexes for simple 
$\hecke[A]$-modules is
\scalebox{0.9}{$\apex=\big(\mathrm{P}\big(\lambda|\charkk\big),<_{d}\big)$}, and
there is precisely one simple  
$\hecke[A]$-module $\lmod[\lambda]$ for $\lambda\in\apex$.

\end{enumerate}

\item Let $W$ be odd dihedral type $I_{2}(n)$.
\begin{enumerate}

\item The $J$-cells of $\hecke[In]$
are given by $\Pcal=\{b<_{lr}m<_{lr}t\}$.

\item $\jcell_{b}$ and $\jcell_{t}$ are left and right cells, 
while the left cells in $\jcell_{m}$ are given by reduced expressions 
starting with either $1$ or $2$, and the right cells by reduced expressions 
that end with either $1$ or $2$.

\item $\jcell_{b}$ is always strictly idempotent, 
$\jcell_{m}$ is strictly idempotent if $[2]$ is invertible 
in $\K$, and $\jcell_{m}$ is strictly idempotent if $[2][n]$ is invertible 
in $\K$. Assume $\K$ is an admissible field. Then we have
\begin{gather*}
\sand[b]\cong\K
,\quad
\sand[m]\cong\K[\Z/(\tfrac{n-1}{2})\Z]
,\quad
\sand[t]\cong\K
.
\end{gather*}

\item Let $\K$ be an admissible field.
The set of apexes for simple 
$\hecke[In]$-modules is
$\apex=\Pcal$, and
there is precisely one simple  
$\hecke[In]$-module for $\lambda\in\{b,t\}\subset\apex$, 
and there are precisely $\tfrac{n-1}{2}$ simple  
$\hecke[In]$-modules for $m\in\apex$.
The dimensions of the simple $\hecke[In]$-modules 
are $1$, $1$ and $2$, the latter 
$\tfrac{n-1}{2}$ times.

\end{enumerate}

\item Let $\K$ be an admissible field.
We have that $\hecke[A]$ and $\hecke[In]$ for $n$ odd are semisimple.

\end{enumerate}

\end{Theorem}

\begin{proof}
\textit{(a).} Proving how the cell 
structure looks like is not trivial, but still well-known 
and follows from \cite{KaLu-reps-coxeter-groups}. Thus, it remains 
to verify part (iv).
That $\apex=\mathrm{P}\big(\lambda|\charkk\big)$ 
follows from the calculation of eigenvalues 
of distinguished involutions in the Hecke algebra, 
see \cite[Conjectures 14.2 and Chapter 15]{Lu-hecke-book} 
(we mean the arXiv version of the book). The rest is then \autoref{T:KLSandwich}
and $H$-reduction \autoref{T:SandwichCMP}.

\textit{(b).} This follows from \autoref{SS:KLDihedral} below. 
That is, the hardest parts are proven in \autoref{L:KLZModN} 
and \autoref{L:KLSandwichMatrix}. The rest 
is easy to prove using the explicit description of the multiplication 
of KL basis elements given in \autoref{SS:KLDihedral}.

\textit{(c).} By the above and \autoref{P:SandwichCellsSemisimple}, 
it remains to show that the $\dmod[\lambda]$ are simple. However, 
we already know the dimensions of the simples by (a).(iv) and (b).(iv), 
and the proof completes.
\end{proof}

\begin{Example}\label{E:SandwichDihedralTwo}
Back to \autoref{E:KLDihedralTwo} for $n=5$. 
The sandwich matrices for the bottom and top $J$-cell 
are $\smatrix[b]=\begin{pmatrix}b_{\emptyset}\end{pmatrix}$ 
and $\smatrix[t]=\begin{pmatrix}10\cdot b_{12121}\end{pmatrix}$, respectively.
We compute the sandwich matrices $\smatrix[{m,1}]$ 
and $\smatrix[{m,2}]$ in the proof of 
\autoref{L:KLSandwichMatrix}. These matrices are of rank $2$.
Thus, the bottom and top $J$-cells give one dimensional 
simple $\algebra$-modules, while the middle $J$-cell gives two 
simple $\algebra$-modules of dimension $2$.
\end{Example}

\begin{Remark}\label{R:KLSandwich}
\autoref{T:KLSandwich} also shows that the KL basis is cellular 
if and only if $W$ is of type $A$. That the KL basis 
of type $A$ is cellular 
follows from \cite{KaLu-reps-coxeter-groups}, as we wrote above. 
The converse is known, 
but hard to find in the literature.
The case of general sandwich cellularity is new.

We point out that all finite type Hecke algebras are sandwich cellular, 
even cellular by \cite{Ge-hecke-cellular}, and the statement 
in \autoref{T:KLSandwich} is that the 
KL basis is not a sandwich cellular basis.

Finally, our sandwich approach to classify simple $\hecke[In]$-modules 
from \autoref{T:KLSimples} is new, 
but the results are certainly not new: In the semisimple 
case the classification follows directly from Tits' deformation theorem 
(explicitly, \cite[Theorem 8.1.7]{GePf-characters-coxeter-groups}) and classical 
theory. We also stress that the approach taken in \autoref{T:KLSimples} 
with a bit more work, see \autoref{R:KLZModN} below, 
also classifies simple $\hecke[In]$-modules 
over arbitrary fields. This classification is again 
well-known, see {\eg} \cite[Chapter 8]{GePf-characters-coxeter-groups} 
for a concise treatment and bibliographical remarks, but our point 
is that it also follows as part of the general theory of sandwich cellular 
algebras.

Note that, by \autoref{P:SandwichJReduction}, $J$-reduction \autoref{T:SandwichCMP}
works for all pairs $(\hecke,\klbasis)$
in arbitrary type (not necessary connected or finite). This recovers an result of {\eg} \cite[Section 4]{Lu-characters-reductive-groups}.
\end{Remark}

%%%%%%%%%%%%%%%%%%%%%%%%%%%%%%%%%%%%%%%%%

\subsection{\texorpdfstring{$p$}{p}-Kazhdan--Lusztig bases and sandwich cellularity}\label{SS:KLp}

%%%%%%%%%%%%%%%%%%%%%%%%%%%%%%%%%%%%%%%%%

For a fixed prime $p$ 
{\eg} \cite{JeWi-p-canonical} defines another distinguished basis 
$\pklbasis=\{b_{w}^{p}\mid w\in W\}$ of 
$\hecke$ that we call the \emph{$p$-KL basis}.
Not much is known about the (sandwich) cellularity of the $p$-KL basis, and all we can say is:

\begin{Theorem}\label{T:KLSandwichp}
\leavevmode
\begin{enumerate}

\item If $W$ is of type $A$, then $(\hecke,\pklbasis)$ is a sandwich pair 
for all primes, and the cell structure is the same as for $(\hecke,\klbasis)$.

\item If $W$ is of type $B_{2}$, $C_{2}$, $G_{2}$, $B_{3}$, $C_{3}$, $B_{4}$, $C_{4}$ or $D_{4}$, then $(\hecke,\pklbasis)$ is never a sandwich pair.

\item Let $W$ be a Weyl group that is not of type $A$.
Varying over all primes $p$, the pair $(\hecke,\klbasis)$ is 
not a sandwich pair up to finitely many potential exceptions.

\end{enumerate}
\end{Theorem}

Note that \autoref{T:KLSandwich}.(c) and \autoref{T:KLSandwichp}.(a) 
give the same classification of simple $\hecke$-modules via $H$-reduction
\autoref{T:SandwichCMP}.

\begin{proof}
\textit{(a).} This is \cite[Theorem 5.14]{Je-cellular-pkl}.

\textit{(b).} 
We will show that $(\hecke,\pklbasis)$ cannot be a sandwich pair by contradicting \autoref{T:SandwichIsSandwich}.(d).
For $p>2$ in types $B_{2}$ and $C_{2}$ and for $p>3$ 
in type $G_{2}$ the $p$-KL basis is the 
usual KL basis, so \autoref{T:KLSandwich} 
applies. For the remaining cases we simply list the $p$-KL cells, which have 
$H$-cells of different sizes:
\begin{gather*}
\text{type $B_{2}$}\colon
\xy
(0,0)*{\begin{gathered}
\begin{tabular}{C}
\cellcolor{mydarkblue!25}
\scalebox{0.9}{$b_{1212}^{2}$}
\end{tabular}
\\[2pt]
\begin{tabular}{C|C}
\arrayrulecolor{tomato}
\scalebox{0.9}{$b_{121}^{2}$} & \cellcolor{spinach!25}\scalebox{0.9}{$b_{12}^{2}$}
\\
\hline
\cellcolor{spinach!25}\scalebox{0.9}{$b_{21}^{2}$} & \cellcolor{mydarkblue!25}\scalebox{0.9}{$b_{2}^{2},b_{212}^{2}$}
\\
\end{tabular}
\\[2pt]
\begin{tabular}{C}
\cellcolor{mydarkblue!25}\scalebox{0.9}{$b_{1}^{2}$}
\end{tabular}
\\[2pt]
\begin{tabular}{C}
\cellcolor{mydarkblue!25}\scalebox{0.9}{$b_{\emptyset}^{2}$}
\end{tabular}
\end{gathered}};
\endxy
\,,\quad
\text{type $C_{2}$}\colon
\xy
(0,0)*{\begin{gathered}
\begin{tabular}{C}
\cellcolor{mydarkblue!25}
\scalebox{0.9}{$b_{1212}^{2}$}
\end{tabular}
\\[2pt]
\begin{tabular}{C|C}
\arrayrulecolor{tomato}
\cellcolor{mydarkblue!25}
\scalebox{0.9}{$b_{1}^{2},b_{121}^{2}$} & \cellcolor{spinach!25}\scalebox{0.9}{$b_{12}^{2}$}
\\
\hline
\cellcolor{spinach!25}\scalebox{0.9}{$b_{21}^{2}$} & \scalebox{0.9}{$b_{212}^{2}$}
\\
\end{tabular}
\\[2pt]
\begin{tabular}{C}
\cellcolor{mydarkblue!25}\scalebox{0.9}{$b_{2}^{2}$}
\end{tabular}
\\[2pt]
\begin{tabular}{C}
\cellcolor{mydarkblue!25}\scalebox{0.9}{$b_{\emptyset}^{2}$}
\end{tabular}
\end{gathered}};
\endxy
\,,\\[0.5cm]
\text{type $G_{2}$}\colon
\xy
(0,0)*{\begin{gathered}
\begin{tabular}{C}
\cellcolor{mydarkblue!25}
\scalebox{0.9}{$b_{1212}^{2}$}
\end{tabular}
\\[2pt]
\begin{tabular}{C|C}
\arrayrulecolor{tomato}
\scalebox{0.9}{$b_{121}^{2},b_{12121}^{2}$} & \cellcolor{spinach!25}\scalebox{0.9}{$b_{1212}^{2}$}
\\
\hline
\cellcolor{spinach!25}\scalebox{0.9}{$b_{2121}^{2}$} & \cellcolor{mydarkblue!25}\scalebox{0.9}{$b_{212}^{2},b_{21212}^{2}$}
\\
\end{tabular}
\\[2pt]
\begin{tabular}{C|C}
\arrayrulecolor{tomato}
\scalebox{0.9}{$b_{1}^{2}$} & \cellcolor{spinach!25}\scalebox{0.9}{$b_{12}^{2}$}
\\
\hline
\cellcolor{spinach!25}\scalebox{0.9}{$b_{21}^{2}$} & \cellcolor{mydarkblue!25}\scalebox{0.9}{$b_{2}^{2}$}
\\
\end{tabular}
\\[2pt]
\begin{tabular}{C}
\cellcolor{mydarkblue!25}\scalebox{0.9}{$b_{\emptyset}^{2}$}
\end{tabular}
\end{gathered}};
\endxy
\,,\quad
\xy
(0,0)*{\begin{gathered}
\begin{tabular}{C}
\cellcolor{mydarkblue!25}
\scalebox{0.9}{$b_{121212}^{3}$}
\end{tabular}
\\[2pt]
\begin{tabular}{C|C}
\arrayrulecolor{tomato}
\scalebox{0.9}{$b_{121}^{3},b_{12121}^{3}$} & \cellcolor{spinach!25}\scalebox{0.9}{$b_{12}^{3},b_{1212}^{3}$}
\\
\hline
\cellcolor{spinach!25}\scalebox{0.9}{$b_{21}^{3},b_{2121}^{3}$} & \cellcolor{mydarkblue!25}\scalebox{0.9}{$b_{2}^{3},b_{212}^{3},b_{21212}^{3}$}
\\
\end{tabular}
\\[2pt]
\begin{tabular}{C}
\cellcolor{mydarkblue!25}\scalebox{0.9}{$b_{1}^{3}$}
\end{tabular}
\\[2pt]
\begin{tabular}{C}
\cellcolor{mydarkblue!25}\scalebox{0.9}{$b_{\emptyset}^{3}$}
\end{tabular}
\end{gathered}};
\endxy
\,.
\end{gather*}
Here we used $\dynkin[]B2$, 
$\dynkin[]C2$ and $\dynkin[reverse arrows]G2$ with the same 
labeling of vertices as in the dihedral case. Similar calculations 
verify the remaining cases. Note hereby that for $p$ big enough 
we can use \autoref{T:KLSandwich} since the $p$-KL basis 
and the KL agree for big enough primes, and there are only finitely 
many cases left to check which we verified with computer help 
(one can verify this without computer, but we did it by computer).

\textit{(c).} This follows from \cite[Proposition 4.2, part (7)]{JeWi-p-canonical} 
and \autoref{T:KLSandwich}.(a).
\end{proof}

\begin{Remark}\label{R:KLpKLConjecture}
For finite Coxeter systems of Dynkin type 
calculations suggest that $(\hecke,\pklbasis)$ 
is a sandwich pair if and only if $W$ is of type $A$.
\end{Remark}

\begin{Remark}\label{R:KLpKL}
The statements in \autoref{T:KLSandwichp} are reformulations 
of the vast literature on $p$-KL cells, see {\eg} 
\cite{JeWi-p-canonical}, \cite{Je-abc-pcells} or \cite{Je-cellular-pkl}, 
into the theory of sandwich cellular algebras.
\end{Remark}

%%%%%%%%%%%%%%%%%%%%%%%%%%%%%%%%%%%%%%%%%

\subsection{Kazhdan--Lusztig structure constants in dihedral type}\label{SS:KLDihedral}

%%%%%%%%%%%%%%%%%%%%%%%%%%%%%%%%%%%%%%%%%

Consider $\dynkin[Coxeter,gonality=\infty]I{}$ with nodes labeled $1$ and $2$.
The Coxeter group associated to this 
Coxeter diagram is the infinite dihedral group $\group[D]_{\infty}=\langle 
1,2|1^2=2^2=\emptyset\rangle$. 
Every element of $\group[D]_{\infty}$ has a unique reduced expression and
we write $k21$ and $k12$ for the reduced expressions $\dots 21$ 
and $\dots 12$ in $k$ symbols. For example $521=12121$.

The associated Hecke algebra $\hecke[I\infty]$ has a KL basis $\{b_{w}|w\in\group[D]_{\infty}\}$ 
(whose precise definition does not matter) with identity $b_{\emptyset}$.
Set $b_{0ab}=0$. The nonidentity multiplication rules are 
given by the \emph{(scaled) Clebsch--Gordan formula} where the steps size is two:
\begin{gather*}
b_{k12}b_{j12}
=
\begin{cases*}
[2]b_{(|k-j|+1)12}+\dots+[2]b_{(|k+j|-1)12}
&\text{j12=21\dots12},
\\
b_{|k-j|12}+2b_{(|k-j|+2)12}+\dots+2b_{(|k+j|-2)12}+b_{|k+j|12}
&\text{j12=12\dots 12}.
\end{cases*}
\end{gather*}
There are also similar formulas with $b_{j21}$ on the right 
and $b_{k21}$ on the left (in total four configurations).

Let $\group[D]_{n}=\langle 
1,2|1^2=2^2=(12)^{n}=\emptyset\rangle$ be the dihedral group of the 
$n$ gon. The longest element is $w_{0}=n12=n21$.
Let $[2]_{i}=(v^{i}+v^{-i})$.
With respect to the KL basis and its multiplication rules, the 
only change compare to $\group[D]_{\infty}$ is that expressions of the form (here $d>0$)
\begin{gather*}
b_{(n-d)12}+b_{(n+d)12}
\mapsto
[2]_{d}b_{w_0}
,\quad
b_{(n-d)21}+b_{(n+d)21}
\mapsto
[2]_{d}b_{w_0}
,
\end{gather*}
are replaced as indicated. This is the 
\emph{(scaled) truncated Clebsch--Gordan formula}.

\begin{Example}\label{E:KLDihedral}
For the infinite dihedral group we get
\begin{align*}
b_{1212}b_{21212}
&=
[2]b_{12}+[2]b_{1212}+[2]b_{121212}+[2]b_{12121212}
,
\\
b_{1212}b_{121212}
&=
b_{12}+2b_{1212}+2b_{121212}
+2b_{12121212}+b_{1212121212}
,
\end{align*}
from the Clebsch--Gordan formula. Moreover, 
\begin{align*}
b_{1212}b_{21212}
&=
[2]b_{12}+\colorbox{spinach!50}{\mystrut$[2]b_{1212}$}+\underline{[2]b_{121212}}+\colorbox{spinach!50}{\mystrut$[2]b_{12121212}$}
=[2]b_{12}+\big([2]_{3}+2[2]\big)b_{121212}
,
\\
b_{1212}b_{121212}
&=
\colorbox{orchid!50}{\mystrut $b_{12}$}+\colorbox{spinach!50}{\mystrut$2b_{1212}$}+\underline{2b_{121212}}+
\colorbox{spinach!50}{\mystrut$2b_{12121212}$}+\colorbox{orchid!50}{\mystrut $b_{1212121212}$}
=\big([2]_{4}+2[2]_{2}+[2]_{0}\big)b_{121212}
,
\end{align*}
are calculations for $n=6$. Here we used $[2]_{2}[2]+[2]=[2]_{3}+2[2]$ and $[2]_{0}=2$.
\end{Example}

The next lemma can be easily proven using the 
truncated Clebsch--Gordan formulas.

\begin{Lemma}\label{L:KLCells}
For the pair $(\hecke,\klbasis)$
the cell structure for $\hecke[In]$ is as in \autoref{E:KLDihedralOne}.\qed
\end{Lemma}

\begin{Lemma}\label{L:KLZModN}
Let $\K$ be an admissible field.
For the involutive sandwich 
pair $(\hecke,\klbasis)$ in odd dihedral type $I_{2}(n)$ we have $\sand[m]\cong\K[\Z/(\tfrac{n-1}{2})\Z]$.
\end{Lemma}

\begin{proof}
We focus on the $H$-cell containing $b_{1}$, and let 
$c_{121}=\frac{1}{[2]}b_{121}$. Note that $n\geq 3$ and $n$ odd.

Let $P_{0}(X)=1$, $P_{1}(X)=X$ and $P_{k+1}(X)=(X-1)P_{k}(X)-P_{k-1}(X)$ for $k>1$. The polynomial $P_{n-1}(X)$, known as 
\emph{multiplication by $[3]$}, is the characteristic polynomial of the graph
\begin{gather*}
\begin{tikzpicture}[anchorbase, scale=1]
\draw [thick] (0,0) to (2.825,0);
\draw [thick] (3.375,0) to (6.2,0);
\draw [thick] (0.95,0) to[out=135,in=180] (1,0.5) to[out=0,in=45] (1.05,0);
\draw [thick] (1.95,0) to[out=135,in=180] (2,0.5) to[out=0,in=45] (2.05,0);
\draw [thick] (4.15,0) to[out=135,in=180] (4.2,0.5) to[out=0,in=45] (4.25,0);
\draw [thick] (5.15,0) to[out=135,in=180] (5.2,0.5) to[out=0,in=45] (5.25,0);
\draw [thick] (6.15,0) to[out=135,in=180] (6.2,0.5) to[out=0,in=45] (6.25,0);
\node at (0,-0.01) {\Large $\bullet$};
\node at (1,-0.01) {\Large $\bullet$};
\node at (2,-0.01) {\Large $\bullet$};
\node at (3.08,-0.01) {\Large $\,\dots$};
\node at (4.2,-0.01) {\Large $\bullet$};
\node at (5.2,-0.01) {\Large $\bullet$};
\node at (6.2,-0.01) {\Large $\bullet$};
\node at (0,-0.3) {\tiny $[1]_{q}$};
\node at (1,-0.3) {\tiny $[3]_{q}$};
\node at (2,-0.3) {\tiny $[7]_{q}$};
\node at (4.15,-0.3) {\tiny $[m{-}4]_{q}$};
\node at (5.25,-0.3) {\tiny $[m{-}2]_{q}$};
\node at (6.2,-0.3) {\tiny $[m]_{q}$};
\end{tikzpicture}
,\quad
m=\tfrac{n-1}{2}.
\end{gather*}
This graph is the fusion graph of $\mathrm{SO}_{q}(3)$, 
the semisimplification of quantum $\mathfrak{so}_{3}(\C)$ tilting modules for $q=\exp(\pi i/n)$. The illustrated number below the vertices are 
the values of the Perron--Frobenius eigenvector of the graph.

The Clebsch--Gordan formulas show that this 
graph is also the action graph of $c_{121}$ on 
the $H$-cell containing $b_{1}$, up to rescaling of the KL basis, and we 
get $P_{(n-1)/2}(c_{121})=0$ via the identification of
$c_{121}$ with the Grothendieck class of the defining 
representation of $\mathfrak{so}_{3}(\C)$.

Hence, we have $\sand[m]\cong\K[X]/\big(P_{(n-1)/2}(X)\big)$.
It is known that the polynomial $P_{k}(X)$ has $k$ distinct roots. Thus, $\sand[m]\cong\K[X]/\big(P_{(n-1)/2}(X)\big)$
and the Chinese reminder theorem imply the lemma.
\end{proof}

\begin{Lemma}\label{L:KLSandwichMatrix}
Let $\K$ be an admissible field. 
For the involutive sandwich 
pair $(\hecke,\klbasis)$ the ranks of the sandwich
matrices are $1$ for the bottom, $2$ for all middle ones, and $1$ for the top.
\end{Lemma}

The involved calculations are a bit ugly and, for the sake 
of brevity, we only sketch the proof:

\begin{proof}[Sketch of a proof of \autoref{L:KLSandwichMatrix}]
The claim is immediate for the bottom and top $J$-cells. 
For $\jcell_{m}$ we first compute the general pairing matrix.
We scale such that we do not need to worry about $[2]$ and write $c_{w}=\frac{1}{[2]}b_{w}$.
Take $C=(C_{1},\dots)=(c_{1},c_{121},\dots,c_{12},c_{1212},\dots)$ 
and $R=(R_{1},\dots)=(c_{1},c_{121},\dots,c_{21},c_{2121},\dots)$
as index sets for columns and rows.
Next, we write down a matrix $M$ with the $i$-$j$ entry being $R_{i}C_{j}\pmod{\K\{b_{w_{0}}\}}$.
For $n=5$ we for example get
\begin{gather*}
M=
\scalebox{0.8}{$\begin{tabular}{C||CC|CC}
& c_{1} & c_{121} & c_{12} & c_{1212} \\
\hline
\hline
c_{1} & c_{1} & c_{121} & c_{12} & c_{1212} \\
c_{121} & c_{121} & c_{1}+c_{121} & c_{12}+c_{1212} & c_{12} \\
\hline
c_{21} & c_{21} & c_{21}+c_{2121} & c_{2}+c_{212} &  c_{212} \\
c_{2121} & c_{2121} & c_{21} & c_{212} & c_{2} \\
\end{tabular}$}
.
\end{gather*}
The isomorphism constructed in the proof of 
\autoref{L:KLZModN} implies that the matrix $M$, up to scaling, 
can be replaced by a matrix $N$ with four blocks corresponding to 
$\K[\Z/(\tfrac{n-1}{2})\Z]$. For example, for $n=5$ one gets
\begin{gather*}
N=
\scalebox{0.8}{$\begin{tabular}{C||CC|CC}
& d_{1} & d_{121} & d_{12} & d_{1212} \\
\hline
\hline
d_{1} & d_{1} & d_{121} & d_{12} & d_{1212} \\
d_{121} & d_{121} & d_{1} & d_{1212} & d_{12} \\
\hline
d_{21} & d_{21} & d_{2121} & d_{2} &  d_{212} \\
d_{2121} & d_{2121} & d_{21} & d_{212} & d_{2} \\
\end{tabular}$}
.
\end{gather*}
In this example, the units in the four blocks 
corresponding to $\K[\Z/2\Z]$ 
are $d_{1}$, $d_{1212}$, $d_{2121}$ and $d_{2}$.
In general the units are $d_{1}$, $d_{(n-1)12}$, $d_{(n-1)21}$ and $d_{2}$.

The idempotents for the simple $\K[\Z/(\tfrac{n-1}{2})\Z]$-modules 
are, of course, easy to write down, and it is then also easy to show that 
the ranks of the associated sandwich matrices are $\mrk(\smatrix[{m,i}])=2$.
For example, for $n=5$ the idempotents for the 
northwest corner are $e_{1}=d_{1}+d_{121}$ and $e_{2}=d_{1}-d_{121}$.
Let $e_{3}=d_{1212}+d_{12}$ be the idempotent for the southwest corner.
For $e_{1}$ we get the matrix
\begin{gather*}
Ne_{1}
=
\begin{pmatrix}
e_{1} & e_{1} & e_{1} & e_{1}
\\
e_{1} & e_{1} & e_{1} & e_{1}
\\
e_{3} & e_{3} & e_{3} & e_{3}
\\
e_{3} & e_{3} & e_{3} & e_{3}
\end{pmatrix}
.
\end{gather*}
The matrix is clearly of rank two, which shows that the sandwich 
matrix for $e_{1}$ is also of rank two.
The rank of the sandwich matrix for $e_{2}$ can be computed similarly.
\end{proof}

\begin{Example}\label{E:KLZModN}
With appropriate care one can use $H$-reduction 
\autoref{T:SandwichCMP} for arbitrary fields to classify 
simple $\group[D]_{n}$-modules for $n$ odd as well as their quantum counterparts.	

For example, take $n=5$ and $v=1$ and let us exclude $\chark=2$.
The polynomial from 
the proof of \autoref{L:KLZModN} defining the quotient is $P_{2}(X)=X^2-X-1$
whose roots 
are the golden ratio $\phi$, interpreted in $\K$, appearing in \autoref{Eq:KLCharTable} 
and its conjugate. (For the general case note that $\phi$ is $[3]_{q}$ 
specialized at $q=\exp(\pi i/n)$ with $n=5$.)
Thus, unless $\chark=5$, we can apply $H$-reduction for the middle $J$-cell and obtain two simple $\group[D]_{5}$-modules of that apex.
For $\chark=5$ one needs to distinguish the cases when 
$\K=\F_{25}$ and $\K=\F_{5}$, but one still gets 
the correct number of simple $\group[D]_{5}$-modules with the middle apex.

In general, we cannot rescale the KL basis and the polynomial 
defining the quotient is given by $P^{\prime}_{0}(X)=1$, $P^{\prime}_{1}(X)=X$ and $[2]P^{\prime}_{k+1}(X)=(X-[2])P^{\prime}_{k}(X)-[2]P^{\prime}_{k-1}(X)$ for $k>1$.
\end{Example}

\begin{Remark}\label{R:KLZModN}
As in \autoref{E:KLZModN}, it is not hard to discuss the case for 
general $\K$, but there is some annoying rescaling involved 
so we decide not to include this into this paper.
\end{Remark}

\begin{Remark}\label{R:KLDihedral}
The multiplication rules for the KL basis of
the dihedral group are well-known. The above, in particular, 
is a reformulation of \cite[Section 4]{duCl-positivity-finite-hecke}.
The application of sandwich cellularity to dihedral representation theory, {\eg} \autoref{L:KLCells}, \autoref{L:KLZModN} and \autoref{L:KLSandwichMatrix}, is new.
\end{Remark}

%%%%%%%%%%%%%%%%%%%%%%%%%%%%%%%%%%%%%%%%%

\section{Diagram algebras and sandwich cellularity}\label{S:DAlgebras}

%%%%%%%%%%%%%%%%%%%%%%%%%%%%%%%%%%%%%%%%%

We now discuss several examples of diagram algebras, all of which 
fit into the theory of sandwich cellular algebras.

%%%%%%%%%%%%%%%%%%%%%%%%%%%%%%%%%%%%%%%%%

\subsection{Some non-involutive diagram algebras}\label{SS:DAlgebrasTMon}

%%%%%%%%%%%%%%%%%%%%%%%%%%%%%%%%%%%%%%%%%

Let $\sym=\Aut\big(\{1,\dots,n\}\big)$ be the \emph{symmetric group} on the set $\{1,\dots,n\}$, and in this paper the \emph{planar symmetric group} is $\psym=\onemon$, independent of $n$. 
We now discuss two examples, which in some sense 
are analogs of the (planar) symmetric group in the theory of monoids.

\begin{Definition}\label{D:DAlgebrasTMonSymmetric}
The \emph{transformation monoid} $\tmon$ on the set $\{1,\dots,n\}$ 
is $\End\big(\{1,\dots,n\}\big)$.
\end{Definition}

Note that the subgroup of invertible 
elements $\group\subset\tmon$ is isomorphic to $\sym$ and we will 
identify $\sym$ as a subgroup of $\tmon$ in this way.

The elements of $\tmon$ can be written in one-line notation 
with $(ijk\dots)$ denoting the map $1\mapsto i$, $2\mapsto j$, 
$3\mapsto k$ {\etc} Alternatively we can model $\tmon$ 
using \emph{string diagrams} as follows.

We consider isotopy classes of diagrams of $2n$ points in the rectangle $[0,1]\times[0,1]$, with $n$ equally 
spaced points at the bottom and top.
For $(ijk\dots)$, connect the first bottom point with the $i$th on the top, 
the second with the $j$th on the top, the third with the $k$th 
on the top and so on.
Two diagrams represent the same element if and only if they 
represent the same map in $\tmon$. For example:
\begin{gather}\label{Eq:DAlgebrasTMonOneLine}
\begin{gathered}
(24138567)
\leftrightsquigarrow
\begin{tikzpicture}[anchorbase,scale=0.55]
\draw[usual] (4,-3) to[out=90,in=270] (3,-1);
\draw[usual] (3,-3) to[out=90,in=270] (2,-1);
\draw[usual] (2,-3) to[out=90,in=270] (1,-1);
\draw[usual] (1,-3) to[out=90,in=270] (4,-1);
\draw[usual] (0,-3) to[out=90,in=270] (-1,-1);
\draw[usual] (-1,-3) to[out=90,in=270] (-3,-1);
\draw[usual] (-2,-3) to[out=90,in=270] (0,-1);
\draw[usual] (-3,-3) to[out=90,in=270] (-2,-1);
\end{tikzpicture}
\,,\quad
(24637158)
\leftrightsquigarrow
\begin{tikzpicture}[anchorbase,scale=0.55]
\draw[usual] (4,-3) to[out=90,in=270] (4,-1);
\draw[usual] (3,-3) to[out=90,in=270] (1,-1);
\draw[usual] (2,-3) to[out=90,in=0] (1,-2.5) to (-2,-2.5) to[out=180,in=270] (-3,-1);
\draw[usual] (1,-3) to[out=90,in=270] (3,-1);
\draw[usual] (0,-3) to[out=90,in=270] (-1,-1);
\draw[usual] (-1,-3) to[out=90,in=270] (2,-1);
\draw[usual] (-2,-3) to[out=90,in=270] (0,-1);
\draw[usual] (-3,-3) to[out=90,in=270] (-2,-1);
\end{tikzpicture}
\,,\\
(23135555)
\leftrightsquigarrow
\begin{tikzpicture}[anchorbase,scale=0.55]
\draw[usual] (4,-3) to[out=90,in=270] (1,-1);
\draw[usual] (3,-3) to[out=90,in=270] (1,-1);
\draw[usual] (2,-3) to[out=90,in=270] (1,-1);
\draw[usual] (1,-3) to[out=90,in=270] (1,-1);
\draw[usual] (0,-3) to[out=90,in=270] (-1,-1);
\draw[usual] (-1,-3) to[out=90,in=270] (-3,-1);
\draw[usual] (-2,-3) to[out=90,in=270] (-1,-1);
\draw[usual] (-3,-3) to[out=90,in=270] (-2,-1);
\draw[usual,dot] (0,-1) to (0,-1.25);
\draw[usual,dot] (2,-1) to (2,-1.25);
\draw[usual,dot] (3,-1) to (3,-1.25);
\draw[usual,dot] (4,-1) to (4,-1.25);
\end{tikzpicture}
\,,\quad
(11335577)
\leftrightsquigarrow
\begin{tikzpicture}[anchorbase,scale=0.55]
\draw[usual] (4,-3) to[out=90,in=270] (3,-1);
\draw[usual] (3,-3) to[out=90,in=270] (3,-1);
\draw[usual] (2,-3) to[out=90,in=270] (1,-1);
\draw[usual] (1,-3) to[out=90,in=270] (1,-1);
\draw[usual] (0,-3) to[out=90,in=270] (-1,-1);
\draw[usual] (-1,-3) to[out=90,in=270] (-1,-1);
\draw[usual] (-2,-3) to[out=90,in=270] (-3,-1);
\draw[usual] (-3,-3) to[out=90,in=270] (-3,-1);
\draw[usual,dot] (-2,-1) to (-2,-1.25);
\draw[usual,dot] (0,-1) to (0,-1.25);
\draw[usual,dot] (2,-1) to (2,-1.25);
\draw[usual,dot] (4,-1) to (4,-1.25);
\end{tikzpicture}
\,.
\end{gathered}
\end{gather}
The dots on strings used in the pictures are reminders 
that there are top point not in the image of the displayed maps.
The first two diagrams are string diagrams for elements of $\sym[8]\subset\tmon[8]$.

\begin{Notation}\label{N:DAlgebraNames}
We have \emph{crossings}, \emph{merges} and \emph{top dots}, which are:
\begin{gather*}
\text{crossings}\colon
\begin{tikzpicture}[anchorbase,scale=0.55]
\draw[usual] (0,0) to[out=90,in=270] (1,1);
\draw[usual] (1,0) to[out=90,in=270] (0,1);
\end{tikzpicture}
\,,\quad
\text{merges}\colon
\begin{tikzpicture}[anchorbase,scale=0.55]
\draw[usual] (0,0) to[out=90,in=270] (0,1);
\draw[usual] (1,0) to[out=90,in=270] (0,1);
\end{tikzpicture}
\,,
\begin{tikzpicture}[anchorbase,scale=0.55]
\draw[usual] (0,0) to[out=90,in=270] (0,1);
\draw[usual] (1,0) to[out=90,in=270] (0,1);
\draw[usual] (2,0) to[out=90,in=270] (0,1);
\end{tikzpicture}
\,,\dots,\quad
\text{top dots}\colon
\begin{tikzpicture}[anchorbase,scale=0.55]
\draw[white] (0,0) to[out=90,in=270] (0,1);
\draw[usual,dot] (0,1) to (0,0.7);
\end{tikzpicture}
\,.
\end{gather*}
We will use this terminology throughout this section.
\end{Notation}

The diagrammatic multiplication $\circ$ in the transformation monoid is 
given by vertical gluing 
(and rescaling), using the convention from 
\autoref{N:SandwichGroundRing}.
For example,
\begin{gather*}
a\circ b
=
\begin{tikzpicture}[anchorbase,scale=0.55]
\draw[usual] (4,-1) to[out=90,in=270] (1,1);
\draw[usual] (3,-1) to[out=90,in=270] (1,1);
\draw[usual] (2,-1) to[out=90,in=270] (1,1);
\draw[usual] (1,-1) to[out=90,in=270] (1,1);
\draw[usual] (0,-1) to[out=90,in=270] (-1,1);
\draw[usual] (-1,-1) to[out=90,in=270] (-3,1);
\draw[usual] (-2,-1) to[out=90,in=270] (-1,1);
\draw[usual] (-3,-1) to[out=90,in=270] (-2,1);
\draw[usual,dot] (0,1) to (0,0.75);
\draw[usual,dot] (2,1) to (2,0.75);
\draw[usual,dot] (3,1) to (3,0.75);
\draw[usual,dot] (4,1) to (4,0.75);
\draw[very thick,densely dotted,tomato] (-3,-1) to (4,-1);
\draw[usual] (4,-3) to[out=90,in=270] (3,-1);
\draw[usual] (3,-3) to[out=90,in=270] (3,-1);
\draw[usual] (2,-3) to[out=90,in=270] (1,-1);
\draw[usual] (1,-3) to[out=90,in=270] (1,-1);
\draw[usual] (0,-3) to[out=90,in=270] (-1,-1);
\draw[usual] (-1,-3) to[out=90,in=270] (-1,-1);
\draw[usual] (-2,-3) to[out=90,in=270] (-3,-1);
\draw[usual] (-3,-3) to[out=90,in=270] (-3,-1);
\draw[usual,dot] (-2,-1) to (-2,-1.25);
\draw[usual,dot] (0,-1) to (0,-1.25);
\draw[usual,dot] (2,-1) to (2,-1.25);
\draw[usual,dot] (4,-1) to (4,-1.25);
\node at (-3.5,0) {$a$};
\node at (-3.5,-2) {$b$};
\end{tikzpicture}
=
\begin{tikzpicture}[anchorbase,scale=0.55]
\draw[usual] (4,-3) to[out=90,in=270] (1,-1);
\draw[usual] (3,-3) to[out=90,in=270] (1,-1);
\draw[usual] (2,-3) to[out=90,in=270] (1,-1);
\draw[usual] (1,-3) to[out=90,in=270] (1,-1);
\draw[usual] (0,-3) to[out=90,in=270] (-1,-2) to[out=90,in=270] (-3,-1);
\draw[usual] (-1,-3) to[out=90,in=270] (-1,-2) to[out=90,in=270] (-3,-1);
\draw[usual] (-2,-3) to[out=90,in=270] (-3,-2) to[out=90,in=270] (-2,-1);
\draw[usual] (-3,-3) to[out=90,in=270] (-3,-2) to[out=90,in=270] (-2,-1);
\draw[usual,dot] (-1,-1) to (-1,-1.25);
\draw[usual,dot] (0,-1) to (0,-1.25);
\draw[usual,dot] (2,-1) to (2,-1.25);
\draw[usual,dot] (3,-1) to (3,-1.25);
\draw[usual,dot] (4,-1) to (4,-1.25);
\end{tikzpicture}
\,.
\end{gather*}

\begin{Definition}\label{D:TMonPlanar}
The \emph{planar transformation monoid} $\ptmon$ on the set $\{1,\dots,n\}$ 
is the submonoid of $\tmon$ for which we can draw planar string diagrams 
without leaving the defining rectangle $[0,1]\times[0,1]$.
\end{Definition}

In other words, $\ptmon$ does not have any crossings.

\begin{Remark}\label{R:TMonPlanar}
In the semigroup literature 
$\tmon$ is also called \emph{full transformation monoid} and 
$\ptmon$ is called \emph{order-preserving transformation monoid}.	
\end{Remark}

\begin{Example}\label{E:TMonPlanar}
From the four elements displayed in \autoref{Eq:DAlgebrasTMonOneLine} 
only the southeast is in $\ptmon[8]$.
Moreover, we have 
$\ptmon[n]\cap\sym[n]=\{\monoidunit\}$, where 
$\monoidunit$ is the unit in $\tmon[n]$.
\end{Example}

Using the diagrammatic description we find:

\begin{Lemma}\label{L:DAlgebrasTMonFactor}
For $a\in\tmon$ there is a unique factorization of the form $a=\tau\circ\sigma_{\lambda}\circ\beta$ such that $\beta$ has a minimal 
number of crossings, $\tau$ 
has no crossings, $\beta$ contains no top dots, 
$\tau$ contains no merges 
and $\sigma_{\lambda}\in\sym[\lambda]$ for minimal $\lambda$.

Similarly for $a\in\ptmon$ but with $\sigma_{\lambda}\in\psym=\onemon$.
\end{Lemma}

\begin{proof}
The rearrangement
\begin{gather*}
a=
\begin{tikzpicture}[anchorbase,scale=0.55]
\draw[usual] (4,-3) to[out=90,in=270] (1,-1);
\draw[usual] (3,-3) to[out=90,in=270] (1,-1);
\draw[usual] (2,-3) to[out=90,in=270] (1,-1);
\draw[usual] (1,-3) to[out=90,in=270] (1,-1);
\draw[usual] (0,-3) to[out=90,in=270] (-1,-1);
\draw[usual] (-1,-3) to[out=90,in=270] (-3,-1);
\draw[usual] (-2,-3) to[out=90,in=270] (-1,-1);
\draw[usual] (-3,-3) to[out=90,in=270] (-2,-1);
\draw[usual,dot] (0,-1) to (0,-1.25);
\draw[usual,dot] (2,-1) to (2,-1.25);
\draw[usual,dot] (3,-1) to (3,-1.25);
\draw[usual,dot] (4,-1) to (4,-1.25);
\end{tikzpicture}
=
\begin{tikzpicture}[anchorbase,scale=0.55]
\draw[usual] (-3,0) to[out=90,in=270] (-3,2);
\draw[usual] (-2,0) to[out=90,in=270] (-1,2);
\draw[usual] (-1,0) to[out=90,in=270] (-2,2);
\draw[usual] (0,0) to[out=90,in=270] (-1,2);
\draw[usual] (1,0) to[out=90,in=270] (1,2);
\draw[usual] (2,0) to[out=90,in=270] (1,2);
\draw[usual] (3,0) to[out=90,in=270] (1,2);
\draw[usual] (4,0) to[out=90,in=270] (1,2);
\draw[very thick,densely dashed,mydarkblue] (-3,2) to (4,2);
\draw[usual] (-3,2) to[out=90,in=270] (-2,3);
\draw[usual] (-2,2) to[out=90,in=270] (-3,3);
\draw[usual] (-1,2) to[out=90,in=270] (-1,3);
\draw[usual] (1,2) to[out=90,in=270] (1,3);
\draw[very thick,densely dashed,mydarkblue] (-3,3) to (4,3);
\draw[usual] (-3,3) to[out=90,in=270] (-3,4);
\draw[usual] (-2,3) to[out=90,in=270] (-2,4);
\draw[usual] (-1,3) to[out=90,in=270] (-1,4);
\draw[usual,dot] (0,4) to (0,3.75);
\draw[usual] (1,3) to[out=90,in=270] (1,4);
\draw[usual,dot] (2,4) to (2,3.75);
\draw[usual,dot] (3,4) to (3,3.75);
\draw[usual,dot] (4,4) to (4,3.75);
\node at (-3.5,3.5) {$\tau$};
\node at (-3.5,2.5) {$\sigma_{4}$};
\node at (-3.5,1) {$\beta$};
\end{tikzpicture}
\end{gather*}
generalizes immediately.
\end{proof}

With respect to the factorization in \autoref{L:DAlgebrasTMonFactor}, 
we call $\lambda$ the number of \emph{through strands}, 
$\beta$ the \emph{bottom}, $\tau$ the \emph{top} and 
$\sigma_{\lambda}$ the \emph{middle} of $a$.
The middle in \autoref{L:DAlgebrasTMonFactor} also motivates:

\begin{Definition}\label{D:DAlgebrasTMonTCells}
We define the \emph{symmetric}
$\K\tmon=\K[\tmon]$ and the \emph{planar diagram algebra} $\K\ptmon=\K[\ptmon]$ 
associated the (planar) transformation monoid.
\end{Definition}

\begin{Notation}\label{N:DAlgebrasTMonTCells}
We write $\xmon$ for either $\tmon$ or $\ptmon$.
\end{Notation}

Let $\stirling{n}{\lambda}$ denote the $(n,\lambda)$th Stirling number (of the second kind), see also \autoref{R:DAlgebrasTMonStirling} below.

\begin{Proposition}\label{P:DAlgebrasTMonTCells}
We have the following for the pair $(\K\xmon,\xmon)$.
\begin{enumerate}

\item The $J$-cells of $\K\xmon$
are given by diagrams with a fixed number of through strands $\lambda$.
The $\leq_{lr}$-order is a total order and increases as the number of through strands decreases. That is,
\begin{gather*}
\Pcal=\{n<_{lr}n-1<_{lr}\dots<_{lr}1\}.
\end{gather*}

\item The left cells of $\K\xmon$ 
are given by diagrams where one fixes 
the bottom of the diagram, and similarly 
right cells are given by diagrams where one fixes 
the top of the diagram.
The $\leq_{l}$ and the $\leq_{r}$-order increases as the number of through strands decreases. Within $\jcell_{\lambda}$ we have
\begin{gather*}
\K\tmon\colon
\#\lcell_{\lambda}=\bstirling{n}{\lambda}
,\quad
\K\ptmon\colon
\#\lcell_{\lambda}=\binom{n-1}{\lambda-1}
,\quad
\K\xmon\colon
\#\rcell_{\lambda}=\binom{n}{\lambda}
.
\end{gather*}

\item Each $J$-cell of $\K\xmon$ is strictly idempotent, and
\begin{gather*}
\K\tmon\colon
\sand\cong\K[\sym[\lambda]]
,\quad
\K\ptmon\colon
\sand\cong\K[\psym[\lambda]]\cong\K.
\end{gather*}

\item The pair $(\K\xmon,\xmon)$ is a sandwich pair, that comes neither 
from a cellular nor an affine cellular algebra.

\item Let $\K$ be a field with $\chark=0$.
For $\lambda\neq 1$ we have
\begin{gather*}
\K\tmon\colon
\mrk(\gmatrix)=\binom{n}{\lambda}
,\quad
\K\ptmon\colon
\mrk(\gmatrix)=\binom{n-1}{\lambda-1}.
\end{gather*}
We also have $\mrk(\gmatrix[1])=1$ for both.

\end{enumerate}
\end{Proposition}

Note that the numerical data given in \autoref{P:DAlgebrasTMonTCells} 
is sufficient to determine the numbers and sizes of 
left, right, $J$ and $H$-cells by \autoref{L:SandwichNumericalCells}.

\begin{Remark}\label{R:DAlgebrasTMonStirling}
Very similar to binomials, the \emph{Stirling numbers} $\stirling{n}{\lambda}$
count certain combinatorial set partitions, and thus appear very often 
for diagram algebras. Explicitly, $\stirling{n}{\lambda}$ counts the number of ways to partition a set of $n$ labeled objects into $\lambda$ nonempty unlabeled subsets.
Not surprisingly, $\stirling{n}{\lambda}$ also admit 
a triangle description:	
\begin{gather*}
\begin{tikzpicture}[anchorbase,scale=1]
\node at (3,0) {\includegraphics[height=3cm]{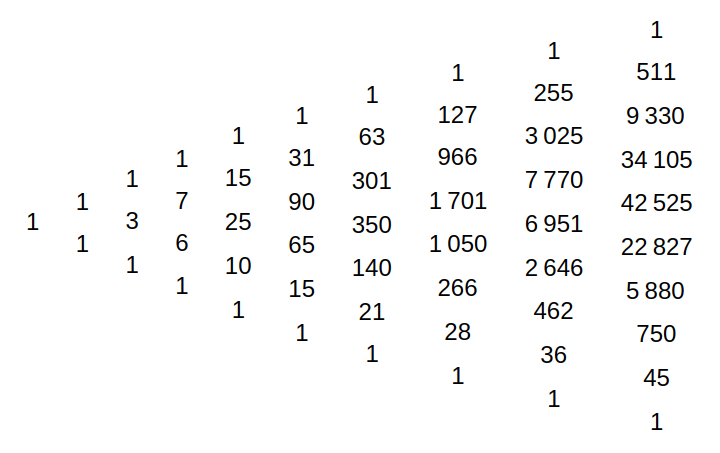}};
\draw[->] (0,0) to node[left]{$\lambda$} (0,-0.5);
\draw[->] (0,0) to node[above]{$n$} (0.5,0);
\draw[very thick,tomato] (2.3,0.275) to (2.5,0.175);
\draw[very thick,tomato] (2.3,0.05) to (2.5,0.15);
\end{tikzpicture}
.
\end{gather*}
The recursion is $\stirling{n}{\lambda}=\lambda\cdot
\stirling{n-1}{\lambda}+\stirling{n-1}{\lambda-1}$ with the 
starting conditions $\stirling{0}{0}=1$ 
and $\stirling{j}{0}=\stirling{0}{j}=0$ for $j>0$. 
For example, $90=3\cdot 25+15$.
\end{Remark}

\begin{proof}[Proof of \autoref{P:DAlgebrasTMonTCells}]
\textit{(a), (b) and (c)}. Except the counts 
for the left and the right cells, this follows from 
\autoref{L:DAlgebrasTMonFactor} and the construction of 
an idempotent for all $J$-cells. The construction of one idempotent 
for $\jcell_{\lambda}$ is easy:
let $e_{\lambda}$ 
be the idempotent that sends the first 
$n-\lambda+1$ entries to $1$ and fixes the remaining ones. For example, 
if $n=7$ and $\lambda=5$, then 
\begin{gather*}
e_{\lambda}=
\begin{tikzpicture}[anchorbase,scale=0.55]
\draw[usual] (0,0) to[out=90,in=270] (0,1);
\draw[usual] (0.5,0) to[out=90,in=270] (0,1);
\draw[usual] (1,0) to[out=90,in=270] (0,1);
\draw[usual] (1.5,0) to[out=90,in=270] (1.5,1);
\draw[usual] (2,0) to[out=90,in=270] (2,1);
\draw[usual] (2.5,0) to[out=90,in=270] (2.5,1);
\draw[usual] (3,0) to[out=90,in=270] (3,1);
\end{tikzpicture}
\,.
\end{gather*}
Clearly, $e_{\lambda}\in\jcell_{\lambda}$ and $e_{\lambda}e_{\lambda}=e_{\lambda}$.

To count the number of right cells, 
observe that they are determined by the number of top dot 
placements at the top of the diagram, and thus, we get 
$\binom{n}{n-\lambda}=\binom{n}{\lambda}$ of 
these. Analogously, for left cells 
we need to count bottom diagrams. For $\K\tmon$ 
this count is given by the Stirling numbers, by the construction of the
latter, see \autoref{R:DAlgebrasTMonStirling}.
In the planar case, the bottom 
diagrams are counted by  
compositions of $n$ into $\lambda$ parts with nonnegative 
components, and the number of
such compositions is known to be $\binom{n-1}{\lambda-1}$.
The claim follows.

\textit{(d).} This follows from (a), (b) and (c), \autoref{P:SandwichMonoid} 
and \autoref{L:SandwichInvolutive}.

\textit{(e).}
Recall from the above that each $H$-cell of $\K\xmon$ is determined by 
fixing a bottom and top diagram.
The idempotent $H$-cells in $\K\xmon$ are easy to describe:

Claim A. An $H$-cell of $\K\xmon$ is idempotent if and only if 
each bottom vertex that gets merged has an top dot right above it. 

Proof of Claim A. One computes
\begin{gather*}
\begin{tikzpicture}[anchorbase,scale=0.55]
\draw[usual] (1,0) to[out=90,in=270] (1,1);
\draw[usual] (2,0) to[out=90,in=270] (1,1);
\draw[usual] (3,0) to[out=90,in=270] (3,1);
\draw[usual,dot] (2,1) to (2,0.75);
\draw[very thick,densely dotted,tomato] (1,1) to (3,1);
\draw[usual] (1,1) to[out=90,in=270] (1,2);
\draw[usual] (2,1) to[out=90,in=270] (1,2);
\draw[usual] (3,1) to[out=90,in=270] (3,2);
\draw[usual,dot] (2,2) to (2,1.75);
\end{tikzpicture}
=
\begin{tikzpicture}[anchorbase,scale=0.55]
\draw[usual] (1,0) to[out=90,in=270] (1,1);
\draw[usual] (2,0) to[out=90,in=270] (1,1);
\draw[usual] (3,0) to[out=90,in=270] (3,1);
\draw[usual,dot] (2,1) to (2,0.75);
\end{tikzpicture}
\,,\quad
\begin{tikzpicture}[anchorbase,scale=0.55]
\draw[usual] (1,0) to[out=90,in=270] (1,1);
\draw[usual] (2,0) to[out=90,in=270] (1,1);
\draw[usual] (3,0) to[out=90,in=270] (2,1);
\draw[usual,dot] (3,1) to (3,0.75);
\draw[very thick,densely dotted,tomato] (1,1) to (3,1);
\draw[usual] (1,1) to[out=90,in=270] (1,2);
\draw[usual] (2,1) to[out=90,in=270] (1,2);
\draw[usual] (3,1) to[out=90,in=270] (2,2);
\draw[usual,dot] (3,2) to (3,1.75);
\end{tikzpicture}
\neq
\begin{tikzpicture}[anchorbase,scale=0.55]
\draw[usual] (1,0) to[out=90,in=270] (1,1);
\draw[usual] (2,0) to[out=90,in=270] (1,1);
\draw[usual] (3,0) to[out=90,in=270] (2,1);
\draw[usual,dot] (3,1) to (3,0.75);
\end{tikzpicture}
\,.
\end{gather*}
The claim follows by the immediate generalizing of these calculations.

Now, if $\lambda=1$, then (e) is 
immediate from the cell structure of 
$\K\xmon$, so let $\lambda>1$.

For $\K\ptmon$ we calculated the rank of a submatrix 
of $\gmatrix$ to be $\binom{n-1}{\lambda-1}$. To this end, 
we use induction on $n+\lambda$. For small values of $n+\lambda$
the claim is easily verified, see also \autoref{E:DAlgebrasTMonTMoreStrands} below. Assume inductively that the claim holds. Decompose the diagrams 
with $n$ strands and $\lambda$ through strands by fixing the leftmost 
strand to be a through strand:
\begin{gather}\label{Eq:DAlgebrasInduction}
\begin{tikzpicture}[anchorbase,scale=0.55]
\draw[usual] (1,0) to[out=90,in=270] (1,2.5);
\draw[mor] (1.9,0.5) rectangle node[pos=0.5]{$\binom{n-2}{\lambda}$}(5.1,1.5);
\draw[usual] (2,0) to[out=90,in=270] (2,0.5);
\draw[usual] (3,0) to[out=90,in=270] (3,0.5);
\draw[usual] (4,0) to[out=90,in=270] (4,0.5);
\draw[usual] (5,0) to[out=90,in=270] (5,0.5);
\draw[usual] (2,1.5) to[out=90,in=270] (1,2.5);
\draw[usual] (3,1.5) to[out=90,in=270] (1,2.5);
\draw[usual] (4,1.5) to[out=90,in=270] (2,2.5);
\draw[usual] (5,1.5) to[out=90,in=270] (3,2.5);
\draw[usual,dot] (4,2.5) to (4,2.25);
\draw[usual,dot] (5,2.5) to (5,2.25);
\end{tikzpicture}
\,,\quad
\begin{tikzpicture}[anchorbase,scale=0.55]
\draw[usual] (1,0) to[out=90,in=270] (1,2.5);
\draw[mor] (1.9,0.5) rectangle node[pos=0.5]{$\binom{n-2}{\lambda-1}$}(5.1,1.5);
\draw[usual] (2,0) to[out=90,in=270] (2,0.5);
\draw[usual] (3,0) to[out=90,in=270] (3,0.5);
\draw[usual] (4,0) to[out=90,in=270] (4,0.5);
\draw[usual] (5,0) to[out=90,in=270] (5,0.5);
\draw[usual] (2,1.5) to[out=90,in=270] (2,2.5);
\draw[usual] (3,1.5) to[out=90,in=270] (3,2.5);
\draw[usual] (4,1.5) to[out=90,in=270] (4,2.5);
\draw[usual] (5,1.5) to[out=90,in=270] (5,2.5);
\end{tikzpicture}
\,.
\end{gather}
(We pull top dots into the boxes if necessary.)
The diagrams in these boxes have $n-1$ strands 
with $\lambda$ and $\lambda-1$ through strands, respectively. 
Claim A shows that if we put these together in a submatrix 
the ranks add.
Thus, we get that a submatrix of rank $\binom{n-1}{\lambda-1}=\binom{n-2}{\lambda-1}+\binom{n-2}{\lambda-2}$ 
by induction.

For $\K\tmon$ the argument is the same, using induction, 
diagrams as in \autoref{Eq:DAlgebrasInduction}, Claim A and the formula 
$\binom{n}{\lambda}=\binom{n-1}{\lambda}+\binom{n-1}{\lambda-1}$ .
\end{proof}

\begin{Example}\label{E:DAlgebrasTMonTThreeStrands}
For $n=3$ we get
\begin{gather*}
\xy
(0,0)*{\begin{gathered}
\begin{tabular}{C}
\arrayrulecolor{tomato}
\cellcolor{mydarkblue!25}
\begin{tikzpicture}[anchorbase,scale=0.55]
\draw[usual] (1,0) to[out=90,in=270] (1,1);
\draw[usual] (2,0) to[out=90,in=270] (1,1);
\draw[usual] (3,0) to[out=90,in=270] (1,1);
\draw[usual,dot] (2,1) to (2,0.75);
\draw[usual,dot] (3,1) to (3,0.75);
\end{tikzpicture}
\\
\hline 
\cellcolor{mydarkblue!25}
\begin{tikzpicture}[anchorbase,scale=0.55]
\draw[usual] (1,0) to[out=90,in=270] (2,1);
\draw[usual] (2,0) to[out=90,in=270] (2,1);
\draw[usual] (3,0) to[out=90,in=270] (2,1);
\draw[usual,dot] (1,1) to (1,0.75);
\draw[usual,dot] (3,1) to (3,0.75);
\end{tikzpicture}
\\
\hline 
\cellcolor{mydarkblue!25}
\begin{tikzpicture}[anchorbase,scale=0.55]
\draw[usual] (1,0) to[out=90,in=270] (3,1);
\draw[usual] (2,0) to[out=90,in=270] (3,1);
\draw[usual] (3,0) to[out=90,in=270] (3,1);
\draw[usual,dot] (1,1) to (1,0.75);
\draw[usual,dot] (2,1) to (2,0.75);
\end{tikzpicture}
\end{tabular}
\\[1pt]
\begin{tabular}{C|C|C}
\arrayrulecolor{tomato}
\cellcolor{mydarkblue!25}
\begin{tikzpicture}[anchorbase,scale=0.55]
\draw[usual] (1,0) to[out=90,in=270] (1,1);
\draw[usual] (2,0) to[out=90,in=270] (2,1);
\draw[usual] (3,0) to[out=90,in=270] (2,1);
\draw[usual,dot] (3,1) to (3,0.75);
\end{tikzpicture}
,
\begin{tikzpicture}[anchorbase,scale=0.55]
\draw[usual] (1,0) to[out=90,in=270] (2,1);
\draw[usual] (2,0) to[out=90,in=270] (1,1);
\draw[usual] (3,0) to[out=90,in=270] (1,1);
\draw[usual,dot] (3,1) to (3,0.75);
\end{tikzpicture}
& \cellcolor{mydarkblue!25}
\begin{tikzpicture}[anchorbase,scale=0.55]
\draw[usual] (1,0) to[out=90,in=270] (1,1);
\draw[usual] (2,0) to[out=90,in=270] (2,1);
\draw[usual] (3,0) to[out=90,in=270] (1,1);
\draw[usual,dot] (3,1) to (3,0.75);
\end{tikzpicture}
,
\begin{tikzpicture}[anchorbase,scale=0.55]
\draw[usual] (1,0) to[out=90,in=270] (2,1);
\draw[usual] (2,0) to[out=90,in=270] (1,1);
\draw[usual] (3,0) to[out=90,in=270] (2,1);
\draw[usual,dot] (3,1) to (3,0.75);
\end{tikzpicture} 
& 
\begin{tikzpicture}[anchorbase,scale=0.55]
\draw[usual] (1,0) to[out=90,in=270] (2,1);
\draw[usual] (2,0) to[out=90,in=270] (2,1);
\draw[usual] (3,0) to[out=90,in=270] (1,1);
\draw[usual,dot] (3,1) to (3,0.75);
\end{tikzpicture}
,
\begin{tikzpicture}[anchorbase,scale=0.55]
\draw[usual] (1,0) to[out=90,in=270] (1,1);
\draw[usual] (2,0) to[out=90,in=270] (1,1);
\draw[usual] (3,0) to[out=90,in=270] (2,1);
\draw[usual,dot] (3,1) to (3,0.75);
\end{tikzpicture}
\\
\hline
\cellcolor{mydarkblue!25}
\begin{tikzpicture}[anchorbase,scale=0.55]
\draw[usual] (1,0) to[out=90,in=270] (1,1);
\draw[usual] (2,0) to[out=90,in=270] (3,1);
\draw[usual] (3,0) to[out=90,in=270] (3,1);
\draw[usual,dot] (2,1) to (2,0.75);
\end{tikzpicture}
,
\begin{tikzpicture}[anchorbase,scale=0.55]
\draw[usual] (1,0) to[out=90,in=270] (3,1);
\draw[usual] (2,0) to[out=90,in=270] (1,1);
\draw[usual] (3,0) to[out=90,in=270] (1,1);
\draw[usual,dot] (2,1) to (2,0.75);
\end{tikzpicture}
& 
\begin{tikzpicture}[anchorbase,scale=0.55]
\draw[usual] (1,0) to[out=90,in=270] (3,1);
\draw[usual] (2,0) to[out=90,in=270] (1,1);
\draw[usual] (3,0) to[out=90,in=270] (3,1);
\draw[usual,dot] (2,1) to (2,0.75);
\end{tikzpicture}
,
\begin{tikzpicture}[anchorbase,scale=0.55]
\draw[usual] (1,0) to[out=90,in=270] (1,1);
\draw[usual] (2,0) to[out=90,in=270] (3,1);
\draw[usual] (3,0) to[out=90,in=270] (1,1);
\draw[usual,dot] (2,1) to (2,0.75);
\end{tikzpicture}
& \cellcolor{mydarkblue!25}
\begin{tikzpicture}[anchorbase,scale=0.55]
\draw[usual] (1,0) to[out=90,in=270] (1,1);
\draw[usual] (2,0) to[out=90,in=270] (1,1);
\draw[usual] (3,0) to[out=90,in=270] (3,1);
\draw[usual,dot] (2,1) to (2,0.75);
\end{tikzpicture}
,
\begin{tikzpicture}[anchorbase,scale=0.55]
\draw[usual] (1,0) to[out=90,in=270] (3,1);
\draw[usual] (2,0) to[out=90,in=270] (3,1);
\draw[usual] (3,0) to[out=90,in=270] (1,1);
\draw[usual,dot] (2,1) to (2,0.75);
\end{tikzpicture}
\\
\hline
\begin{tikzpicture}[anchorbase,scale=0.55]
\draw[usual] (1,0) to[out=90,in=270] (2,1);
\draw[usual] (2,0) to[out=90,in=270] (3,1);
\draw[usual] (3,0) to[out=90,in=270] (3,1);
\draw[usual,dot] (1,1) to (1,0.75);
\end{tikzpicture}
,
\begin{tikzpicture}[anchorbase,scale=0.55]
\draw[usual] (1,0) to[out=90,in=270] (3,1);
\draw[usual] (2,0) to[out=90,in=270] (2,1);
\draw[usual] (3,0) to[out=90,in=270] (2,1);
\draw[usual,dot] (1,1) to (1,0.75);
\end{tikzpicture}
& \cellcolor{mydarkblue!25}
\begin{tikzpicture}[anchorbase,scale=0.55]
\draw[usual] (1,0) to[out=90,in=270] (3,1);
\draw[usual] (2,0) to[out=90,in=270] (2,1);
\draw[usual] (3,0) to[out=90,in=270] (3,1);
\draw[usual,dot] (1,1) to (1,0.75);
\end{tikzpicture}
,
\begin{tikzpicture}[anchorbase,scale=0.55]
\draw[usual] (1,0) to[out=90,in=270] (2,1);
\draw[usual] (2,0) to[out=90,in=270] (3,1);
\draw[usual] (3,0) to[out=90,in=270] (2,1);
\draw[usual,dot] (1,1) to (1,0.75);
\end{tikzpicture}
& \cellcolor{mydarkblue!25}
\begin{tikzpicture}[anchorbase,scale=0.55]
\draw[usual] (1,0) to[out=90,in=270] (2,1);
\draw[usual] (2,0) to[out=90,in=270] (2,1);
\draw[usual] (3,0) to[out=90,in=270] (3,1);
\draw[usual,dot] (1,1) to (1,0.75);
\end{tikzpicture}
,
\begin{tikzpicture}[anchorbase,scale=0.55]
\draw[usual] (1,0) to[out=90,in=270] (3,1);
\draw[usual] (2,0) to[out=90,in=270] (3,1);
\draw[usual] (3,0) to[out=90,in=270] (2,1);
\draw[usual,dot] (1,1) to (1,0.75);
\end{tikzpicture}
\end{tabular}
\\[1pt]
\begin{tabular}{C}
\arrayrulecolor{tomato}
\cellcolor{mydarkblue!25}
\begin{gathered}
\phantom{.}
\\[-0.35cm]
\begin{tikzpicture}[anchorbase,scale=0.55]
\draw[usual] (1,0) to[out=90,in=270] (1,1);
\draw[usual] (2,0) to[out=90,in=270] (2,1);
\draw[usual] (3,0) to[out=90,in=270] (3,1);
\end{tikzpicture}
,
\begin{tikzpicture}[anchorbase,scale=0.55]
\draw[usual] (1,0) to[out=90,in=270] (2,1);
\draw[usual] (2,0) to[out=90,in=270] (1,1);
\draw[usual] (3,0) to[out=90,in=270] (3,1);
\end{tikzpicture}
,
\begin{tikzpicture}[anchorbase,scale=0.55]
\draw[usual] (1,0) to[out=90,in=270] (1,1);
\draw[usual] (2,0) to[out=90,in=270] (3,1);
\draw[usual] (3,0) to[out=90,in=270] (2,1);
\end{tikzpicture}
\\[1pt]
\begin{tikzpicture}[anchorbase,scale=0.55]
\draw[usual] (1,0) to[out=90,in=270] (2,1);
\draw[usual] (2,0) to[out=90,in=270] (3,1);
\draw[usual] (3,0) to[out=90,in=270] (1,1);
\end{tikzpicture}
,
\begin{tikzpicture}[anchorbase,scale=0.55]
\draw[usual] (1,0) to[out=90,in=270] (3,1);
\draw[usual] (2,0) to[out=90,in=270] (1,1);
\draw[usual] (3,0) to[out=90,in=270] (2,1);
\end{tikzpicture}
,
\begin{tikzpicture}[anchorbase,scale=0.55]
\draw[usual] (1,0) to[out=90,in=270] (3,1);
\draw[usual] (2,0) to[out=90,in=270] (1,0.5) to[out=90,in=270] (2,1);
\draw[usual] (3,0) to[out=90,in=270] (1,1);
\end{tikzpicture}
\end{gathered}
\end{tabular}
\end{gathered}};
(-56,19)*{\jt};
(-56,-2.5)*{\jcell_{m}};
(-56,-21)*{\jb};
(58,19)*{\sand[t]\cong\K[\sym[1]]};
(58,-2.5)*{\sand[m]\cong\K[\sym[2]]};
(58,-21)*{\sand[b]\cong\K[\sym[3]]};
\endxy
\quad
,
\\[0.27cm]
\hline
\\[-0.27cm]
\xy
(0,0)*{\begin{gathered}
\begin{tabular}{C}
\arrayrulecolor{tomato}
\cellcolor{mydarkblue!25}(111) \\
\hline 
\cellcolor{mydarkblue!25}(222) \\
\hline 
\cellcolor{mydarkblue!25}(333)
\end{tabular}
\\[1pt]
\begin{tabular}{C|C|C}
\arrayrulecolor{tomato}
\cellcolor{mydarkblue!25}(122),(211) & \cellcolor{mydarkblue!25}(121),(212) & (221),(112)
\\
\hline
\cellcolor{mydarkblue!25}(133),(311) & (313),(131) & \cellcolor{mydarkblue!25}(113),(331)
\\
\hline
(233),(322) & \cellcolor{mydarkblue!25}(323),(232) & \cellcolor{mydarkblue!25}(223),(332)
\end{tabular}
\\[1pt]
\begin{tabular}{C}
\arrayrulecolor{tomato}
\cellcolor{mydarkblue!25}
\begin{gathered}
\phantom{.}
\\[-0.35cm]
(123),(213),(132)
\\[-1pt]
(231),(312),(321)
\end{gathered}
\end{tabular}
\end{gathered}};
(-46,13.8)*{\jt};
(-46,-2.5)*{\jcell_{m}};
(-46,-16.5)*{\jb};
(46,13.8)*{\sand[t]\cong\K[\sym[1]]};
(46,-2.5)*{\sand[m]\cong\K[\sym[2]]};
(46,-16.5)*{\sand[b]\cong\K[\sym[3]]};
\endxy
\quad
,
\end{gather*}
which display the cells of $\K\tmon[3]$ in diagrammatic notation 
and in one-line notation. For $\K\ptmon[3]$ one gets
\begin{gather*}
\xy
(0,0)*{\begin{gathered}
\begin{tabular}{C}
\arrayrulecolor{tomato}
\cellcolor{mydarkblue!25}
\begin{tikzpicture}[anchorbase,scale=0.55]
\draw[usual] (1,0) to[out=90,in=270] (1,1);
\draw[usual] (2,0) to[out=90,in=270] (1,1);
\draw[usual] (3,0) to[out=90,in=270] (1,1);
\draw[usual,dot] (2,1) to (2,0.75);
\draw[usual,dot] (3,1) to (3,0.75);
\end{tikzpicture}
\\
\hline 
\cellcolor{mydarkblue!25}
\begin{tikzpicture}[anchorbase,scale=0.55]
\draw[usual] (1,0) to[out=90,in=270] (2,1);
\draw[usual] (2,0) to[out=90,in=270] (2,1);
\draw[usual] (3,0) to[out=90,in=270] (2,1);
\draw[usual,dot] (1,1) to (1,0.75);
\draw[usual,dot] (3,1) to (3,0.75);
\end{tikzpicture}
\\
\hline 
\cellcolor{mydarkblue!25}
\begin{tikzpicture}[anchorbase,scale=0.55]
\draw[usual] (1,0) to[out=90,in=270] (3,1);
\draw[usual] (2,0) to[out=90,in=270] (3,1);
\draw[usual] (3,0) to[out=90,in=270] (3,1);
\draw[usual,dot] (1,1) to (1,0.75);
\draw[usual,dot] (2,1) to (2,0.75);
\end{tikzpicture}
\end{tabular}
\\[1pt]
\begin{tabular}{C|C}
\arrayrulecolor{tomato}
\cellcolor{mydarkblue!25}
\begin{tikzpicture}[anchorbase,scale=0.55]
\draw[usual] (1,0) to[out=90,in=270] (1,1);
\draw[usual] (2,0) to[out=90,in=270] (2,1);
\draw[usual] (3,0) to[out=90,in=270] (2,1);
\draw[usual,dot] (3,1) to (3,0.75);
\end{tikzpicture}
&
\begin{tikzpicture}[anchorbase,scale=0.55]
\draw[usual] (1,0) to[out=90,in=270] (1,1);
\draw[usual] (2,0) to[out=90,in=270] (1,1);
\draw[usual] (3,0) to[out=90,in=270] (2,1);
\draw[usual,dot] (3,1) to (3,0.75);
\end{tikzpicture}
\\
\hline
\cellcolor{mydarkblue!25}
\begin{tikzpicture}[anchorbase,scale=0.55]
\draw[usual] (1,0) to[out=90,in=270] (1,1);
\draw[usual] (2,0) to[out=90,in=270] (3,1);
\draw[usual] (3,0) to[out=90,in=270] (3,1);
\draw[usual,dot] (2,1) to (2,0.75);
\end{tikzpicture}
& \cellcolor{mydarkblue!25}
\begin{tikzpicture}[anchorbase,scale=0.55]
\draw[usual] (1,0) to[out=90,in=270] (1,1);
\draw[usual] (2,0) to[out=90,in=270] (1,1);
\draw[usual] (3,0) to[out=90,in=270] (3,1);
\draw[usual,dot] (2,1) to (2,0.75);
\end{tikzpicture}
\\
\hline
\begin{tikzpicture}[anchorbase,scale=0.55]
\draw[usual] (1,0) to[out=90,in=270] (2,1);
\draw[usual] (2,0) to[out=90,in=270] (3,1);
\draw[usual] (3,0) to[out=90,in=270] (3,1);
\draw[usual,dot] (1,1) to (1,0.75);
\end{tikzpicture}
& \cellcolor{mydarkblue!25}
\begin{tikzpicture}[anchorbase,scale=0.55]
\draw[usual] (1,0) to[out=90,in=270] (2,1);
\draw[usual] (2,0) to[out=90,in=270] (2,1);
\draw[usual] (3,0) to[out=90,in=270] (3,1);
\draw[usual,dot] (1,1) to (1,0.75);
\end{tikzpicture}
\end{tabular}
\\[1pt]
\begin{tabular}{C}
\arrayrulecolor{tomato}
\cellcolor{mydarkblue!25}
\begin{gathered}
\phantom{.}
\\[-0.35cm]
\begin{tikzpicture}[anchorbase,scale=0.55]
\draw[usual] (1,0) to[out=90,in=270] (1,1);
\draw[usual] (2,0) to[out=90,in=270] (2,1);
\draw[usual] (3,0) to[out=90,in=270] (3,1);
\end{tikzpicture}
\end{gathered}
\end{tabular}
\end{gathered}};
(-30,16)*{\jt};
(-30,-6)*{\jcell_{m}};
(-30,-22)*{\jb};
(33,16)*{\sand[t]\cong\K[\psym[1]]\cong\K};
(33,-6)*{\sand[m]\cong\K[\psym[2]]\cong\K};
(33,-22)*{\sand[b]\cong\K[\psym[3]]\cong\K};
\endxy
\quad
.
\end{gather*}
(We have not illustrated the cell structure in one-line notation.)
\end{Example}

Note that there is redundant information in the 
cell pictures in \autoref{E:DAlgebrasTMonTThreeStrands}. 
For example, the $H$-cells of $\K\tmon$ within $\jcell_{\lambda}$ are always 
of size $\lambda!$, so we actually only need to remember the idempotent $H$-cells. GAP's package Semigroups does exactly that:

\begin{Example}\label{E:DAlgebrasTMonTMoreStrands}
The cell structures of $\K\tmon$ and 
$\K\ptmon$ are of the form
\begin{gather*}
\tmon[3]\colon
\begin{tikzpicture}[anchorbase,scale=1]
\node at (0,0) {\raisebox{0.5cm}{\scalebox{-1}[-1]{\includegraphics[height=3.5cm]{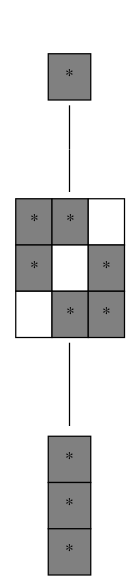}}}};
\end{tikzpicture}
,\quad
\tmon[4]\colon
\begin{tikzpicture}[anchorbase,scale=1]
\node at (0,0) {\raisebox{0.5cm}{\scalebox{-1}[-1]{\includegraphics[height=3.5cm]{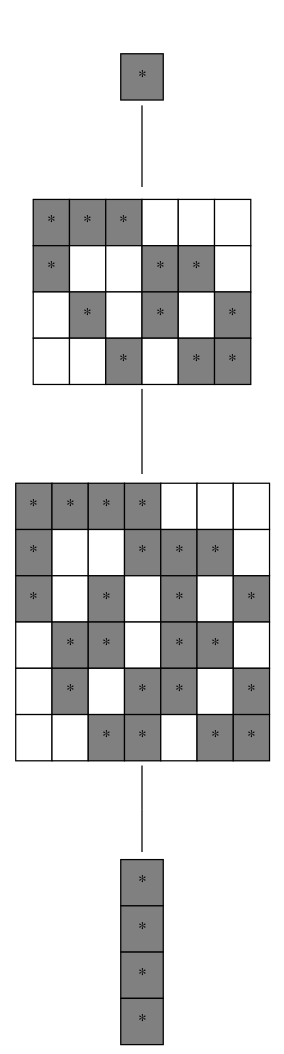}}}};
\end{tikzpicture}
,\quad
\tmon[5]\colon
\begin{tikzpicture}[anchorbase,scale=1]
\node at (0,0) {\raisebox{0.5cm}{\scalebox{-1}[-1]{\includegraphics[height=3.5cm]{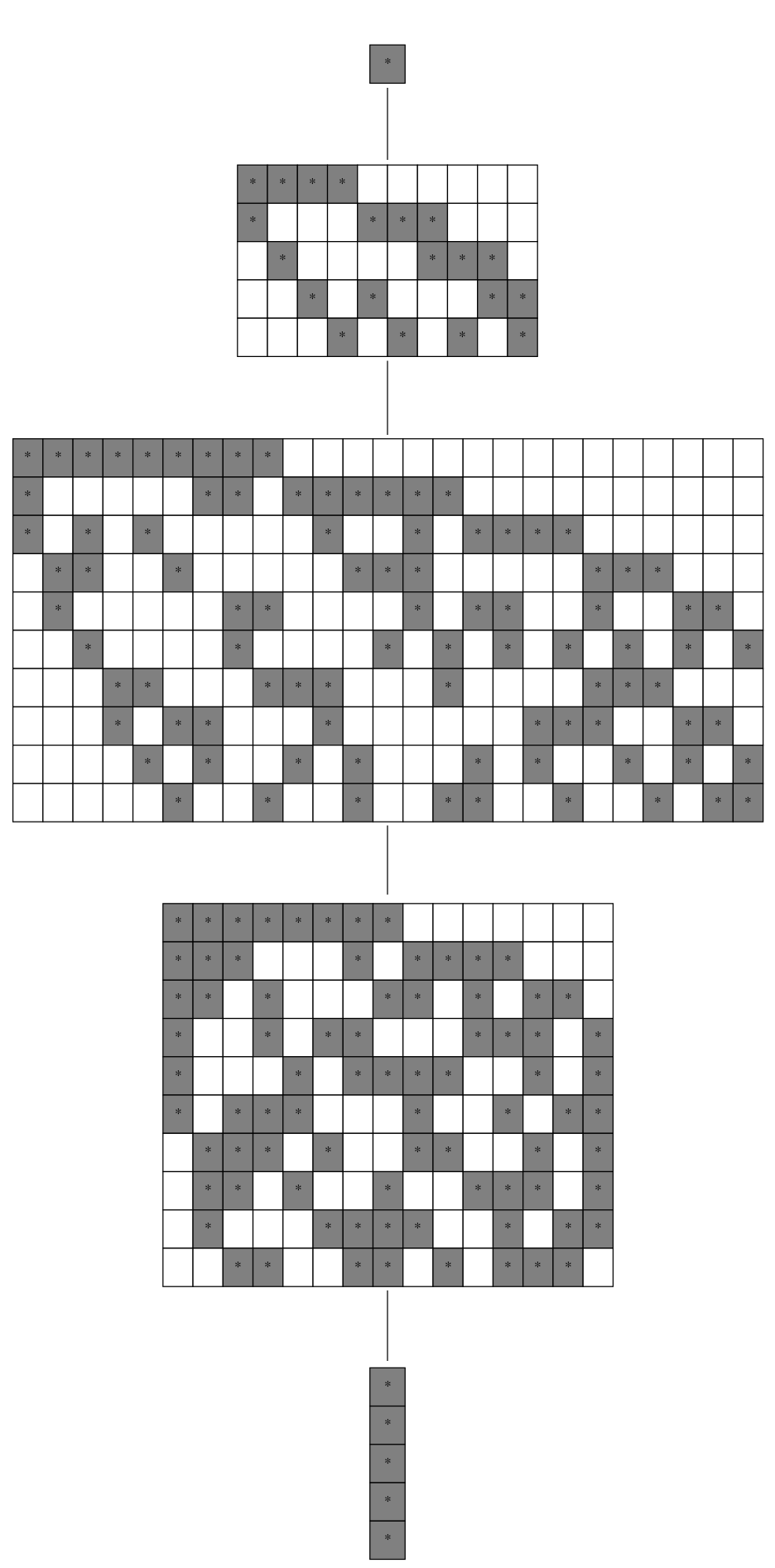}}}};
\end{tikzpicture}
,\quad
\ptmon[3]\colon
\begin{tikzpicture}[anchorbase,scale=1]
\node at (0,0) {\raisebox{0.5cm}{\scalebox{-1}[-1]{\includegraphics[height=3.5cm]{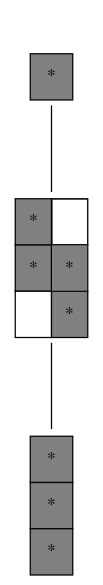}}}};
\end{tikzpicture}
,\quad
\ptmon[4]\colon
\begin{tikzpicture}[anchorbase,scale=1]
\node at (0,0) {\raisebox{0.5cm}{\scalebox{-1}[-1]{\includegraphics[height=3.5cm]{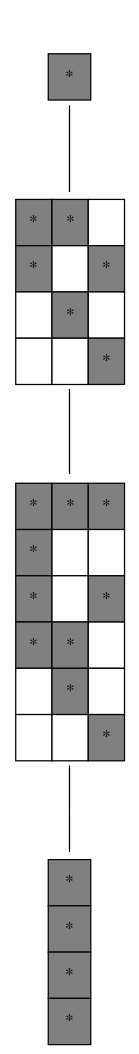}}}};
\end{tikzpicture}
,\quad
\ptmon[5]\colon
\begin{tikzpicture}[anchorbase,scale=1]
\node at (0,0) {\raisebox{0.5cm}{\scalebox{-1}[-1]{\includegraphics[height=3.5cm]{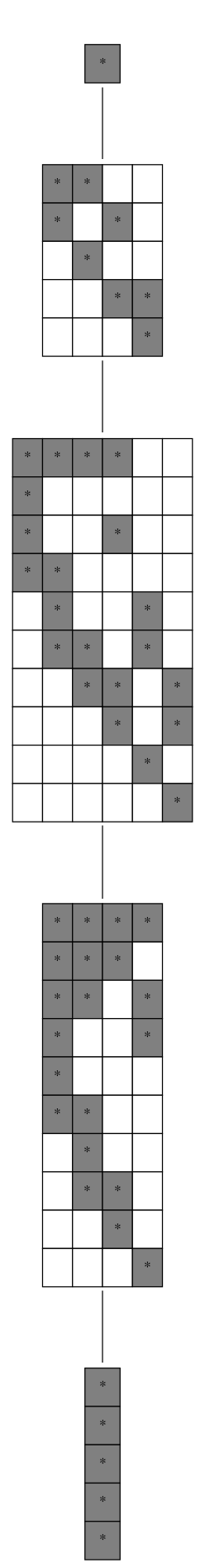}}}};
\end{tikzpicture}
,
\end{gather*}
where $n=3$, $n=4$ and $n=5$, respectively. The Gram matrices, except the topmost, 
have rank given by the number of rows, {\eg} the $(5\times 10)$-matrix for $\K\tmon[5]$ has rank $5$.

Note also that we need $\chark=0$, since, for example, 
the Gram matrix $\gmatrix[2]$ for $\K\tmon[3]$ has determinant $-2$. 
\end{Example}

\begin{Theorem}\label{P:DAlgebrasTMonTSimples}
Assume that $\K$ is a field, and consider 
the sandwich pair $(\K\xmon,\xmon)$.
The set of apexes for simple 
$\K\xmon$-modules is
$\apex=\{n<_{lr}n-1<_{lr}\dots<_{lr}1\}$, and
there are precisely $\big|\mathrm{P}\big(\lambda|\charkk\big)\big|$ simple  
$\K\tmon$-modules and one simple $\K\ptmon$-module for $\lambda\in\apex$.
\end{Theorem}

\begin{proof}
By \autoref{P:DAlgebrasTMonTCells} and $H$-reduction \autoref{T:SandwichCMP},
it suffices to identify the simple $\sand$-modules. This is 
immediate for $\K\ptmon$. For $\K\tmon$ we use the classical 
result that one can index simple $\sym[\lambda]$-modules by 
$\mathrm{P}\big(\lambda|\charkk\big)$, see for example \cite[Theorem 3.43]{Ma-hecke-schur}.
\end{proof}

\begin{Example}\label{E:DAlgebrasTMonTMoreStrandsTwo}
Continuing \autoref{E:DAlgebrasTMonTMoreStrands}, and 
let $\K$ be a field with
$\chark\nmid 5!$.
Then the number of simple $\K\tmon$-modules
for $n=5$ are given by $(7,5,3,2,1)$, where 
we ordered them by apex reading from the bottom to the top. 
For $\K\ptmon$ the sequence is $(1,1,1,1,1)$.
\end{Example}

By convention, the trivial and the sign $\sym[1]$-modules
are the same. This is relevant for the top cell in the following 
theorem.

\begin{Theorem}\label{T:DAlgebrasTMonTSimples}
Assume that $\K$ is a field with $\chark=0$ or $\chark\gg 0$ 
(meaning up to finitely many characteristics), and consider the 
sandwich pair $(\K\xmon,\xmon)$.
\begin{enumerate}

\item For $\K\tmon$, if $K$ is simple but not the sign $\sand$-module, 
then we have 
\begin{gather*}
\dim\big(\lmod[{\lambda,K}]\big)=\binom{n}{\lambda}\cdot\dim(K)
\end{gather*}
for $\lambda\in\apex$.
If $K$ is the sign $\sand$-module, 
then 
\begin{gather*}
\dim\big(\lmod[{\lambda,K}]\big)=\binom{n-1}{\lambda-1}
\end{gather*}
for $\lambda\in\apex$.

\item For $\K\ptmon$ and $\lambda\in\apex$ we have
\begin{gather*}
\dim\big(\lmod[\lambda]\big)=\binom{n-1}{\lambda-1}.
\end{gather*}

\end{enumerate}
\end{Theorem}

\begin{Remark}\label{R:DAlgebrasTMonTSimples}
Note that the dimensions of the simple $\sym[\lambda]$-modules
that appear in \autoref{T:DAlgebrasTMonTSimples} can be 
computed using well-known formulas since we assume that 
$\chark=0$.
Thus, \autoref{T:DAlgebrasTMonTSimples} gives an explicit 
description of the dimensions of the simple $\K\xmon[n]$-modules.

In contrast, the characteristic $p$ version of \autoref{T:DAlgebrasTMonTSimples} is much more difficult. In fact, to the best of our knowledge, the problem of determining the dimensions of simple $\K\xmon[n]$-modules in characteristic $p$ is still open
and for $\xmon[n]=\tmon$ that result would 
generalize the problem of finding the dimensions of the simple 
$\K\sym[\lambda]$-modules.
\end{Remark}

\begin{Remark}\label{R:DAlgebrasTMonTSimplesCrypto}
The dimensions of the simple $\K\xmon$-modules 
computed in \autoref{T:DAlgebrasTMonTSimples} are fairly large.
This could potentially be relevant for 
cryptographical purposes following \cite{KhSiTu-monoidal-cryptography}.
\end{Remark}

\begin{proof}[Proof of \autoref{T:DAlgebrasTMonTSimples}]
The calculations below use invertibility of certain, finitely many, numbers. These numbers are always invertible if $\chark=0$ but also for $\chark\gg 0$ since we only need to invert finitely many numbers.

\textit{(a).} The proof splits into two parts.
We will see why we get two different cases using 
the diagrammatic description of 
the primitive idempotents of $\K[\sym[\lambda]]$ due to (Gyoja--)Aiston(--Morton) 
\cite{AiMo-hecke-idempotents}, see 
also \cite[Definition 2.26]{TuVaWe-super-howe}.
Fix $\mu\in\mathrm{P}(\lambda|\infty)$, seen as a Young diagram,
with $\mu_{1}^{r},\dots,\mu_{k}^{r}$ rows and 
$\mu_{1}^{c},\dots,\mu_{l}^{c}$ columns. Note that $\mu_{1}^{r}+\dots+\mu_{k}^{r}=n=\mu_{1}^{c}+\dots+\mu_{l}^{c}$. The full symmetrizer for $\mu_{i}^{r}$ is the sum of all symmetric group elements 
on the respective strands, and the full antisymmetrizer for $\mu_{i}^{c}$ is the signed sum of all symmetric group elements on the respective strands.
Define $e_{row}(\mu)$ and $e_{col}(\mu)$ to be the tensor product (horizontal juxtaposition of strands) 
of the full symmetrizers and full antisymmetrizers associated to the 
rows respectively columns of $\mu$. Let $w$ be any shortest presentation 
(with respect to simple transpositions) of
the permutation that permutes 
the row standard filling of $\mu$ to the column standard filling of $\mu$.
Then, for some nonzero scalar $s\in\K$, $e_{\mu}=\frac{1}{s}\cdot e_{row}(\mu)\cdot(\id\otimes w)\cdot e_{col}(\mu)\cdot(\id\otimes w^{-1})$ is a primitive idempotent 
in the group algebra of $\sym[\lambda]$ projecting to the simple $\K[\sym[\lambda]]$-module $K=K(\mu)$ associated to $\mu$.

As not unusual in these types of diagrammatics, we draw symmetrizers as reddish shaded boxes, and antisymmetrizer as greenish shaded boxes labeled by their number of strands and an $s$ or an $a$ to distinguish the boxes.
Let $w$ be the permutation 
$(23)(45)(34)$ in this presentation. Illustrating $w$ and $w^{-1}$ also as shaded boxes, we get
\begin{gather*}
\mu=\ydiagram{2,2,1}\colon
e_{\mu}=\frac{1}{s}\cdot
\begin{tikzpicture}[anchorbase,scale=0.55]
\draw[mor2] (-0.1,0) rectangle node[pos=0.5]{\scalebox{0.7}{$3s$}}(1.6,0.5);
\draw[mor2] (1.9,0) rectangle node[pos=0.5]{\scalebox{0.7}{$2s$}}(3.1,0.5);
\draw[usual] (0,-0.5) to (0,0);
\draw[mor] (0.4,-0.5) rectangle node[pos=0.5]{\scalebox{0.7}{$w$}}(3.1,0);
\draw[mor3] (-0.1,-1) rectangle node[pos=0.5]{\scalebox{0.7}{$2a$}}(0.9,-0.5);
\draw[mor3] (1.1,-1) rectangle node[pos=0.5]{\scalebox{0.7}{$2a$}}(2.4,-0.5);
\draw[mor3] (2.6,-1) rectangle node[pos=0.5]{\scalebox{0.7}{$1a$}}(3.1,-0.5);
\draw[usual] (0,-1.5) to (0,-1);
\draw[mor] (0.4,-1.5) rectangle node[pos=0.5]{\scalebox{0.7}{$w^{-1}$}}(3.1,-1);
\end{tikzpicture}
\,,
\end{gather*}
as the primitive idempotent associated to $\mu$.
We have the following relations:
\begin{gather}\label{Eq:DAlgebraRels}
\begin{tikzpicture}[anchorbase,scale=0.55]
\draw[usual] (0,0.5) to[out=90,in=270] (1,1.5);
\draw[usual] (1,0.5) to[out=90,in=270] (0,1.5);
\draw[white] (0,0) to (1,-1);
\draw[white] (1,0) to (0,-1);
\draw[mor2] (-0.2,0) rectangle node[pos=0.5]{\scalebox{0.7}{$s$}}(1.2,0.5);
\end{tikzpicture}
=
\begin{tikzpicture}[anchorbase,scale=0.55]
\draw[mor2] (-0.2,0) rectangle node[pos=0.5]{\scalebox{0.7}{$s$}}(1.2,0.5);
\end{tikzpicture}
\,,\quad
\begin{tikzpicture}[anchorbase,scale=0.55]
\draw[usual] (0,0.5) to[out=90,in=270] (0,1.5);
\draw[usual] (1,0.5) to[out=90,in=270] (0,1.5);
\draw[white] (0,0) to (1,-1);
\draw[white] (1,0) to (0,-1);
\draw[mor2] (-0.2,0) rectangle node[pos=0.5]{\scalebox{0.7}{$s$}}(1.2,0.5);
\end{tikzpicture}
=
\begin{tikzpicture}[anchorbase,scale=0.55]
\draw[usual] (0,0) to[out=90,in=270] (0,1);
\draw[usual] (1,0) to[out=90,in=270] (0,1);
\draw[white] (0,0) to (1,-1);
\draw[white] (1,0) to (0,-1);
\end{tikzpicture}
\,,\quad
\begin{tikzpicture}[anchorbase,scale=0.55]
\draw[usual] (0,0.5) to[out=90,in=270] (1,1.5);
\draw[usual] (1,0.5) to[out=90,in=270] (0,1.5);
\draw[white] (0,0) to (1,-1);
\draw[white] (1,0) to (0,-1);
\draw[mor3] (-0.2,0) rectangle node[pos=0.5]{\scalebox{0.7}{$a$}}(1.2,0.5);
\end{tikzpicture}
=
-
\begin{tikzpicture}[anchorbase,scale=0.55]
\draw[mor3] (-0.2,0) rectangle node[pos=0.5]{\scalebox{0.7}{$a$}}(1.2,0.5);
\end{tikzpicture}
,\quad
\begin{tikzpicture}[anchorbase,scale=0.55]
\draw[usual] (0,0.5) to[out=90,in=270] (0,1.5);
\draw[usual] (1,0.5) to[out=90,in=270] (0,1.5);
\draw[white] (0,0) to (1,-1);
\draw[white] (1,0) to (0,-1);
\draw[mor3] (-0.2,0) rectangle node[pos=0.5]{\scalebox{0.7}{$a$}}(1.2,0.5);
\end{tikzpicture}
=
0.
\end{gather}
That is, the (anti)symmetrizers eat crossings, while merges eat symmetrizers 
but annihilate antisymmetrizers. As we will see, the case where there 
is no symmetrizer box is special as all merges are annihilated.

\textit{Case 1: $K$ is not the sign $\sand$-module.} Assume that $K$ is simple but not the sign $\sand$-module.
We will show that $\dmod[\lambda]e_{K}$ is a simple $\K\tmon$-module.
Note that this implies 
$\dim\big(\lmod[{\lambda,K}]\big)=\binom{n}{\lambda}\cdot\dim(K)$, which is what we wanted to prove.

To prove that $\dmod[\lambda]e_{K}$ is a simple $\K\tmon$-module, 
we use induction on $n-\lambda$. If $n-\lambda=0$, 
then we are in the bottom cell and the claim is clear. So let $n-\lambda=k$ 
and assume that we have proven the claim for $n-\lambda<k$. 
Define idempotents
\begin{gather*}
e_{i,j}=
\begin{tikzpicture}[anchorbase,scale=0.55]
\draw[usual] (0,0)node[below]{$i$} to[out=90,in=270] (0,1)node[above]{$i$};
\draw[usual] (2,0)node[below]{$j$} to[out=90,in=270] (0,1)node[above]{$i$};
\draw[usual] (-0.5,0) to[out=90,in=270] (-0.5,1);
\draw[usual] (0.5,0) to[out=90,in=270] (0.5,1);
\draw[usual] (1,0) to[out=90,in=270] (1,1);
\draw[usual] (1.5,0) to[out=90,in=270] (1.5,1);
\draw[usual] (2.5,0) to[out=90,in=270] (2.5,1);
\end{tikzpicture}
\,,\quad
1\leq i<j\leq n,
\end{gather*}
that merge the $i$th and $j$th strands and are the identity everywhere else.
Note that the top of the diagrams for the $e_{i,j}$ have $n-1$ strands.

For $x\in\dmod[\lambda]e_{K}$ suppose that $e_{i,j}x\neq 0$ for some 
$e_{i,j}$. Using the permutation $\sigma_{i,j}$ of the $i$th and $j$th strand, this implies that $y=e_{1,n}\sigma_{i,j}x\in e_{1,n}\dmod[\lambda]e_{K}$ 
is nonzero. Since $\K\tmon[{n-1}]\cong e_{1,n}\K\tmon e_{1,n}$, we know by induction that $e_{1,n}\dmod[\lambda]e_{K}$ is 
simple as a $\K\tmon[{n-1}]$-module and
we get that $\dmod[\lambda]e_{K}\subset\K\tmon e_{1,n}\dmod[\lambda]e_{K}\subset\K\tmon x$, which shows that 
$\dmod[\lambda]e_{K}$ is a simple $\K\tmon$-module.

From \autoref{Eq:DAlgebraRels} it follows that at least one $e_{i,j}$ does not annihilate $\dmod[\lambda]e_{K}$: in the case 
where $K=K(\mu)$ is not the sign $\sand$-module the idempotent
$e_{K}$ has at least one nontrivial symmetrizer box at the top and hitting it 
will produce a nonzero diagram.

\textit{Case 2: $K$ is the sign $\sand$-module.}
The case where $K$ is the sign $\sand$-module is special and 
better to be analyzed by different means. 
Let $e_{\lambda}$ be the idempotent as in the proof of 
\autoref{P:DAlgebrasTMonTCells}.
One easily sees that $e_{\lambda}\tmon e_{\lambda}\cong\sym[\lambda]$.
Let $\module[M]=
\{(x_{1},\dots,x_{n})\in\K^{n}|x_{1}+\dots+x_{n}=0\}$ denote the 
$\K\tmon$-module with action induced from the natural action of $\K\tmon$ 
on $\K^{n}$. The space $\module[M]$ is a $e_{\lambda}\tmon e_{\lambda}\cong\sym[\lambda]$-module by restriction and hence, we get a 
$e_{\lambda}\K\tmon e_{\lambda}\cong\K[\sym[\lambda]]$-module
$\bigwedge^{\lambda-1}\module[M]$. It is known that $\bigwedge^{\lambda-1}\module[M]$ is a simple 
$\K[\sym[\lambda]]$-module, and thus, is 
also a simple $\K\tmon$-module. Its dimension is $\binom{n-1}{\lambda-1}$ 
and $\bigwedge^{\lambda-1}\module[M]$ has apex $\lambda$, by construction. 
Since we already know all the other simple $\K\tmon$-module with apex $\lambda$, 
and they are of bigger dimensions, we conclude that 
$\lmod[{\lambda,K}]\cong\bigwedge^{\lambda-1}\module[M]$.

\textit{(b).} This follows directly from \autoref{P:SandwichCellsGram} 
and \autoref{P:DAlgebrasTMonTCells}
\end{proof}	

In the proof of \autoref{T:DAlgebrasTMonTSimples} we could have 
alternatively computed the sandwich matrices 
as the following example shows.

\begin{Example}\label{E:DAlgebrasTMonTMoreStrandsThree}
Continuing \autoref{E:DAlgebrasTMonTMoreStrandsTwo}, and let $\K$ 
be such that $\chark\nmid 3!$.
Then the dimensions of simple $\K\tmon[3]$-modules
are given by $1$, $2$, $1$ for the bottom apex, 
$3$, $2$ for the middle apex and $1$ for the top apex.
The two sandwich matrices of ranks $3$ and $2$ for the middle apex are
\begin{gather*}
\smatrix[{m,triv}]
=
\begin{psmallmatrix}
e_{triv} & e_{triv} & e_{triv} & e_{triv} & 0 & 0
\\
e_{triv} & e_{triv} & e_{triv} & e_{triv} & 0 & 0
\\
e_{triv} & e_{triv} & 0 & 0 & e_{triv} & e_{triv}
\\
e_{triv} & e_{triv} & 0 & 0 & e_{triv} & e_{triv}
\\
0 & 0 & e_{triv} & e_{triv} & e_{triv} & e_{triv}
\\
0 & 0 & e_{triv} & e_{triv} & e_{triv} & e_{triv}
\end{psmallmatrix}
,\quad
\smatrix[{m,sign}]
=
\begin{psmallmatrix}
e_{sign} & -e_{sign} & e_{sign} & -e_{sign} & 0 & 0
\\
-e_{sign} & e_{sign} & -e_{sign} & e_{sign} & 0 & 0
\\
e_{sign} & -e_{sign} & 0 & 0 & -e_{sign} & e_{sign}
\\
-e_{sign} & e_{sign} & 0 & 0 & e_{sign} & -e_{sign}
\\
0 & 0 & -e_{sign} & e_{sign} & -e_{sign} & e_{sign}
\\
0 & 0 & e_{sign} & -e_{sign} & e_{sign} & -e_{sign}
\end{psmallmatrix}
,
\end{gather*}
where $e_{triv}$ and $e_{sign}$ are the idempotents for the 
trivial and the sign $\sand[2]$-module, respectively.

For $\K\ptmon[3]$ there is actually no restriction on the field
since the Gram matrices are
\begin{gather*}
\gmatrix[b]=\begin{pmatrix}1\end{pmatrix}
,\quad
\gmatrix[m]=\begin{pmatrix}1 & 0\\1 & 1\\0 & 1\end{pmatrix}
,\quad
\gmatrix[t]=\begin{pmatrix}1\\1\\1\end{pmatrix}
\end{gather*} 
and we get $1$, $2$, $1$ as the simple dimensions.
\end{Example}

%\begin{Remark}\label{R:DAlgebrasTMonPTMon}
%With respect to $\K\ptmon$, recall 
%from \autoref{P:DAlgebrasTMonTCells} that we 
%have $\binom{n-1}{\lambda-1}$ left cells and 
%$\binom{n}{\lambda}$ right cells within $\jcell_{\lambda}$. 
%The identity $\binom{n}{\lambda}=\binom{n-1}{\lambda}+\binom{n-1}{\lambda-1}$
%gives numerically that the number of rows in the matrix for $\jcell_{\lambda}$ equals the number of columns in the matrix for $\jcell_{\lambda+1}$ (the cell at the bottom of $\jcell_{\lambda}$) plus the number of columns in the matrix for $\jcell_{\lambda}$.
%This suggests that the cell module $\dmod[\lambda]$ has a Jordan--H{\"o}lder 
%filtration of length two with $\lmod[{\lambda+1}]$ as a submodule and $\lmod[{\lambda}]$ as a quotient. 
%One can check in small examples that this is indeed the case.
%\end{Remark}

\begin{Remark}\label{R:DAlgebraQuantum}
The quantum version of the above discussion simply replaces 
crossings by \emph{over} and \emph{undercrossings}:
\begin{gather*}
\begin{tikzpicture}[anchorbase,scale=0.55]
\draw[usual] (0,0) to[out=90,in=270] (1,1);
\draw[usual] (1,0) to[out=90,in=270] (0,1);
\end{tikzpicture}
\rightsquigarrow
\begin{tikzpicture}[anchorbase,scale=0.55]
\draw[usual] (1,0) to[out=90,in=270] (0,1);
\draw[usual,crossline] (0,0) to[out=90,in=270] (1,1);
\end{tikzpicture}
\text{ or }
\begin{tikzpicture}[anchorbase,scale=0.55]
\draw[usual] (0,0) to[out=90,in=270] (1,1);
\draw[usual,crossline] (1,0) to[out=90,in=270] (0,1);
\end{tikzpicture}
.
\end{gather*}
In this case the sandwiched algebras are the Hecke algebras 
of type $A$ as in \autoref{S:KL}. Otherwise there is no difference.
\end{Remark}

\begin{Remark}\label{R:DAlgebrasTMonHistory}
The transformation monoid and its 
planar counterpart are around for donkey's 
years, and so are their 
diagrammatic descriptions. The diagrammatic 
incarnations appear, for example, in the formulation of 
categories of (planar) transformations, 
see {\eg} \cite[Figure 8]{EaRu-congruences-categories}.

The representation theory of $\K\tmon$ was studied 
since the early days of monoid representation theory.
Often relying on versions of \autoref{T:SandwichCMP}, but in a less 
general setting with the focus on monoids.
For example, the classification of simple $\K\tmon$-modules 
is given in \cite[Section 5.3]{St-rep-monoid}, but in a different language 
and not using sandwich cellularity. 
That section also lists a few original references, starting with 
\cite{HeZu-tmon}.

We do not know any reference 
for the representation theory of 
the planar transformation monoid.
This discussion appears to be new, and 
is easy from the perspective 
of sandwich cellularity.
\end{Remark}

%%%%%%%%%%%%%%%%%%%%%%%%%%%%%%%%%%%%%%%%%

\subsection{Some involutive diagram algebras}\label{SS:DAlgebras}

%%%%%%%%%%%%%%%%%%%%%%%%%%%%%%%%%%%%%%%%%

In the previous section we found 
the cells and parameterized 
the simple $\K\xmon$-modules, for $\xmon$ being $\tmon$ or $\ptmon$,
using their diagrammatic incarnation.
The same strategy works, by its very construction, for a wide range of 
diagram monoids and algebras.

We list now a few such algebras 
(the reader unfamiliar with these 
is referred to {\eg} \cite{HaRa-partition-algebras} 
for more details). These algebras have in common 
that their multiplication is defined via an underlying monoid, up to evaluations of closed components.

\begin{Notation}\label{N:DAlgebrasPlanar}
As before, \emph{planar} means 
the submonoids having only planar diagrams of the same type, while 
the other listed diagrammatic 
descriptions are \emph{symmetric}. Below we will give one prototypical example of the diagrammatics for these monoids.
\end{Notation}

\begin{enumerate}[label=$\bullet$]

\item The \emph{partition monoid} $\pamon$ of all diagrams of partitions of a $2n$-element set. The \emph{planar partition monoid} $\ppamon$ is the 
respective planar submonoid
of $\pamon$.
\begin{gather*}
\begin{tikzpicture}[anchorbase]
\draw[usual] (0.5,0) to[out=90,in=180] (1.25,0.45) to[out=0,in=90] (2,0);
\draw[usual] (0.5,0) to[out=90,in=180] (1,0.35) to[out=0,in=90] (1.5,0);
\draw[usual] (0,1) to[out=270,in=180] (0.75,0.55) to[out=0,in=270] (1.5,1);
\draw[usual] (1.5,1) to[out=270,in=180] (2,0.55) to[out=0,in=270] (2.5,1);
\draw[usual] (0,0) to (0.5,1);
\draw[usual] (1,0) to (1,1);
\draw[usual] (2.5,0) to (2.5,1);
\draw[usual,dot] (2,1) to (2,0.8);
\end{tikzpicture}
\in\pamon
,\quad
\begin{tikzpicture}[anchorbase]
\draw[usual] (0.5,0) to[out=90,in=180] (1.25,0.45) to[out=0,in=90] (2,0);
\draw[usual] (0.5,0) to[out=90,in=180] (1,0.35) to[out=0,in=90] (1.5,0);
\draw[usual] (0.5,1) to[out=270,in=180] (1,0.55) to[out=0,in=270] (1.5,1);
\draw[usual] (1.5,1) to[out=270,in=180] (2,0.55) to[out=0,in=270] (2.5,1);
\draw[usual] (0,0) to (0,1);
\draw[usual] (2.5,0) to (2.5,1);
\draw[usual,dot] (1,0) to (1,0.2);
\draw[usual,dot] (1,1) to (1,0.8);
\draw[usual,dot] (2,1) to (2,0.8);
\end{tikzpicture}
\in\ppamon
.
\end{gather*}

\item The \emph{rook-Brauer monoid} $\robrmon$ consisting of all diagrams with components of size $1$ or $2$. The planar rook-Brauer monoid 
$\probrmon=\momon$ is also called \emph{Motzkin monoid}.
\begin{gather*}
\begin{tikzpicture}[anchorbase]
\draw[usual] (1,0) to[out=90,in=180] (1.25,0.25) to[out=0,in=90] (1.5,0);
\draw[usual] (1,1) to[out=270,in=180] (1.75,0.55) to[out=0,in=270] (2.5,1);
\draw[usual] (0,0) to (0.5,1);
\draw[usual] (2.5,0) to (2,1);
\draw[usual,dot] (0.5,0) to (0.5,0.2);
\draw[usual,dot] (2,0) to (2,0.2);
\draw[usual,dot] (0,1) to (0,0.8);
\draw[usual,dot] (1.5,1) to (1.5,0.8);
\end{tikzpicture}
\in\robrmon
,\quad
\begin{tikzpicture}[anchorbase]
\draw[usual] (0.5,0) to[out=90,in=180] (1.25,0.5) to[out=0,in=90] (2,0);
\draw[usual] (1,0) to[out=90,in=180] (1.25,0.25) to[out=0,in=90] (1.5,0);
\draw[usual] (2,1) to[out=270,in=180] (2.25,0.75) to[out=0,in=270] (2.5,1);
\draw[usual] (0,0) to (1,1);
\draw[usual,dot] (2.5,0) to (2.5,0.2);
\draw[usual,dot] (0,1) to (0,0.8);
\draw[usual,dot] (0.5,1) to (0.5,0.8);
\draw[usual,dot] (1.5,1) to (1.5,0.8);
\end{tikzpicture}
\in\momon
.
\end{gather*}

\item The \emph{Brauer monoid} $\brmon$ consisting of all diagrams with components of size $2$. The planar Brauer monoid $\pbrmon=\tlmon$ is known as the \emph{Temperley--Lieb monoid} (sometimes $\tlmon$ is called \emph{Jones monoid} or \emph{Kauffman monoid}).
\begin{gather*}
\begin{tikzpicture}[anchorbase]
\draw[usual] (0.5,0) to[out=90,in=180] (1.25,0.45) to[out=0,in=90] (2,0);
\draw[usual] (1,0) to[out=90,in=180] (1.25,0.25) to[out=0,in=90] (1.5,0);
\draw[usual] (0,1) to[out=270,in=180] (0.75,0.55) to[out=0,in=270] (1.5,1);
\draw[usual] (1,1) to[out=270,in=180] (1.75,0.55) to[out=0,in=270] (2.5,1);
\draw[usual] (0,0) to (0.5,1);
\draw[usual] (2.5,0) to (2,1);
\end{tikzpicture}
\in\brmon
,\quad
\begin{tikzpicture}[anchorbase]
\draw[usual] (0.5,0) to[out=90,in=180] (1.25,0.5) to[out=0,in=90] (2,0);
\draw[usual] (1,0) to[out=90,in=180] (1.25,0.25) to[out=0,in=90] (1.5,0);
\draw[usual] (0,1) to[out=270,in=180] (0.25,0.75) to[out=0,in=270] (0.5,1);
\draw[usual] (2,1) to[out=270,in=180] (2.25,0.75) to[out=0,in=270] (2.5,1);
\draw[usual] (0,0) to (1,1);
\draw[usual] (2.5,0) to (1.5,1);
\end{tikzpicture}
\in\tlmon
.
\end{gather*}

\item The \emph{rook monoid} or 
\emph{symmetric inverse semigroup} $\romon$ consisting of all diagrams with components of size $1$ or $2$, and all partitions have at most one component 
at the bottom and at most one at the top. The \emph{planar rook monoid}
$\promon$ is the corresponding submonoid.
\begin{gather*}
\begin{tikzpicture}[anchorbase]
\draw[usual] (0,0) to (1,1);
\draw[usual] (0.5,0) to (0,1);
\draw[usual] (2,0) to (2,1);
\draw[usual] (2.5,0) to (0.5,1);
\draw[usual,dot] (1,0) to (1,0.2);
\draw[usual,dot] (1.5,0) to (1.5,0.2);
\draw[usual,dot] (1.5,1) to (1.5,0.8);
\draw[usual,dot] (2.5,1) to (2.5,0.8);
\end{tikzpicture}
\in\romon
,\quad
\begin{tikzpicture}[anchorbase]
\draw[usual] (0,0) to (0.5,1);
\draw[usual] (0.5,0) to (1,1);
\draw[usual] (2,0) to (1.5,1);
\draw[usual] (2.5,0) to (2.5,1);
\draw[usual,dot] (1,0) to (1,0.2);
\draw[usual,dot] (1.5,0) to (1.5,0.2);
\draw[usual,dot] (0,1) to (0,0.8);
\draw[usual,dot] (2,1) to (2,0.8);
\end{tikzpicture}
\in\promon
.
\end{gather*}

\item The \emph{symmetric group} $\sym$ consisting of all matchings with components of size $1$. The \emph{planar symmetric group} is trivial $\psym\cong\onemon$.
\begin{gather*}
\begin{tikzpicture}[anchorbase]
\draw[usual] (0,0) to (1,1);
\draw[usual] (0.5,0) to (0,1);
\draw[usual] (1,0) to (1.5,1);
\draw[usual] (1.5,0) to (2.5,1);
\draw[usual] (2,0) to (2,1);
\draw[usual] (2.5,0) to (0.5,1);
\end{tikzpicture}
\in\sym
,\quad
\begin{tikzpicture}[anchorbase]
\draw[usual] (0,0) to (0,1);
\draw[usual] (0.5,0) to (0.5,1);
\draw[usual] (1,0) to (1,1);
\draw[usual] (1.5,0) to (1.5,1);
\draw[usual] (2,0) to (2,1);
\draw[usual] (2.5,0) to (2.5,1);
\end{tikzpicture}
\in\psym
.
\end{gather*}
The monoids $\sym$ and $\psym$ are, of course, the same groups as in \autoref{SS:DAlgebrasTMon}.
$\sym$ is isomorphic to a subgroup of any of the symmetric monoids, and 
$\psym$ is isomorphic to a subgroup of any of the planar monoids, both 
in the evident way. We will use this below.

\end{enumerate}

\begin{Notation}\label{N:DAlgebrasTMonTCellsTwo}
Below we write $\xmon$ for any of the monoids listed above.
\end{Notation}

\begin{Notation}\label{N:DAlgebraNamesTwo}
Additionally to \autoref{N:DAlgebraNames} we now also have
\begin{gather*}
\text{caps}\colon
\begin{tikzpicture}[anchorbase,scale=0.55]
\draw[white] (0,0) to[out=90,in=270] (1,1);
\draw[usual] (0,0) to[out=90,in=180] (0.5,0.5) to[out=0,in=90] (1,0);
\end{tikzpicture}
\,,\quad
\text{cups}\colon
\begin{tikzpicture}[anchorbase,scale=0.55]
\draw[white] (0,0) to[out=90,in=270] (0,1);
\draw[usual] (0,1) to[out=270,in=180] (0.5,0.5) to[out=0,in=270] (1,1);
\end{tikzpicture}
\,,\quad
\text{splits}\colon
\begin{tikzpicture}[anchorbase,scale=0.55]
\draw[usual] (0,0) to[out=90,in=270] (0,1);
\draw[usual] (0,0) to[out=90,in=270] (1,1);
\end{tikzpicture}
\,,
\begin{tikzpicture}[anchorbase,scale=0.55]
\draw[usual] (0,0) to[out=90,in=270] (0,1);
\draw[usual] (0,0) to[out=90,in=270] (1,1);
\draw[usual] (0,0) to[out=90,in=270] (2,1);
\end{tikzpicture}
\,,\dots,\quad
\text{bottom dots}\colon
\begin{tikzpicture}[anchorbase,scale=0.55]
\draw[white] (0,0) to[out=90,in=270] (0,1);
\draw[usual,dot] (0,0) to (0,0.3);
\end{tikzpicture}
\,,
\end{gather*}
which we, as indicated, give the names \emph{caps}, \emph{cups}, 
\emph{splits} and \emph{bottom dots}.
\end{Notation}

\begin{Definition}\label{D:DAlgebrasTheAlgebras}
Fix $\para\in\K$.
For $\xmon$ we define its associated 
(symmetric or planar) 
\emph{diagram algebra} $\xmon(\para)$ to be the algebra with 
basis $\xmon$ and multiplication of basis elements given by the monoid 
multiplication except that all closed components are evaluated to $\para$.
\end{Definition}

Note that $\xmon(\para)$ is $\K$-linear, but we decided to use 
$\xmon(\para)$ as the notation instead of {\eg} 
$\K[\xmon](\para)$.

\begin{Example}\label{E:DAlgebrasTheAlgebras}
Typical examples how the multiplication in diagram algebras works is
\begin{gather*}
\begin{tikzpicture}[anchorbase]
\draw[usual] (0.5,0) to[out=90,in=180] (1.25,0.5) to[out=0,in=90] (2,0);
\draw[usual] (1,0) to[out=90,in=180] (1.25,0.25) to[out=0,in=90] (1.5,0);
\draw[usual] (2,1) to[out=270,in=180] (2.25,0.75) to[out=0,in=270] (2.5,1);
\draw[usual] (0,0) to (1,1);
\draw[usual,dot] (2.5,0) to (2.5,0.2);
\draw[usual,dot] (0,1) to (0,0.8);
\draw[usual,dot] (0.5,1) to (0.5,0.8);
\draw[usual,dot] (1.5,1) to (1.5,0.8);
\draw[very thick,densely dotted,tomato] (0,1) to (2.5,1);
\draw[usual] (0.5,1) to[out=90,in=180] (1.25,1.5) to[out=0,in=90] (2,1);
\draw[usual] (1,1) to[out=90,in=180] (1.25,1.25) to[out=0,in=90] (1.5,1);
\draw[usual] (2,2) to[out=270,in=180] (2.25,1.75) to[out=0,in=270] (2.5,2);
\draw[usual] (0,1) to (1,2);
\draw[usual,dot] (2.5,1) to (2.5,1.2);
\draw[usual,dot] (0,2) to (0,1.8);
\draw[usual,dot] (0.5,2) to (0.5,1.8);
\draw[usual,dot] (1.5,2) to (1.5,1.8);
\end{tikzpicture}
=\para\cdot
\begin{tikzpicture}[anchorbase]
\draw[usual] (0.5,0) to[out=90,in=180] (1.25,0.5) to[out=0,in=90] (2,0);
\draw[usual] (1,0) to[out=90,in=180] (1.25,0.25) to[out=0,in=90] (1.5,0);
\draw[usual] (2,1) to[out=270,in=180] (2.25,0.75) to[out=0,in=270] (2.5,1);
\draw[usual,dot] (0,0) to (0,0.2);
\draw[usual,dot] (1,1) to (1,0.8);
\draw[usual,dot] (2.5,0) to (2.5,0.2);
\draw[usual,dot] (0,1) to (0,0.8);
\draw[usual,dot] (0.5,1) to (0.5,0.8);
\draw[usual,dot] (1.5,1) to (1.5,0.8);
\end{tikzpicture}
,\quad
\begin{tikzpicture}[anchorbase]
\draw[usual] (0.5,1) to[out=270,in=180] (1.25,0.5) to[out=0,in=270] (2,1);
\draw[usual] (1,1) to[out=270,in=180] (1.25,0.75) to[out=0,in=270] (1.5,1);
\draw[usual] (2,0) to[out=90,in=180] (2.25,0.25) to[out=0,in=90] (2.5,0);
\draw[usual] (0,1) to (1,0);
\draw[usual,dot] (2.5,1) to (2.5,0.8);
\draw[usual,dot] (0,0) to (0,0.2);
\draw[usual,dot] (0.5,0) to (0.5,0.2);
\draw[usual,dot] (1.5,0) to (1.5,0.2);
\draw[very thick,densely dotted,tomato] (0,1) to (2.5,1);
\draw[usual] (0.5,1) to[out=90,in=180] (1.25,1.5) to[out=0,in=90] (2,1);
\draw[usual] (1,1) to[out=90,in=180] (1.25,1.25) to[out=0,in=90] (1.5,1);
\draw[usual] (2,2) to[out=270,in=180] (2.25,1.75) to[out=0,in=270] (2.5,2);
\draw[usual] (0,1) to (1,2);
\draw[usual,dot] (2.5,1) to (2.5,1.2);
\draw[usual,dot] (0,2) to (0,1.8);
\draw[usual,dot] (0.5,2) to (0.5,1.8);
\draw[usual,dot] (1.5,2) to (1.5,1.8);
\end{tikzpicture}
=\para^{3}\cdot
\begin{tikzpicture}[anchorbase]
\draw[usual] (2,0) to[out=90,in=180] (2.25,0.25) to[out=0,in=90] (2.5,0);
\draw[usual] (1,0) to (1,1);
\draw[usual,dot] (0,0) to (0,0.2);
\draw[usual,dot] (0.5,0) to (0.5,0.2);
\draw[usual,dot] (1.5,0) to (1.5,0.2);
\draw[usual] (2,1) to[out=270,in=180] (2.25,0.75) to[out=0,in=270] (2.5,1);
\draw[usual,dot] (0,1) to (0,0.8);
\draw[usual,dot] (0.5,1) to (0.5,0.8);
\draw[usual,dot] (1.5,1) to (1.5,0.8);
\end{tikzpicture}
,
\end{gather*}
which illustrate the multiplication of $\momon[6](\para)$.
\end{Example}

\begin{Remark}\label{R:DAlgebrasTheAlgebras}
For some of the $\xmon$ one can define 
associated multiparameter 
diagram algebras. For example, for $\momon$ one could evaluate 
circles to $\para_{1}$ and intervals to $\para_{2}$. Our discussion 
below works {\muta} for these multiparameter diagram algebras as well.
\end{Remark}

A main difference between $\xmon$ and $\tmon$, $\ptmon$ is that 
$\xmon$ is involutive using 
the \emph{diagrammatic antiinvolution} ${\placeholder}^{\star}$, {\eg}:
\begin{gather*}
\left(\begin{tikzpicture}[anchorbase]
\draw[usual] (0.5,0) to[out=90,in=180] (1.25,0.45) to[out=0,in=90] (2,0);
\draw[usual] (0.5,0) to[out=90,in=180] (1,0.35) to[out=0,in=90] (1.5,0);
\draw[usual] (0,1) to[out=270,in=180] (0.75,0.55) to[out=0,in=270] (1.5,1);
\draw[usual] (1.5,1) to[out=270,in=180] (2,0.55) to[out=0,in=270] (2.5,1);
\draw[usual] (0,0) to (0.5,1);
\draw[usual] (1,0) to (1,1);
\draw[usual] (2.5,0) to (2.5,1);
\draw[usual,dot] (2,1) to (2,0.8);
\end{tikzpicture}\right)^{\star}
=
\begin{tikzpicture}[anchorbase,yscale=-1]
\draw[usual] (0.5,0) to[out=90,in=180] (1.25,0.45) to[out=0,in=90] (2,0);
\draw[usual] (0.5,0) to[out=90,in=180] (1,0.35) to[out=0,in=90] (1.5,0);
\draw[usual] (0,1) to[out=270,in=180] (0.75,0.55) to[out=0,in=270] (1.5,1);
\draw[usual] (1.5,1) to[out=270,in=180] (2,0.55) to[out=0,in=270] (2.5,1);
\draw[usual] (0,0) to (0.5,1);
\draw[usual] (1,0) to (1,1);
\draw[usual] (2.5,0) to (2.5,1);
\draw[usual,dot] (2,1) to (2,0.8);
\end{tikzpicture}
\,.
\end{gather*}
The $\K$-linear extensions of the diagrammatic antiinvolution 
endows $\xmon(\para)$ 
with the structure of an involutive sandwich cellular algebra, as 
we will see in \autoref{P:DAlgebrasDiagram} below. 
(We will always use the diagrammatic antiinvolution in this section.)
Hence, by \autoref{L:SandwichInvolutive} the $J$-cells are squares 
and we can focus on describing left cells and right cells 
come for free.

\begin{Example}\label{E:DAlgebrasCells}
Let us give two examples how the cells structure of $\tlmon(\para)$ 
looks like.
\begin{gather}\label{Eq:TLCells}
\begin{gathered}
\xy
(0,0)*{\begin{gathered}
\begin{tabular}{C|C}
\arrayrulecolor{tomato}
\cellcolor{mydarkblue!25}
\begin{tikzpicture}[anchorbase]
\draw[usual] (0,0) to[out=90,in=180] (0.25,0.2) to[out=0,in=90] (0.5,0);
\draw[usual] (0,0.5) to[out=270,in=180] (0.25,0.3) to[out=0,in=270] (0.5,0.5);
\draw[usual] (1,0) to (1,0.5);
\end{tikzpicture} &
\cellcolor{mydarkblue!25}
\begin{tikzpicture}[anchorbase]
\draw[usual] (0.5,0) to[out=90,in=180] (0.75,0.2) to[out=0,in=90] (1,0);
\draw[usual] (0,0.5) to[out=270,in=180] (0.25,0.3) to[out=0,in=270] (0.5,0.5);
\draw[usual] (0,0) to (1,0.5);
\end{tikzpicture}
\\
\hline
\cellcolor{mydarkblue!25}
\begin{tikzpicture}[anchorbase,xscale=-1]
\draw[usual] (0.5,0) to[out=90,in=180] (0.75,0.2) to[out=0,in=90] (1,0);
\draw[usual] (0,0.5) to[out=270,in=180] (0.25,0.3) to[out=0,in=270] (0.5,0.5);
\draw[usual] (0,0) to (1,0.5);
\end{tikzpicture} &
\cellcolor{mydarkblue!25}
\begin{tikzpicture}[anchorbase,xscale=-1]
\draw[usual] (0,0) to[out=90,in=180] (0.25,0.2) to[out=0,in=90] (0.5,0);
\draw[usual] (0,0.5) to[out=270,in=180] (0.25,0.3) to[out=0,in=270] (0.5,0.5);
\draw[usual] (1,0) to (1,0.5);
\end{tikzpicture}
\end{tabular}
\\[3pt]
\begin{tabular}{C}
\arrayrulecolor{tomato}
\cellcolor{mydarkblue!25}
\begin{tikzpicture}[anchorbase]
\draw[usual] (0,0) to (0,0.5);
\draw[usual] (0.5,0) to (0.5,0.5);
\draw[usual] (1,0) to (1,0.5);
\end{tikzpicture}
\end{tabular}
\end{gathered}};
(-45,4)*{\jcell_{1}};
(-45,-7)*{\jcell_{3}};
(45,4)*{\sand[1]\cong\onemon};
(45,-7)*{\sand[3]\cong\onemon};
(-51,0)*{\phantom{a}};
\endxy
\quad,
\\[0.27cm]
\hline
\\[-0.27cm]
\xy
(0,0)*{\begin{gathered}
\begin{tabular}{C|C}
\arrayrulecolor{tomato}
\cellcolor{mydarkblue!25}
\begin{tikzpicture}[anchorbase]
\draw[usual] (0,0) to[out=90,in=180] (0.25,0.2) to[out=0,in=90] (0.5,0);
\draw[usual] (0,0.5) to[out=270,in=180] (0.25,0.3) to[out=0,in=270] (0.5,0.5);
\draw[usual] (1,0) to[out=90,in=180] (1.25,0.2) to[out=0,in=90] (1.5,0);
\draw[usual] (1,0.5) to[out=270,in=180] (1.25,0.3) to[out=0,in=270] (1.5,0.5);
\end{tikzpicture} &
\cellcolor{mydarkblue!25}
\begin{tikzpicture}[anchorbase]
\draw[usual] (0,0) to[out=45,in=180] (0.75,0.20) to[out=0,in=135] (1.5,0);
\draw[usual] (0.5,0) to[out=90,in=180] (0.75,0.1) to[out=0,in=90] (1,0);
\draw[usual] (0,0.5) to[out=270,in=180] (0.25,0.3) to[out=0,in=270] (0.5,0.5);
\draw[usual] (1,0.5) to[out=270,in=180] (1.25,0.3) to[out=0,in=270] (1.5,0.5);
\end{tikzpicture}
\\
\hline
\cellcolor{mydarkblue!25}
\begin{tikzpicture}[anchorbase]
\draw[usual] (0,0) to[out=90,in=180] (0.25,0.2) to[out=0,in=90] (0.5,0);
\draw[usual] (1,0) to[out=90,in=180] (1.25,0.2) to[out=0,in=90] (1.5,0);
\draw[usual] (0,0.5) to[out=315,in=180] (0.75,0.3) to[out=0,in=225] (1.5,0.5);
\draw[usual] (0.5,0.5) to[out=270,in=180] (0.75,0.4) to[out=0,in=270] (1,0.5);
\end{tikzpicture} &
\cellcolor{mydarkblue!25}
\begin{tikzpicture}[anchorbase]
\draw[usual] (0,0) to[out=45,in=180] (0.75,0.20) to[out=0,in=135] (1.5,0);
\draw[usual] (0,0.5) to[out=315,in=180] (0.75,0.3) to[out=0,in=225] (1.5,0.5);
\draw[usual] (0.5,0) to[out=90,in=180] (0.75,0.1) to[out=0,in=90] (1,0);
\draw[usual] (0.5,0.5) to[out=270,in=180] (0.75,0.4) to[out=0,in=270] (1,0.5);
\end{tikzpicture}
\end{tabular}
\\[3pt]
\begin{tabular}{C|C|C}
\arrayrulecolor{tomato}
\cellcolor{mydarkblue!25}
\begin{tikzpicture}[anchorbase]
\draw[usual] (0,0) to[out=90,in=180] (0.25,0.2) to[out=0,in=90] (0.5,0);
\draw[usual] (0,0.5) to[out=270,in=180] (0.25,0.3) to[out=0,in=270] (0.5,0.5);
\draw[usual] (1,0) to (1,0.5);
\draw[usual] (1.5,0) to (1.5,0.5);
\end{tikzpicture} & 
\cellcolor{mydarkblue!25}
\begin{tikzpicture}[anchorbase]
\draw[usual] (0.5,0) to[out=90,in=180] (0.75,0.2) to[out=0,in=90] (1,0);
\draw[usual] (0,0.5) to[out=270,in=180] (0.25,0.3) to[out=0,in=270] (0.5,0.5);
\draw[usual] (0,0) to (1,0.5);
\draw[usual] (1.5,0) to (1.5,0.5);
\end{tikzpicture} &
\begin{tikzpicture}[anchorbase]
\draw[usual] (1,0) to[out=90,in=180] (1.25,0.2) to[out=0,in=90] (1.5,0);
\draw[usual] (0,0.5) to[out=270,in=180] (0.25,0.3) to[out=0,in=270] (0.5,0.5);
\draw[usual] (0,0) to (1,0.5);
\draw[usual] (0.5,0) to (1.5,0.5);
\end{tikzpicture}
\\
\hline
\cellcolor{mydarkblue!25}
\begin{tikzpicture}[anchorbase]
\draw[usual] (0,0) to[out=90,in=180] (0.25,0.2) to[out=0,in=90] (0.5,0);
\draw[usual] (0.5,0.5) to[out=270,in=180] (0.75,0.3) to[out=0,in=270] (1,0.5);
\draw[usual] (1,0) to (0,0.5);
\draw[usual] (1.5,0) to (1.5,0.5);
\end{tikzpicture} & 
\cellcolor{mydarkblue!25}
\begin{tikzpicture}[anchorbase]
\draw[usual] (0,0) to (0,0.5);
\draw[usual] (0.5,0) to[out=90,in=180] (0.75,0.2) to[out=0,in=90] (1,0);
\draw[usual] (0.5,0.5) to[out=270,in=180] (0.75,0.3) to[out=0,in=270] (1,0.5);
\draw[usual] (1.5,0) to (1.5,0.5);
\end{tikzpicture} &
\cellcolor{mydarkblue!25}
\begin{tikzpicture}[anchorbase]
\draw[usual] (0,0) to (0,0.5);
\draw[usual] (0.5,0) to (1.5,0.5);
\draw[usual] (1,0) to[out=90,in=180] (1.25,0.2) to[out=0,in=90] (1.5,0);
\draw[usual] (0.5,0.5) to[out=270,in=180] (0.75,0.3) to[out=0,in=270] (1,0.5);
\end{tikzpicture}
\\
\hline
\begin{tikzpicture}[anchorbase]
\draw[usual] (0,0) to[out=90,in=180] (0.25,0.2) to[out=0,in=90] (0.5,0);
\draw[usual] (1,0.5) to[out=270,in=180] (1.25,0.3) to[out=0,in=270] (1.5,0.5);
\draw[usual] (1,0) to (0,0.5);
\draw[usual] (1.5,0) to (0.5,0.5);
\end{tikzpicture} & 
\cellcolor{mydarkblue!25}
\begin{tikzpicture}[anchorbase]
\draw[usual] (0,0) to (0,0.5);
\draw[usual] (1.5,0) to (0.5,0.5);
\draw[usual] (0.5,0) to[out=90,in=180] (0.75,0.2) to[out=0,in=90] (1,0);
\draw[usual] (1,0.5) to[out=270,in=180] (1.25,0.3) to[out=0,in=270] (1.5,0.5);
\end{tikzpicture} &
\cellcolor{mydarkblue!25}
\begin{tikzpicture}[anchorbase]
\draw[usual] (0,0) to (0,0.5);
\draw[usual] (0.5,0) to (0.5,0.5);
\draw[usual] (1,0) to[out=90,in=180] (1.25,0.2) to[out=0,in=90] (1.5,0);
\draw[usual] (1,0.5) to[out=270,in=180] (1.25,0.3) to[out=0,in=270] (1.5,0.5);
\end{tikzpicture}
\end{tabular}
\\[3pt]
\begin{tabular}{C}
\arrayrulecolor{tomato}
\cellcolor{mydarkblue!25}
\begin{tikzpicture}[anchorbase]
\draw[usual] (0,0) to (0,0.5);
\draw[usual] (0.5,0) to (0.5,0.5);
\draw[usual] (1,0) to (1,0.5);
\draw[usual] (1.5,0) to (1.5,0.5);
\end{tikzpicture}
\end{tabular}
\end{gathered}};
(-45,13)*{\jcell_{0}};
(-45,-2.5)*{\jcell_{2}};
(-45,-16.5)*{\jcell_{4}};
(45,13)*{\sand[0]\cong\onemon};
(45,-2.5)*{\sand[2]\cong\onemon};
(45,-16.5)*{\sand[4]\cong\onemon};
(-51,0)*{\phantom{a}};
\endxy
\quad.
\end{gathered}
\end{gather}
These are the cells of 
$\tlmon[3](\para)$ and $\tlmon[4](\para)$ for invertible $\para$. If $\para$ is not invertible, then the picture is the same but 
with uncolored diagonal $H$-cells for 
$\jcell_{1}$ and $\jcell_{2}$, and $\jcell_{0}$ having no 
idempotent at all.
\end{Example}

The following analog of \autoref{L:DAlgebrasTMonFactor} is easy to verify.

\begin{Lemma}\label{L:DAlgebrasDMonFactor}
Let $\xmon$ be symmetric.
For $a\in\xmon$ there is a unique factorization of the form $a=\tau\circ\sigma_{\lambda}\circ\beta$ such that $\beta$ 
and $\tau$ have a minimal 
number of crossings, $\beta$ contains no cups, splits or top dots, 
$\tau$ contains no caps, merges or bottom dots
and $\sigma_{\lambda}\in\sym[\lambda]$ for minimal $\lambda$.

Similarly for $a\in\xmon$ when $\xmon$ is planar, but with $\sigma_{\lambda}\in\psym=\onemon$.\qed
\end{Lemma}

As before, we get the notions of through strands {\etc}

\begin{Proposition}\label{P:DAlgebrasDiagram}
We have the following for the pair $\big(\xmon(\para),\xmon\big)$.
\begin{enumerate}

\item The $J$-cells of $\xmon(\para)$
are given by diagrams with a fixed number of through strands $\lambda$.
The $\leq_{lr}$-order is a total order and increases as the number of through strands decreases. See \autoref{Eq:DAlgebrasTable} for a summary.

\item The left cells of $\xmon(\para)$ 
are given by diagrams where one fixes 
the bottom of the diagram, and similarly 
right cells are given by diagrams where one fixes 
the top of the diagram.
The $\leq_{l}$ and the $\leq_{r}$-order increases as the number of through strands decreases. For $\#\lcell_{\lambda}$ see \autoref{Eq:DAlgebrasTable}.

\item The idempotency of the $J$-cells are as follows.
Assume $\para$ is invertible in $\K$.
Then all $J$-cells 
of $\xmon(\para)$ are strictly idempotent.
Second, assume $\para$ is not invertible in $\K$. Then the following 
$J$-cells for $\xmon(\para)$ are strictly idempotent, while all other $J$-cells are not idempotent (we list the various $\xmon$):

\begin{enumerate}

\item All $J$-cells except the top for $\pamon$, $\ppamon$, $\robrmon$ and $\momon$, and $\brmon$, $\tlmon$ for even $n$.

\item Only the bottom $J$-cell for $\romon$ and 
$\promon$.

\item All $J$-cells for $\sym$ and $\psym$, and $\brmon$, $\tlmon$ for $n$ odd.

\end{enumerate}
See also \autoref{Eq:DAlgebrasTable}.
The sandwiched algebras are given in \autoref{Eq:DAlgebrasTable}.

\item The pair $\big(\xmon(\para),\xmon\big)$ is an involutive sandwich pair,
that for $\xmon$ symmetric comes neither from a cellular nor an affine 
cellular algebra, but can be refined into a cellular pair.

\end{enumerate}
\end{Proposition}

The following table summarizes the cell structure of these diagram monoids where $\lambda\in\Pcal$ (the column for $\sand$ only applies if the parent $J$-cell is strictly idempotent):
\begin{gather}\label{Eq:DAlgebrasTable}
\scalebox{0.9}{\begin{tabular}{c||c|c|c|c}
Monoid & $\Pcal=\apex$ for $\para^{-1}\in\K$ & $\apex$ for $\para$ not inv. & $\#\lcell_{\lambda}$ & $\sand$ \\
\hline
\hline
$\pamon$ & $\{n{<_{lr}}n{-}1{<_{lr}}\dots{<_{lr}}0\}$ & $\{n{<_{lr}}n{-}1{<_{lr}}\dots{<_{lr}}1\}$ & $\sum_{t=0}^{n}\begin{Bsmallmatrix}n\\t\end{Bsmallmatrix}\binom{t}{\lambda}$ & $\cong\sym[\lambda]$ \\
\hline
$\ppamon$ & $\{n{<_{lr}}n{-}1{<_{lr}}\dots{<_{lr}}0\}$ & $\{n{<_{lr}}n{-}1{<_{lr}}\dots{<_{lr}}1\}$ & $\tfrac{4\lambda+2}{2n+2\lambda+2}\binom{2n}{(2n{-}2\lambda)/2}$ & $\cong\onemon$ \\
\hline
$\robrmon$ & $\{n{<_{lr}}n{-}1{<_{lr}}\dots{<_{lr}}0\}$ & $\{n{<_{lr}}n{-}1{<_{lr}}\dots{<_{lr}}1\}$ & $\sum_{t=0}^{n}\binom{n}{\lambda}\binom{n-\lambda}{2t}(2t-1)!!$ &  $\cong\sym[\lambda]$ \\
\hline
$\momon$ & $\{n{<_{lr}}n{-}1{<_{lr}}\dots{<_{lr}}0\}$ & $\{n{<_{lr}}n{-}1{<_{lr}}\dots{<_{lr}}1\}$ & $\sum_{t=0}^{n}\tfrac{\lambda+1}{\lambda+t+1}\binom{n}{\lambda+2t}\binom{\lambda+2t}{t}$ & $\cong\onemon$ \\
\hline
$\brmon$ & $\{n{<_{lr}}n{-}2{<_{lr}}\dots{<_{lr}}0,1\}$ & $\{n{<_{lr}}n{-}2{<_{lr}}\dots{<_{lr}}2,1\}$ & $\binom{n}{\lambda}(n-\lambda-1)!!$ & $\cong\sym[\lambda]$ \\
\hline
$\tlmon$ & $\{n{<_{lr}}n{-}2{<_{lr}}\dots{<_{lr}}0,1\}$ & $\{n{<_{lr}}n{-}2{<_{lr}}\dots{<_{lr}}2,1\}$ & $\tfrac{2\lambda+2}{n+\lambda+2}\binom{n}{(n-\lambda)/2}$ & $\cong\onemon$ \\
\hline
$\romon$ & $\{n{<_{lr}}n{-}1{<_{lr}}\dots{<_{lr}}0\}$ & $\{n\}$ & $\binom{n}{\lambda}$ &  $\cong\sym[\lambda]$ \\
\hline
$\promon$ & $\{n{<_{lr}}n{-}1{<_{lr}}\dots{<_{lr}}0\}$ & $\{n\}$ & $\binom{n}{\lambda}$ &  $\cong\onemon$ \\
\hline
$\sym$ & $\{n\}$ & $\{n\}$ & $1$ & $\cong\sym[\lambda]$ \\
\hline
$\psym$ & $\{n\}$ & $\{n\}$ & $1$ & $\cong\onemon$ \\
\end{tabular}}
.
\end{gather}
For $\tlmon$ and 
$\brmon$ the last entry of $\Pcal$ is either 
$0$ or $1$, depending on the parity of $n$. Similarly, the 
last entry of $\apex$ for $\para$ not invertible is either $2$ or $1$, again 
depending on the parity of $n$.

\begin{Example}\label{E:DAlgebrasDiagramGap}
GAP produces the following outputs for the 
underlying monoids:
\begin{gather*}
\pamon[4]\colon
\begin{tikzpicture}[anchorbase,scale=1]
\node at (0,0) {\raisebox{0.5cm}{\scalebox{-1}[-1]{\includegraphics[height=3.5cm]{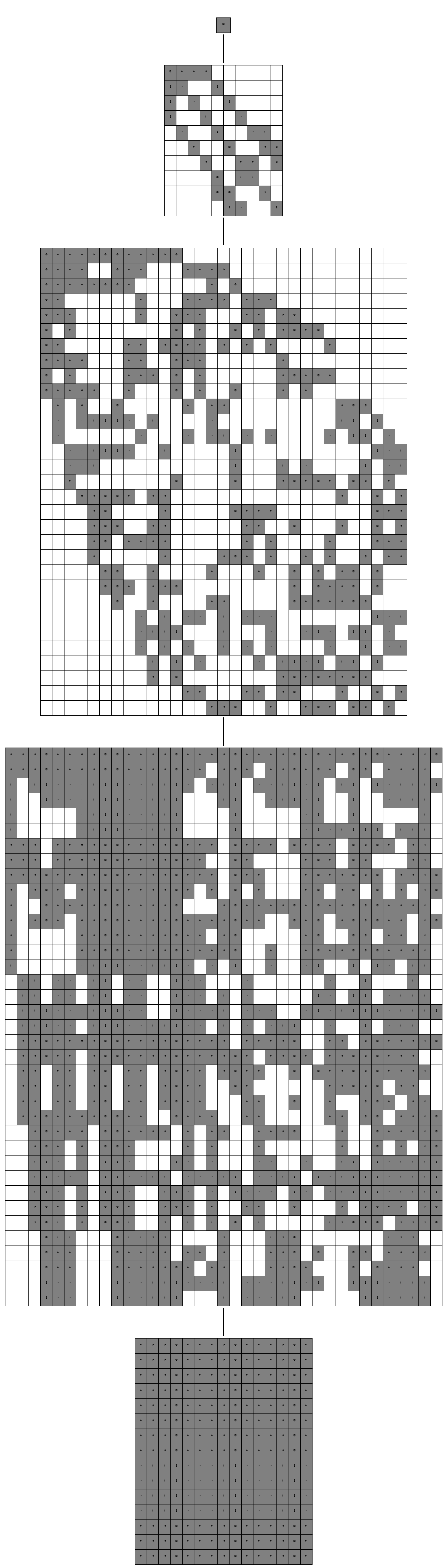}}}};
\end{tikzpicture}
,\quad
\ppamon[4]\colon
\begin{tikzpicture}[anchorbase,scale=1]
\node at (0,0) {\raisebox{0.5cm}{\scalebox{-1}[-1]{\includegraphics[height=3.5cm]{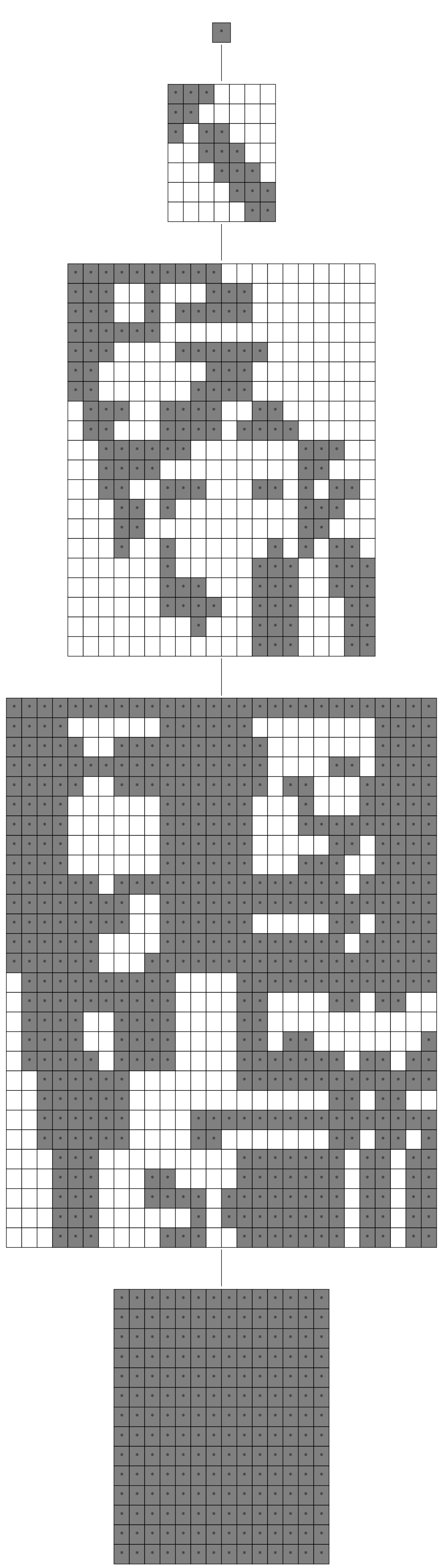}}}};
\end{tikzpicture}
,\quad
\robrmon[4]\colon
\begin{tikzpicture}[anchorbase,scale=1]
\node at (0,0) {\raisebox{0.5cm}{\scalebox{-1}[-1]{\includegraphics[height=3.5cm]{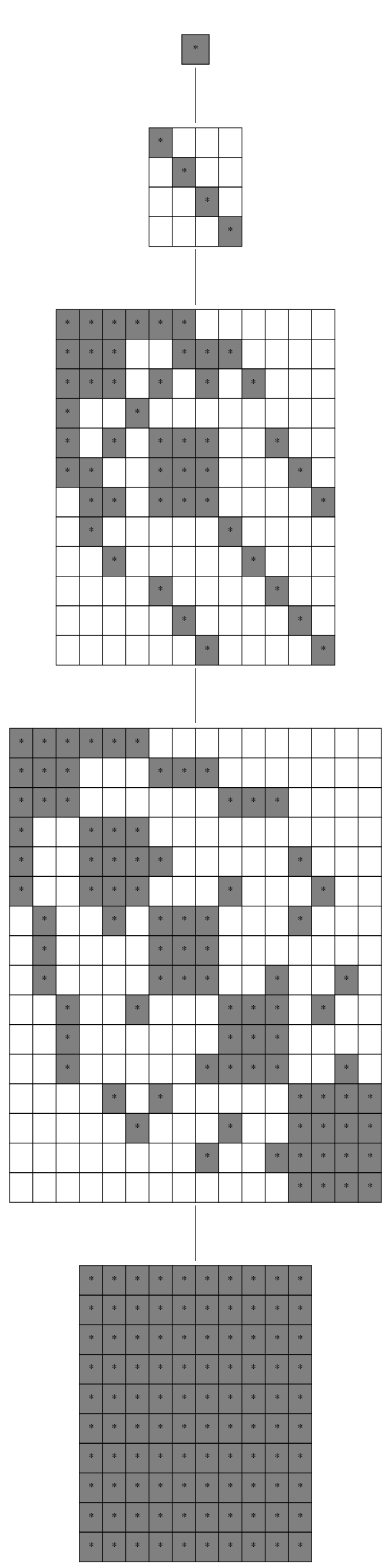}}}};
\end{tikzpicture}
,\quad
\momon[4]\colon
\begin{tikzpicture}[anchorbase,scale=1]
\node at (0,0) {\raisebox{0.5cm}{\scalebox{-1}[-1]{\includegraphics[height=3.5cm]{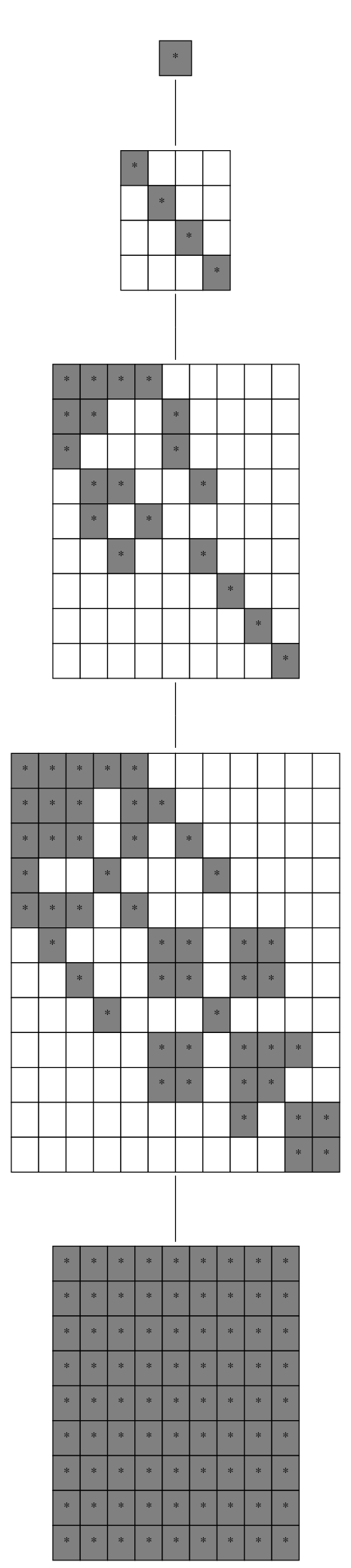}}}};
\end{tikzpicture}
,
\\
\brmon[4]\colon
\begin{tikzpicture}[anchorbase,scale=1]
\node at (0,0) {\raisebox{0.5cm}{\scalebox{-1}[-1]{\includegraphics[height=3.5cm]{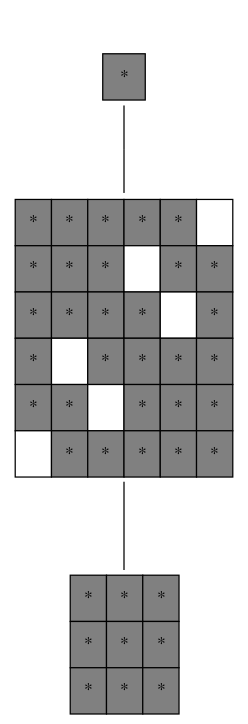}}}};
\end{tikzpicture}
,\quad
\tlmon[4]\colon
\begin{tikzpicture}[anchorbase,scale=1]
\node at (0,0) {\raisebox{0.5cm}{\scalebox{-1}[-1]{\includegraphics[height=3.5cm]{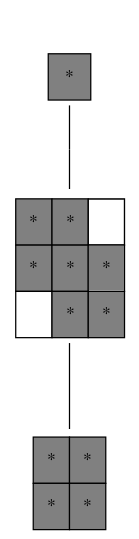}}}};
\end{tikzpicture}
,\quad
\romon[4]\colon
\begin{tikzpicture}[anchorbase,scale=1]
\node at (0,0) {\raisebox{0.5cm}{\scalebox{-1}[-1]{\includegraphics[height=3.5cm]{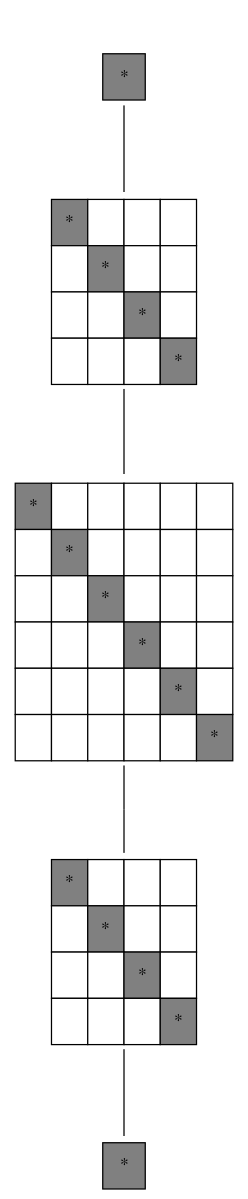}}}};
\end{tikzpicture}
,\quad
\promon[4]\colon
\begin{tikzpicture}[anchorbase,scale=1]
\node at (0,0) {\raisebox{0.5cm}{\scalebox{-1}[-1]{\includegraphics[height=3.5cm]{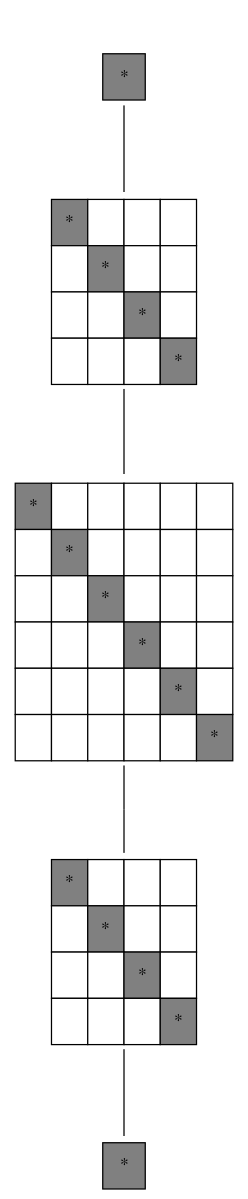}}}};
\end{tikzpicture}
.
\end{gather*}
The monoids $\sym$ and $\psym$ are not illustrated as they 
just have one cell.
\end{Example}

\begin{proof}[Proof of \autoref{P:DAlgebrasDiagram}]
\textit{(a) and (b).} 
Only the counts for $\#\lcell_{\lambda}$ are not immediate. 
Counting them is a combinatorial exercise that has 
been solved several times in the literature. 
The respective triangles, 
{\cf} \autoref{R:DAlgebrasTMonStirling}, 
are A049020, A008313, A096713,
A064189, A111062, and A007318 on OEIS.

\textit{(c).} This is easy to see and omitted. The main trick is to use
variations of
\begin{gather*}
\begin{tikzpicture}[anchorbase]
\draw[usual] (0,0.5) to[out=270,in=180] (0.25,0.3) to[out=0,in=270] (0.5,0.5);
\draw[usual] (0.5,0) to[out=90,in=270] (1,0.5);
\draw[usual] (0.5,0.5) to[out=90,in=180] (0.75,0.7) to[out=0,in=90] (1,0.5);
\draw[usual] (0,0.5) to[out=90,in=270] (0.5,1);
\end{tikzpicture}
=
\begin{tikzpicture}[anchorbase]
\draw[usual] (0,0) to (0,1);
\end{tikzpicture}
\end{gather*}
to ensure that certain cells stay strictly idempotent
even when $\para$ is not invertible.

\textit{(d).} Since $\xmon$ is an admissible monoid 
by the above, \autoref{P:SandwichMonoid} applies for $\big(\K[\xmon],\xmon\big)$. 
This result can be pulled over to $\big(\xmon(\para),\xmon\big)$. 
The second claim about cellularity follows then from 
\autoref{L:SandwichInvolutive} and
\autoref{P:SandwichRefine} (and some care with the antiinvolution).
\end{proof}

From now on we use that $\big(\xmon(\para),\xmon\big)$ is an involutive sandwich pair.

\begin{Theorem}\label{T:DAlgebrasDiaSimples}
Assume that $\K$ is a field, and consider the 
involutive sandwich pair $\big(\xmon(\para),\xmon\big)$.
\begin{enumerate}

\item Assume that $\para\neq 0$. The set of apexes for simple 
$\xmon(\para)$-modules is $\apex$
as in \autoref{Eq:DAlgebrasTable}, and
there are precisely $\big|\mathrm{P}\big(\lambda|\charkk\big)\big|$ 
(symmetric) or one (planar) simple  
$\xmon(\para)$-modules for $\lambda\in\apex$.

\item Assume that $\para=0$. Then the same statement holds 
for the restricted set of apexes as detailed in \autoref{Eq:DAlgebrasTable}.

\end{enumerate}
\end{Theorem}

\begin{proof}
This is clear by \autoref{P:DAlgebrasDiagram} 
and $H$-reduction \autoref{T:SandwichCMP}.
\end{proof}

Below we compute ranks of Gram matrices. 
To get started, here is an example:

\begin{Example}\label{E:DAlgebrasDGram}
For $\tlmon[5](\para)$ we have
\begin{gather*}
\begin{tikzpicture}[anchorbase,scale=1]
\node at (0,0) {\raisebox{0.5cm}{\scalebox{-1}[-1]{\includegraphics[height=5cm]{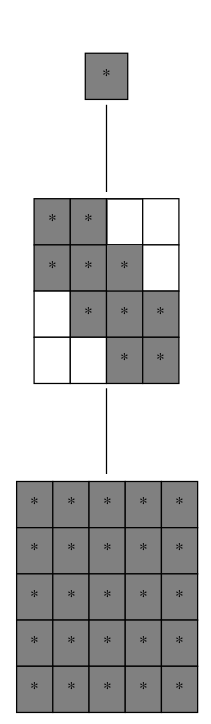}}}};
\end{tikzpicture}
,\quad
\begin{gathered}
\gmatrix[1]=\begin{psmallmatrix}
\para^{2} & \para & 1 & \para & 1 
\\
\para & \para^{2} & \para & 1 & \para 
\\
1 & \para & \para^{2} & \para & 1 
\\
\para & 1 & \para & \para^{2} & \para 
\\
1 & \para & 1 & \para & \para^{2} 
\\
\end{psmallmatrix}
,\quad
\det(\gmatrix[1])=(\para-1)^{4}(\para+1)^{4}(\para^{2}-2)
,
\\[0.8cm]
\gmatrix[3]=\begin{psmallmatrix}
\para & 1 & 0 & 0
\\
1 & \para & 1 & 0
\\
0 & 1 & \para & 1
\\
0 & 0 & 1 & \para
\end{psmallmatrix}
,\quad
\det(\gmatrix[3])=
(\para^{2}+\para-1)(\para^{2}-\para-1)
,
\\[0.8cm]
\gmatrix[5]=\begin{pmatrix}
1
\end{pmatrix}
,\quad
\det(\gmatrix[1])=1
,
\end{gathered}
\end{gather*}
as cells, Gram matrices and their determinants. The ranks 
of $\gmatrix[1]$ and $\gmatrix[3]$ depend on $\para$ 
and are easy to determine from the determinants. For example, 
if $\para=1$ or $\para=\sqrt{2}$, then $\mrk(\gmatrix[1])=1$ 
and $\mrk(\gmatrix[1])=4$, respectively.
\end{Example}

The following proposition computes ranks of Gram matrices 
for $\xmon\in\{\sym,\psym,\romon,\promon,\tlmon,\ppamon\}$ explicitly
and implicitly for $\xmon=\momon$ in $\chark=0$ as well. 
We do not know 
nice formulas for the ranks of the remaining monoids, but we give some partial results by computing the ranks of the cells close to the bottom, 
including explicit formulas for $\xmon=\momon$.

\begin{Proposition}\label{P:DAlgebrasGram}
Assume that $\K$ is a field, and consider the 
involutive sandwich pair $\big(\xmon(\para),\xmon\big)$.
The ranks of (some of the) Gram matrices are as follows.
\begin{enumerate}

\item For $\xmon\in\{\sym,\psym\}$ we have $\mrk(\gmatrix[n])=1$ for $n\in\apex$.

\item For $\xmon\in\{\romon,\promon\}$ we have $\mrk(\gmatrix)=\binom{n}{\lambda}$ 
for $\lambda\in\apex$.

\item For $\xmon=\tlmon$ we have $\mrk(\gmatrix)=tl_{n}^{\lambda}$ for $\lambda\in\apex$, with $tl_{n}^{\lambda}$ explicitly given in \autoref{R:DAlgebrasTL} below. 
For $\xmon=\ppamon$ we have $\mrk(\gmatrix)=tl_{2n}^{2\lambda}$ for $\lambda\in\apex$.

\item For $\xmon=\momon$ and $\chark=0$ the rank $\mrk(\gmatrix)$ depends on 
the multiplicity of $\para$ in an explicit multiset $\mathrm{Roots}_{n,all}$ 
explained in the proof below.

\item For $\xmon\in\{\brmon,\robrmon,\momon,\pamon\}$ 
we have $\mrk(\gmatrix[n])=1$ for $n\in\apex$.

\item For $\xmon=\brmon$, $\chark=0$ and $n-2\in\apex$ we have
\begin{gather*}
\mrk(\gmatrix[{n-2}])=
\begin{cases}
\tfrac{1}{2}n(n-1)&\text{if }\para\not\in\{2,-n+4,-2n+4\},
\\
n&\text{if }\para=2,
\\
\tfrac{1}{2}n(n-3)+1&\text{if }\para=-n+4,
\\
\tfrac{1}{2}(n+1)(n-2)&\text{if }\para=-2n+4.
\end{cases}
\end{gather*}

\item For $\xmon=\{\robrmon,\momon\}$ and $n-1\in\apex$ 
we have $\mrk(\gmatrix[{n-1}])=n-1$.

\item For $\xmon=\momon$, $\chark=0$ and $n-2\in\apex$ we have
\begin{gather*}
\mrk(\gmatrix[{n-2}])=
\begin{cases}
\frac{1}{2}(n^{2}+n-2)&\text{if }\para\not\in\{0\}\cup\mathrm{Roots}_{n-1},
\\
\frac{1}{2}(n^{2}+n-4)&\text{if }\para\in\mathrm{Roots}_{n-1}\setminus\{0\},
\\
n-1&\text{if }\para=0,n\not\equiv 0\bmod 2,
\\
n-2&\text{if }\para=0,n\equiv 0\bmod 2.
\end{cases}
\end{gather*}
The set $\mathrm{Roots}_{n-1}$ can be explicitly computed as explained in the proof below.

\item For $\xmon=\robrmon$, $\chark=0$ and $n-2\in\apex$ we have
\begin{gather*}
\mrk(\gmatrix[{n-2}])=
\begin{cases}
n(n-1)&\text{if }\para\not\in\{0,3,5-n,5-2n\},
\\
\tfrac{1}{2}n(n-1)&\text{if }\para=0,
\\
\tfrac{1}{2}n(n+1)&\text{if }\para=3,
\\
(n-1)^{2}&\text{if }\para=5-n,
\\
n(n-1)-1&\text{if }\para=5-2n.
\end{cases}
\end{gather*}

\end{enumerate}
\end{Proposition}

Note that \autoref{P:DAlgebrasGram} always assumes that 
$\lambda\in\apex$. Hence, depending on $\para$, some cases might not appear.

\begin{Remark}\label{R:DAlgebrasTL}
The number $tl_{n}^{\lambda}$ is given as follows.

Let $l\in\N$ be minimal such that $U_{l+1}(\para)=0$ where 
$U_{k}(X)$ is the 
(normalized) \emph{Chebyshev polynomial} (of the second kind) defined 
by $U_{0}(X)=1$, $U_{1}(X)=X$ and $U_{k}(X)=XU_{k-1}(X)-U_{k-2}(X)$ for $k>1$.
If no such $l\in\N$ exists we set $l=\infty$. The number $l$ 
is often called the \emph{quantum characteristic} of $(\K,\para)$.
Similarly, let further $p\in\N$ 
be minimal such that $p\cdot 1=0\in\K$, and let 
$p=\infty$ if no such $p$ exists, so $p=\charkk$, see \autoref{SS:KLClassical}.

Let $\nu_{p}$ 
denote the $p$-adic valuation. 
Let $\nu_{l,p}(x)=0$ if $x\not\equiv 0\bmod l$, 
and $\nu_{l,p}(x)=\nu_{p}(\tfrac{x}{l})$ otherwise. Let further 
$x=[\dots,x_{1},x_{0}]$ denote the $(l,p)$-adic expansion of $x$ 
given by
\begin{gather*}
[\dots,x_{1},x_{0}]=\sum_{i=1}^{\infty}lp^{i-1}x_{i}+x_{0}=x,
\quad
x_{i>0}\in\{0,\dots,p-1\},x_{0}\in\{0,\dots,l-1\}
.
\end{gather*}
Write $x\triangleleft y$ if $[\dots,x_{1},x_{0}]$ is digit-wise
smaller or equal to $[\dots,y_{1},y_{0}]$. 
We also use $x\triangleleft^{\prime}y$ if 
$x\triangleleft y$, $\nu_{l,p}(x)=\nu_{l,p}(y)$ and the $\nu_{l,p}(x)$th digit 
of $x$ and $y$ agree. Set
\begin{gather*}
e_{n,k}=
\begin{cases}
1&\text{if }n\equiv k\bmod 2,\ \nu_{l,p}(k)=\nu_{l,p}(\tfrac{n+\lambda}{2}),\ k\triangleleft^{\prime}\tfrac{n+\lambda}{2},
\\
-1&\text{if }n\equiv k\bmod 2,\ \nu_{l,p}(k)<\nu_{l,p}(\tfrac{n+\lambda}{2}),
k\triangleleft\tfrac{n+\lambda}{2}-1,
\\
0&\text{else}.
\end{cases}
\end{gather*}
Then $tl_{n}^{\lambda}=
\sum_{r=0}^{(n-\lambda)/2}e_{n-2\lambda+1,\lambda+1}\left(\tfrac{2\lambda+2}{n+\lambda+2}\binom{n}{(n-\lambda)/2}\right)$.
\end{Remark}

\begin{proof}[Proof of \autoref{P:DAlgebrasGram}]
We compute the Gram matrices using a case-by-case 
argument. (The arguments for the various cases are similar 
but differ in details.) We will also use \autoref{P:DAlgebrasDiagram} 
since we need the various numerical data computed therein.

\textit{Case $\romon$ and $\promon$.} In this case 
all Gram matrices are $\para^{n-\lambda}$ multiples 
of identity matrices. This can be seen 
as follows. We use the same order for left and right cells and then 
one calculates that the middle (where rows and columns are swapped) 
is of the form 
\begin{gather}\label{Eq:DAlgebrasRoMon}
\xy
(0,0)*{\begin{gathered}
\begin{tabular}{C||C|C|C}
\arrayrulecolor{tomato}
\rcell/\lcell & \begin{tikzpicture}[anchorbase]
\draw[usual,dot] (0,0) to (0,0.15);
\draw[usual,dot] (0.5,0) to (0.5,0.15);
\draw[usual] (1,0) to (1,0.25);
\end{tikzpicture} & \begin{tikzpicture}[anchorbase]
\draw[usual,dot] (0,0) to (0,0.15);
\draw[usual,dot] (1,0) to (1,0.15);
\draw[usual] (0.5,0) to (0.5,0.25);
\end{tikzpicture} & \begin{tikzpicture}[anchorbase]
\draw[usual,dot] (1,0) to (1,0.15);
\draw[usual,dot] (0.5,0) to (0.5,0.15);
\draw[usual] (0,0) to (0,0.25);
\end{tikzpicture}
\\
\hline
\hline
\begin{tikzpicture}[anchorbase]
\draw[usual,dot] (0,0) to (0,-0.15);
\draw[usual,dot] (0.5,0) to (0.5,-0.15);
\draw[usual] (1,0) to (1,-0.25);
\end{tikzpicture}	& \begin{tikzpicture}[anchorbase]
\draw[usual,dot] (0,0) to (0,0.15);
\draw[usual,dot] (0.5,0) to (0.5,0.15);
\draw[usual] (1,0) to (1,0.25);
\draw[usual,dot] (0,0) to (0,-0.15);
\draw[usual,dot] (0.5,0) to (0.5,-0.15);
\draw[usual] (1,0) to (1,-0.25);
\end{tikzpicture} & \begin{tikzpicture}[anchorbase]
\draw[usual,dot] (0,0) to (0,0.15);
\draw[usual,dot] (1,0) to (1,0.15);
\draw[usual] (0.5,0) to (0.5,0.25);
\draw[usual,dot] (0,0) to (0,-0.15);
\draw[usual,dot] (0.5,0) to (0.5,-0.15);
\draw[usual] (1,0) to (1,-0.25);
\end{tikzpicture} & \begin{tikzpicture}[anchorbase]
\draw[usual,dot] (1,0) to (1,0.15);
\draw[usual,dot] (0.5,0) to (0.5,0.15);
\draw[usual] (0,0) to (0,0.25);
\draw[usual,dot] (0,0) to (0,-0.15);
\draw[usual,dot] (0.5,0) to (0.5,-0.15);
\draw[usual] (1,0) to (1,-0.25);
\end{tikzpicture}
\\
\hline
\begin{tikzpicture}[anchorbase]
\draw[usual,dot] (0,0) to (0,-0.15);
\draw[usual,dot] (1,0) to (1,-0.15);
\draw[usual] (0.5,0) to (0.5,-0.25);
\end{tikzpicture}	& \begin{tikzpicture}[anchorbase]
\draw[usual,dot] (0,0) to (0,0.15);
\draw[usual,dot] (0.5,0) to (0.5,0.15);
\draw[usual] (1,0) to (1,0.25);
\draw[usual,dot] (0,0) to (0,-0.15);
\draw[usual,dot] (1,0) to (1,-0.15);
\draw[usual] (0.5,0) to (0.5,-0.25);
\end{tikzpicture} & \begin{tikzpicture}[anchorbase]
\draw[usual,dot] (0,0) to (0,0.15);
\draw[usual,dot] (1,0) to (1,0.15);
\draw[usual] (0.5,0) to (0.5,0.25);
\draw[usual,dot] (0,0) to (0,-0.15);
\draw[usual,dot] (1,0) to (1,-0.15);
\draw[usual] (0.5,0) to (0.5,-0.25);
\end{tikzpicture} & \begin{tikzpicture}[anchorbase]
\draw[usual,dot] (1,0) to (1,0.15);
\draw[usual,dot] (0.5,0) to (0.5,0.15);
\draw[usual] (0,0) to (0,0.25);
\draw[usual,dot] (0,0) to (0,-0.15);
\draw[usual,dot] (1,0) to (1,-0.15);
\draw[usual] (0.5,0) to (0.5,-0.25);
\end{tikzpicture}
\\
\hline
\begin{tikzpicture}[anchorbase]
\draw[usual,dot] (1,0) to (1,-0.15);
\draw[usual,dot] (0.5,0) to (0.5,-0.15);
\draw[usual] (0,0) to (0,-0.25);
\end{tikzpicture}	& \begin{tikzpicture}[anchorbase]
\draw[usual,dot] (0,0) to (0,0.15);
\draw[usual,dot] (0.5,0) to (0.5,0.15);
\draw[usual] (1,0) to (1,0.25);
\draw[usual,dot] (1,0) to (1,-0.15);
\draw[usual,dot] (0.5,0) to (0.5,-0.15);
\draw[usual] (0,0) to (0,-0.25);
\end{tikzpicture} & \begin{tikzpicture}[anchorbase]
\draw[usual,dot] (0,0) to (0,0.15);
\draw[usual,dot] (1,0) to (1,0.15);
\draw[usual] (0.5,0) to (0.5,0.25);
\draw[usual,dot] (1,0) to (1,-0.15);
\draw[usual,dot] (0.5,0) to (0.5,-0.15);
\draw[usual] (0,0) to (0,-0.25);
\end{tikzpicture} & \begin{tikzpicture}[anchorbase]
\draw[usual,dot] (1,0) to (1,0.15);
\draw[usual,dot] (0.5,0) to (0.5,0.15);
\draw[usual] (0,0) to (0,0.25);
\draw[usual,dot] (1,0) to (1,-0.15);
\draw[usual,dot] (0.5,0) to (0.5,-0.15);
\draw[usual] (0,0) to (0,-0.25);
\end{tikzpicture}
\end{tabular}
\end{gathered}};
\endxy
\quad
\raisebox{-0.3cm}{$\leftrightsquigarrow
\gmatrix[1]=
\begin{pmatrix}
\para^{2} & 0 & 0
\\
0 & \para^{2} & 0
\\
0 & 0 & \para^{2}
\end{pmatrix}$}
.
\end{gather}
This calculation generalizes without problems, hence, 
$\gmatrix$ is of full rank.
Note hereby, as for the rest of the proof, 
that we only care about apexes in the theorem. In particular, 
for $\para=0$ there is only one apex, the bottom $J$-cell, and 
the relevant Gram matrix is the identity.

\textit{Case $\tlmon$ and $\ppamon$.} The calculation of the Gram 
matrices for $\tlmon[\para]$ is known, but not easy, so we will not 
recall it here. See {\eg} \cite{Sp-modular-tl}
(general case using 
Temperley--Lieb combinatorics), \cite{An-simple-tl} (general case using 
tilting modules) or \cite{RiSaAu-temperley-lieb} (characteristic zero) for details. The monoid isomorphism $\ppamon\cong\tlmon[2n]$ 
from \cite[Section 1]{HaRa-partition-algebras} can then 
be used to prove the statement for the planar partition algebra from the 
Temperley--Lieb algebra case.

\textit{Case $\xmon\in\{\brmon,\momon,\robrmon,\pamon\}$, bottom cell.}
This is clear.

\textit{Case $\brmon$, $\lambda=n-2$.}
In this case there is one cap and one cup. Thus, it suffices to remember the 
endpoints of the middle. We do this by using
\begin{gather}\label{Eq:DAlgebrasBrauer}
\begin{tikzpicture}[anchorbase]
\draw[usual] (0.5,0) to (0.5,1);
\draw[usual] (1,0) to (1,1)node[above]{$\phantom{k}$};
\draw[usual] (1.5,0) to (1.5,1);
\draw[usual] (0,0)node[below]{$i$} to[out=90,in=180] (1,0.5) to[out=0,in=90] (2,0)node[below]{$j$};
\end{tikzpicture}
\leftrightsquigarrow
b[i,j]
,
1\leq i<j\leq n
,\quad
\begin{tikzpicture}[anchorbase]
\draw[usual] (0,0) to (0,1);
\draw[usual] (1,0) to (1,1);
\draw[usual] (1.5,0)node[below]{$\phantom{k}$} to (1.5,1);
\draw[usual] (0.5,1)node[above]{$k$} to[out=270,in=180] (1.25,0.5) to[out=0,in=270] (2,1)node[above]{$l$};
\end{tikzpicture}
\leftrightsquigarrow
t[k,l]
,
1\leq k<l\leq n
.
\end{gather}
The numbers in these pictures are the endpoints of the strings, read 
from left to right.
The Gram matrix becomes symbolically $\gmatrix[n-2]=(b[i,j]t[k,l])_{i,j,k,l}$ for $1\leq i<j\leq n$ and $1\leq k<l\leq n$. The entries of the Gram matrix are determined by
\begin{gather}\label{Eq:DAlgebrasSymbolic}
\begin{cases}
b[i,j]t[k,l]=\para & \text{if both endpoints match},
\\
b[i,j]t[k,l]=1 & \text{if one endpoint matches},
\\
b[i,j]t[k,l]=0 & \text{else},
\end{cases}
\quad
\scalebox{0.8}{$\begin{tikzpicture}[anchorbase]
\draw[usual] (0.5,0) to (0.5,1)node[above]{$\phantom{k}$};
\draw[usual] (1,0) to (1,1);
\draw[usual] (1.5,0) to (1.5,1);
\draw[usual] (0,0)node[below]{$1$} to[out=90,in=180] (1,0.5) to[out=0,in=90] (2,0)node[below]{$4$};
\end{tikzpicture}
\circ
\begin{tikzpicture}[anchorbase]
\draw[usual] (0,0) to (0,1);
\draw[usual] (1,0) to (1,1);
\draw[usual] (1.5,0)node[below]{$\phantom{k}$} to (1.5,1);
\draw[usual] (0.5,1)node[above]{$2$} to[out=270,in=180] (1.25,0.5) to[out=0,in=270] (2,1)node[above]{$4$};
\end{tikzpicture}
=
\begin{tikzpicture}[anchorbase]
\draw[usual] (0.5,1) to (0.5,2);
\draw[usual] (1,1) to (1,2);
\draw[usual] (1.5,1) to (1.5,2);
\draw[usual] (0,1)to[out=90,in=180] (1,1.5) to[out=0,in=90] (2,1);
\draw[usual] (0,0) to (0,1);
\draw[usual] (1,0) to (1,1);
\draw[usual] (1.5,0) to (1.5,1);
\draw[usual] (0.5,1) to[out=270,in=180] (1.25,0.5) to[out=0,in=270] (2,1);
\end{tikzpicture}$}
\rightsquigarrow
1.
\end{gather}
This can be seen as indicated above. The matrices are then easy to write down, for example
\begin{gather*}
n=4\colon\gmatrix[n-2]
=\begin{psmallmatrix}
\para & 1 & 1 & 1 & 1 & 0 \\
1 & \para & 1 & 1 & 0 & 1 \\
1 & 1 & \para & 0 & 1 & 1 \\
1 & 1 & 0 & \para & 1 & 1 \\
1 & 0 & 1 & 1 & \para & 1 \\
0 & 1 & 1 & 1 & 1 & \para \\
\end{psmallmatrix}
.
\end{gather*}
Induction verifies that the determinant is
\begin{gather*}
\det(\gmatrix[n-2])=(\para-2)^{\frac{1}{2}n(n-3)}(\para+n-4)^{n-1}(\para+2n-4).
\end{gather*}
Thus, unless $\para\in\{2,-n+4,-2n+4\}$, we get $\rk(\gmatrix[n-2])=\frac{1}{2}n(n-1)$.
For the reaming cases we get $\rk(\gmatrix[n-2])=n$ if $\para=2$, 
$\rk(\gmatrix[n-2])=\frac{1}{2}n(n-3)+1$ if $\para=-n+4$ and 
$\rk(\gmatrix[n-2])=\frac{1}{2}(n+1)(n-2)$ if $\para=-2n+4$.

\textit{Case $\robrmon$ and $\momon$, $\lambda=n-1$.}
The only way to reduce the number of through strands by one is to 
have a start and an top dot. Thus, this cell is exactly as in 
\autoref{Eq:DAlgebrasRoMon}, up to permutations, and we get the same 
Gram matrix as for $\romon$ and $\promon$, namely $\para$ times identity.

\textit{Case $\momon$, $\lambda=n-2$.}
In order to get $n-2$ through strands one either needs to have one cap 
or two bottom dots at the bottom, 
and one cup or two top dots at the top of the diagram.
Having one cap-cup pair is the Temperley--Lieb case, having two start 
and two top dots is the planar rook monoid case, and then there are the mixed cases. In the following illustration 
we again mark the endpoints of strings by their positions read from left to right. 
As for $\brmon$ above, that is \autoref{Eq:DAlgebrasBrauer}, we use 
a translation from diagrams to symbols:
\begin{gather*}
\begin{tikzpicture}[anchorbase]
\draw[white] (0,0) to (0,0.5);
\draw[usual] (0,0)node[below]{$i$} to[out=90,in=180] (0.25,0.25) to[out=0,in=90] (0.5,0)node[below]{$i{+}1$};
\end{tikzpicture}
\leftrightsquigarrow
b[i,i+1]
,1\leq i\leq n
,\;
\begin{tikzpicture}[anchorbase]
\draw[white] (0,0) to (0,0.5);
\draw[usual,dot] (0,0)node[below]{$i$} to (0,0.2);
\draw[usual,dot] (0.5,0)node[below]{$i{+}1$} to (0.5,0.2);
\end{tikzpicture}
\leftrightsquigarrow
c[i,i+1]
,1\leq i\leq n
,\;
\begin{tikzpicture}[anchorbase]
\draw[usual] (0.5,0) to (0.5,0.5);
\draw[usual,dot] (0,0)node[below]{$j$} to (0,0.2);
\draw[usual,dot] (1,0)node[below]{$k$} to (1,0.2);
\end{tikzpicture}
\leftrightsquigarrow
d[j,k]
,1\leq j<k\leq n.
\end{gather*}
We use the notation $t[i,i+1]$, $u[i,i+1]$ and $v[j,k]$ for right cells.
Additionally to \autoref{Eq:DAlgebrasSymbolic} we have:
\begin{gather*}
\begin{cases}
b[i,i+1]u[j,j+1]=\para & \text{if }i=j,
\\
b[i,i+1]u[j,j+1]=0 & \text{else},
\\
b[i,i+1]v[j,k]=0 & \text{always},
\\
u[i,i+1]u[j,j+1]=\para^{2} & \text{if }i=j,
\end{cases}
\quad
\begin{cases}
u[i,i+1]u[j,j+1]=0 & \text{else},
\\
u[i,i+1]v[j,k]=0 & \text{always},
\\
v[i,j]v[k,l]=\para^{2} & \text{if both endpoints match},
\\
v[i,j]v[k,l]=0 & \text{else}.
\end{cases}
\end{gather*}
Using these formulas and 
the corresponding symbolic Gram matrix, an analysis as for $\brmon$ 
shows that the Gram 
matrix of square size $(n-1)+(n-1)+\big(\frac{1}{2}(n-1)(n-2)\big)$ is 
\begin{gather}\label{Eq:DAlgebrasMotzkin}
\renewcommand{\arraystretch}{1.25}
\gmatrix[n-2]=
\left(
\begin{array}{c|c|c}
A_{11} & A_{12} & 0 \\
\hline
A_{12} & A_{22} & 0 \\
\hline
0 & 0 & A_{33} \\
\end{array}
\right)
,
\end{gather}
\begin{gather*}
\renewcommand{\arraystretch}{1}
A_{11}=
\begin{psmallmatrix}
\para & 1 & 0 & 0 & 0 & \dots & 0
\\
1 & \para & 1 & 0 & 0 & \dots & 0
\\
0 & 1 & \para & 1 & 0 & \dots & 0
\\
\vdots & \ddots & \ddots & \ddots & \ddots & \ddots & \vdots
\end{psmallmatrix}
,\quad
A_{12}=\para\id
,\quad
A_{22}=\para^{2}\id
,\quad
A_{33}=\para^{2}\id.
\end{gather*}
The northwest corner is the same as for $\tlmon$. An example is
\begin{gather*}
n=4\colon\gmatrix[n-2]=
\begin{psmallmatrix}
\para & 1 & 0 & \para & 0 & 0 & 0 & 0 & 0 \\
1 & \para & 1 & 0 & \para & 0 & 0 & 0 & 0 \\
0 & 1 & \para & 0 & 0 & \para & 0 & 0 & 0 \\
\para & 0 & 0 & \para^{2} & 0 & 0 & 0 & 0 & 0 \\
0 & \para & 0 & 0 & \para^{2} & 0 & 0 & 0 & 0 \\
0 & 0 & \para & 0 & 0 & \para^{2} & 0 & 0 & 0 \\
0 & 0 & 0 & 0 & 0 & 0 & \para^{2} & 0 & 0 \\
0 & 0 & 0 & 0 & 0 & 0 & 0 & \para^{2} & 0 \\
0 & 0 & 0 & 0 & 0 & 0 & 0 & 0 & \para^{2} \\
\end{psmallmatrix}
.
\end{gather*}
Let $U_{k}(X)$ again be the Chebyshev polynomial, see \autoref{R:DAlgebrasTL}. 
It is not hard to see (and well-known) that 
$\det(A_{11})=U_{n-1}(\para)$. Using this determinant formula 
and \autoref{Eq:DAlgebrasMotzkin} we get
\begin{gather*}
\det(\gmatrix[n-2])=\para^{n(n-1)}U_{n-1}(\para-1)
.
\end{gather*}
The Chebyshev 
polynomial $U_{n-1}(X-1)$ is a polynomial of degree $n-1$ 
and has distinct roots (given by explicit formulas) 
and we denote the set of these roots by $\mathrm{Roots}_{n-1}$. Thus, for $\para\not\in\{0\}\cup\mathrm{Roots}_{n-1}$, we get $\rk(\gmatrix[n-2])=\frac{1}{2}(n^{2}+n-2)$. If we have $\para\neq 0$ and
$\para\in\mathrm{Roots}_{n-1}$, 
then we get $\rk(\gmatrix[n-2])=\frac{1}{2}(n^{2}+n-4)$. For the final case
$\para=0$ we recall that $0\in\mathrm{Roots}_{n-1}$ holds if and only if $n\equiv 0\bmod 2$.
Thus, we get 
$\rk(\gmatrix[n-2])=n-1$ if $n\not\equiv 0\bmod 2$, and $\rk(\gmatrix[n-2])=n-2$ if $n\equiv 0\bmod 2$.

\textit{Case $\momon$, general $\lambda$.}
In general one gets a recursion for 
$\det(\gmatrix[\lambda])$. To explain it, let us
use $n$ together with $\lambda$ in our notation.
The recursion is
\begin{gather}\label{Eq:DAlgebraMotzkinDet}
\scalebox{0.9}{$\det(\gmatrix[\lambda])=
\det(\gmatrix[{n-1,\lambda-1}])
\big(\para^{\#\lcell_{n-1,\lambda}}
\det(\gmatrix[{n-1,\lambda}])\big)
\big(U_{\lambda+1}(\para{-}1)/U_{\lambda}(\para{-}1)\big)^{\#\lcell_{n-1,\lambda+1}}
\det(\gmatrix[{n-1,\lambda+1}])\big)$}
,
\end{gather}
where one omits a factor when the respective Gram matrix is empty.
Here $\#\lcell_{n,\lambda}$ is the number of left cells in $\jcell_{\lambda}$ 
given in \autoref{Eq:DAlgebrasTable}.
This is not 
hard to show using the same strategy as above, see 
for example \cite[Section 5]{BeHa-motzkin} 
(the paper \cite{BeHa-motzkin} has different 
parameters then we do but the arguments given therein work {\muta}).
One can then inductively solve \autoref{Eq:DAlgebraMotzkinDet}
and gets
\begin{gather}\label{Eq:DAlgebraMotzkinDetTwo}
\det(\gmatrix[\lambda])=\para^{x(n,\lambda)}\prod_{t=1}^{\lfloor(n-\lambda)/2\rfloor}
\big(U_{\lambda+t}(\para-1)/U_{t-1}(\para-1)\big)^{\#\lcell_{\lambda+2t}}
,
\end{gather}
with $x(n,\lambda)$ recursively determined by $x(n,\lambda)=x(n-1,\lambda-1)\para^{\#\lcell_{n-1,\lambda}}x(n-1,\lambda)x(n-1,\lambda+1)$ under the same conditions as in \autoref{Eq:DAlgebraMotzkinDet}.
From this one gets different ranks depending on the multiset 
(counting multiplicities as well) 
$\mathrm{Roots}_{n,all}=\{0\}\cup\mathrm{Roots}_{(n+\lambda)/2}\cup\dots\cup\mathrm{Roots}_{\lambda+1}$ for $n+\lambda$ even respectively
$\mathrm{Roots}_{n,all}=\{0\}\cup\mathrm{Roots}_{(n+\lambda-1)/2}\cup\dots\cup\mathrm{Roots}_{\lambda+1}$ for $n+\lambda$ odd.

\textit{Case $\robrmon$, $\lambda=n-2$.}
This case is very similar to $\brmon$ and 
$\momon$ for $\lambda=n-2$ with some small differences. 
That is, we still use the shorthand notation from 
\autoref{Eq:DAlgebrasBrauer} but now also
\begin{gather}\label{Eq:DAlgebrasBrauerTwo}
\begin{tikzpicture}[anchorbase]
\draw[usual] (0.5,0) to (0.5,0.5);
\draw[usual] (1,0) to (1,0.5)node[above]{$\phantom{k}$};
\draw[usual] (1.5,0) to (1.5,0.5);
\draw[usual,dot] (0,0)node[below]{$i$} to (0,0.2);
\draw[usual,dot] (2,0)node[below]{$j$} to (2,0.2);
\end{tikzpicture}
\leftrightsquigarrow
c[i,j]
,\quad
\begin{tikzpicture}[anchorbase]
\draw[usual] (0,0) to (0,0.5);
\draw[usual] (1,0) to (1,0.5);
\draw[usual] (1.5,0)node[below]{$\phantom{k}$} to (1.5,0.5);
\draw[usual,dot] (0.5,0.5)node[above]{$k$} to (0.5,0.3);
\draw[usual,dot] (2,0.5)node[above]{$l$} to (2,0.3);
\end{tikzpicture}
\leftrightsquigarrow
u[k,l]
.
\end{gather}
As one easily checks by drawing the relevant pictures, we then get
\begin{gather*}
\begin{cases}
b[i,j]t[k,l]=\para & \text{if both endpoints match},
\\
b[i,j]t[k,l]=1 & \text{if one endpoint matches},
\\
b[i,j]t[k,l]=0 & \text{else},
\\
c[i,j]u[k,l]=\para^{2} & \text{if both endpoints match},
\\
c[i,j]u[k,l]=0 & \text{else},
\end{cases}
\quad
\begin{cases}
b[i,j]u[k,l]=\para & \text{if both endpoints match},
\\
b[i,j]u[k,l]=0 & \text{else},
\\
c[i,j]t[k,l]=\para & \text{if both endpoints match},
\\
c[i,j]t[k,l]=0 & \text{else},
\end{cases}
\end{gather*}
as entries of the Gram matrix. 
Thus, we have the same block decomposition as in 
\autoref{Eq:DAlgebrasMotzkin} but with
northwest corner corresponding to $\brmon$ and $A_{33}$ being empty.
For example,
\begin{gather*}
n=4\colon
\gmatrix[n-2]=
\begin{psmallmatrix}
\para & 1 & 1 & 1 & 1 & 0 & \para & 0 & 0 & 0 & 0 & 0 \\
1 & \para & 1 & 1 & 0 & 1 & 0 & \para & 0 & 0 & 0 & 0 \\
1 & 1 & \para & 0 & 1 & 1 & 0 & 0 & \para & 0 & 0 & 0 \\
1 & 1 & 0 & \para & 1 & 1 & 0 & 0 & 0 & \para & 0 & 0 \\
1 & 0 & 1 & 1 & \para & 1 & 0 & 0 & 0 & 0 & \para & 0 \\
0 & 1 & 1 & 1 & 1 & \para & 0 & 0 & 0 & 0 & 0 & \para \\
\para & 0 & 0 & 0 & 0 & 0 & \para^{2} & 0 & 0 & 0 & 0 & 0 \\
0 & \para & 0 & 0 & 0 & 0 & 0 & \para^{2} & 0 & 0 & 0 & 0 \\
0 & 0 & \para & 0 & 0 & 0 & 0 & 0 & \para^{2} & 0 & 0 & 0 \\
0 & 0 & 0 & \para & 0 & 0 & 0 & 0 & 0 & \para^{2} & 0 & 0 \\
0 & 0 & 0 & 0 & \para & 0 & 0 & 0 & 0 & 0 & \para^{2} & 0 \\
0 & 0 & 0 & 0 & 0 & \para & 0 & 0 & 0 & 0 & 0 & \para^{2} \\
\end{psmallmatrix}
.
\end{gather*}
From the previous formulas one easily gets that
\begin{gather*}
\det(\gmatrix[n-2])=\para^{\frac{1}{2}n(n-1)}(\para-3)^{\frac{1}{2}n(n-3)}
(\para+n-5)^{n-1}(\para+2n-5)
.
\end{gather*}
Analyzing the rank
using this determinant formula gives the claimed result. 
That is, unless $\para\in\{0,3,5-n,5-2n\}$ we have 
$\rk(\gmatrix[n-2])=n(n-1)$. Otherwise, reading $\{0,3,5-n,5-2n\}$ 
left to right, we get $\rk(\gmatrix[n-2])=\frac{1}{2}n(n-1)$, $\rk(\gmatrix[n-2])=\frac{1}{2}n(n+1)$,
$\rk(\gmatrix[n-2])=(n-1)^{2}$ and $\rk(\gmatrix[n-2])=n(n-1)-1$.
\end{proof}

\begin{Remark}\label{T:DAlgebrasMultiMotzkin}
For a multiparameter version of the Motzkin diagram algebra the formula 
\autoref{Eq:DAlgebraMotzkinDetTwo}
holds after replacing $\para^{x(n,\lambda)}$ with the 
parameter for intervals, keeping the same exponent. 
\end{Remark}

We now compute the dimensions of the 
simple $\xmon$-modules for 
$\xmon\in\{\sym,\psym,\romon,\promon,\tlmon,\ppamon\}$ explicitly, 
and implicitly for $\xmon=\momon$.
We also give some partial results for the remaining ones.

\begin{Theorem}\label{T:DAlgebrasDSimples}
Assume that $\K$ is a field with $\chark=0$ 
for the symmetric monoids, and consider the 
involutive sandwich pair $\big(\xmon(\para),\xmon\big)$.
\begin{enumerate}

\item For $\xmon\in\{\sym,\psym,\romon,\promon,\tlmon,\momon,\ppamon\}$ 
the dimension of the simple $\xmon(\para)$-modules for $\lambda\in\apex$, 
and additionally a simple $\sand[\lambda]$-module $K$ for $\sym$ and 
$\romon$, are
\begin{gather*}
\dim\big(\lmod[{\lambda}]\big)
=\mrk(\gmatrix)
,\quad
\dim\big(\lmod[{\lambda,K}]\big)
=\mrk(\gmatrix)\cdot\dim(K)
.
\end{gather*}
(The ranks of the Gram matrices are computed in \autoref{P:DAlgebrasGram}.)

\item Let $\chark=0$ for $\xmon=\momon$, and let $\K$ be an arbitrary field otherwise.
For $\xmon\in\{\sym,\psym,\romon,\promon,\tlmon,\momon,\ppamon\}$ the algebra $\xmon(\para)$ is semisimple if and only if we are in the following cases:
\begin{enumerate}

\item All cases for $\xmon\in\{\sym,\psym\}$.

\item The parameter $\para$ is nonzero for $\xmon\in\{\romon,\promon\}$.

\item The number $l$ from \autoref{R:DAlgebrasTL} is in $\Z_{\geq n+1}\cup\{\infty\}$ 
for $\xmon\in\{\tlmon,\ppamon[{\lfloor\frac{n}{2}\rfloor}]\}$.

\item We have $\para\not\in\mathrm{Roots}_{n,all}$, for all 
$\lambda\in\Pcal$, for $\xmon=\momon$.

\end{enumerate}

\item For $\lambda\in\apex$ 
we have the following lower bounds, 
where we always want to bound the left side by the right side:
\begin{gather*}
\scalebox{0.9}{$\dim\big(\lmod[{K_{triv},\lambda}]_{\pamon}\big)\geq tl_{2n}^{2\lambda}
,\quad
\dim\big(\lmod[{K_{triv},\lambda}]_{\robrmon}\big)\geq\dim\big(\lmod[\lambda]_{\momon}\big)\geq tl_{n}^{\lambda}
,\quad
\dim\big(\lmod[{K_{triv},\lambda}]_{\brmon}\big)\geq tl_{n}^{\lambda}$}
.
\end{gather*}
Here $K_{triv}$ always denotes the trivial $\sand[\lambda]$-module
and we use subscripts to indicate what kind of simple 
$\xmon(\para)$-modules we consider.

\end{enumerate}
\end{Theorem}

\begin{proof}
\textit{(a).} 
For the planar diagram algebras this is 
just \autoref{P:SandwichCellsGram}, and for $\xmon=\sym$ the statement 
is clear. For $\romon(\para)$ and $\lambda\in\apex$ 
we define a 
pseudo-idempotent $e_{\lambda}$ that have a start-top dot pair for all 
positions $1$ to $n-\lambda$ and is the identity otherwise. For example, for 
$n=5$ and $\lambda=3$ we get
\begin{gather*}
e_{3}=
\begin{tikzpicture}[anchorbase,scale=0.55]
\draw[usual,dot] (0,0) to (0,0.3);
\draw[usual,dot] (0,1) to (0,0.7);
\draw[usual,dot] (0.5,0) to (0.5,0.3);
\draw[usual,dot] (0.5,1) to (0.5,0.7);
\draw[usual,dot] (1,0) to (1,0.3);
\draw[usual,dot] (1,1) to (1,0.7);
\draw[usual] (1.5,0) to (1.5,1);
\draw[usual] (2,0) to (2,1);
\end{tikzpicture}
\,.
\end{gather*}
Hence, we get idempotents $e_{\lambda}^{\prime}=\frac{1}{\para^{n-\lambda}}e_{\lambda}$
in $\jcell_{\lambda}$.
Note also that $\dmod[\lambda]\cong\romon(\para)e_{\lambda}$ 
as $\romon(\para)$-modules.
Let $e_{\sigma(\lambda)}$ denote the projection 
to the component of $\sand^{\oplus k}$ 
determined by the top through strand points of $\sigma e_{\lambda}$, 
where we write $k=\binom{n}{\lambda}$.
The $\K$-linear map $f_{\lambda}\colon\dmod[\lambda]\to\sand^{\oplus k}$ 
given by $f_{\lambda}(\sigma e_{\lambda}^{\prime})=\sigma e_{\sigma(\lambda)}$
is an isomorphism of $\sand$-modules.
Thus, $\dmod[{\lambda,K}]\cong K^{\oplus k}$ 
for all simple $\sand$-modules $K$. Using now the permutation action 
of $\sand[n]$, which connects the various direct summands, one sees that 
$\dmod[{\lambda,K}]$ is a simple $\romon(\delta)$-module.

\textit{(b).} This follows from (a) and \autoref{P:SandwichCellsSemisimple}.

\textit{(c).} Note that 
we have algebra embeddings $\ppamon(\para)\hookrightarrow\pamon(\para)$, 
$\tlmon(\para)\hookrightarrow\momon(\para)\hookrightarrow\robrmon(\para)$
and $\tlmon(\para)\hookrightarrow\brmon(\para)$
given by sending diagrams to themselves interpreted in the bigger monoids.
Moreover, by \autoref{P:DAlgebrasDiagram}, 
these embeddings behave well with the cell order 
in the sense that diagrams in $\jcell_{k}^{source}$ (for the source algebra) are send to diagrams in $\jcell_{k}^{target}$ (for the target algebra). This enables us to compare the sandwich matrix of $\jcell_{k}^{source}$ with 
the one for $\jcell_{k}^{target}$.
Recalling the picture that defines the sandwich matrices, see \autoref{Eq:SandwichMatrixPicture}, since the embedding sends diagrams to themselves, the computations for each entry remain the same. 
Thus, using an appropriate order, the sandwich matrices for the small monoids are submatrices of the sandwich matrices for the bigger monoids and 
\autoref{P:DAlgebrasGram} completes the proof.
\end{proof}

\begin{Remark}\label{R:DAlgebraQuantumTwo}
Similarly as in \autoref{R:DAlgebraQuantum}, 
the quantum versions of the diagram algebras 
$\xmon(\para)$ 
({\eg} \emph{BMW algebras}) can be studied {\ver}, and the Hecke algebras 
of type $A$ appear again as the sandwiched algebras.
\end{Remark}

\begin{Remark}\label{R:DAlgebrasHistory}
The study of diagram algebras has a long history, which we cannot cover 
here in any satisfying way. Nevertheless, let us give a few references 
for some of the results above that have appeared in the literature, but not under the umbrella of sandwich cellularity, which is new in this paper.

The representation theories of $\sym$ and $\psym$ is, of course, well-studied, and we do not comment on these any further.

Statements such as \autoref{P:DAlgebrasDiagram} 
and \autoref{P:DAlgebrasGram} appear throughout 
the semigroup literature, see {\eg} \cite{HaRa-partition-algebras}, 
\cite{Do++-idempotents-diagram-monoids}, \cite{HaJa-representations-diagram-algebras}, sometimes as 
disguised lemmas such as \cite[Remark 9.21 and Lemma 9.19]{EaMiRuTo-congruence-lattices},
and many more. They have very similar
appearances for diagram algebras with the catch that the parameter 
$\para$ is not necessarily $1$ as for the monoids.

Computations of dimensions of simple modules 
of diagram algebras such as in \autoref{T:DAlgebrasDSimples}
have also appeared in many works (although some of the results above appear 
to be new). For example, for $\romon$ see \cite{Mu-sym-inv-semigroup} and
\cite{So-rook-monoid}, and for $\promon$ see \cite{FlHaHe-planar-rook}.
The formulas from \autoref{R:DAlgebrasTL} for $\tlmon$ go back to \cite{RiSaAu-temperley-lieb}, \cite{An-simple-tl} and \cite{Sp-modular-tl}, 
and the ones for $\ppamon$ follows then from the isomorphism 
to $\tlmon[2n]$ from \cite{HaRa-partition-algebras}. We do not know 
any reference for the dimensions of simple modules for any of the remaining 
symmetric monoids, but the Gram determinants 
for $\momon$ appear in \cite{BeHa-motzkin}. The bounds however have been studied in \cite{KhSiTu-monoidal-cryptography}. Let us also note that a lot of these 
diagram algebras and similar diagram algebra can be studied using 
Schur--Weyl--Brauer duality, see {\eg} \cite[Section 3]{AnStTu-semisimple-tilting} for a collection of these dualities.

The quantum versions, which also fit into the theory, {\cf} \autoref{R:DAlgebraQuantumTwo}, have also been studied extensively, 
see {\eg} \cite{Xi-bmw-cellular} or \cite{En-brauer-bmw-cellular} for the BMW case. Note also 
that \autoref{P:DAlgebrasDiagram}.(d) shows that all the 
diagram algebras in this section are cellular. This is known, 
see {\eg} \cite{GrLe-cellular}, \cite{Xi-partition-cellular}
and \cite{En-brauer-bmw-cellular}. But again our point is that we get all of this by using the 
theory of sandwich cellular algebras.
\end{Remark}

\newcommand{\etalchar}[1]{$^{#1}$}

\end{document}